\documentclass[3p,sort&compress]{elsarticle}
\journal{Journal of Computational Physics (Elsevier)}

\usepackage{framed,multirow}
\usepackage{latexsym}
\usepackage{amsmath}
\usepackage{amssymb}
\usepackage{amsthm}
\usepackage{mathtools}
\usepackage{pifont}
\usepackage[shortlabels]{enumitem}
\usepackage{booktabs}
\usepackage{subcaption}
\usepackage[dvipsnames]{xcolor}
\usepackage{yfonts}
\usepackage[scr=boondoxo,scrscaled=1.05]{mathalfa}
\usepackage{appendix}
\usepackage{enumitem}
\usepackage{url}
\usepackage{ifthen}
\usepackage{tensor}

\theoremstyle{plain}

\theoremstyle{definition}
\newtheorem{definition}{Definition}

\theoremstyle{remark}
\newtheorem{remark}{Remark}

\begin{document}

\begin{frontmatter}
\title{Multi-element flow-driven spectral chaos (ME-FSC) method for uncertainty quantification of dynamical systems}

\author[1,2]{Hugo Esquivel}
\ead{hesquive@purdue.edu, hesquive@cuc.edu.co}

\author[1]{Arun Prakash}
\ead{aprakas@purdue.edu}

\author[3]{Guang Lin}
\ead{guanglin@purdue.edu}

\address[1]{Lyles School of Civil Engineering, Purdue University, West Lafayette, Indiana, USA}
\address[2]{Department of Civil and Environmental Engineering, Universidad de la Costa, Barranquilla, Colombia}
\address[3]{Department of Mathematics, School of Mechanical Engineering, Purdue University, West Lafayette, Indiana, USA}

\begin{abstract}
The flow-driven spectral chaos (FSC) is a recently developed method for tracking and quantifying uncertainties in the long-time response of stochastic dynamical systems using the spectral approach.
The method uses a novel concept called \emph{enriched stochastic flow maps} as a means to construct an evolving finite-dimensional random function space that is both accurate and computationally efficient in time.
In this paper, we present a multi-element version of the FSC method (the ME-FSC method for short) to tackle (mainly) those dynamical systems that are inherently discontinuous over the probability space.
In ME-FSC, the random domain is partitioned into several elements, and then the problem is solved separately on each random element using the FSC method.
Subsequently, results are aggregated to compute the probability moments of interest using the law of total probability.
To demonstrate the effectiveness of the ME-FSC method in dealing with discontinuities and long-time integration of stochastic dynamical systems, four representative numerical examples are presented in this paper, including the Van-der-Pol oscillator problem and the Kraichnan-Orszag three-mode problem.
Results show that the ME-FSC method is capable of solving problems that have strong nonlinear dependencies over the probability space, both reliably and at low computational cost.
\end{abstract}

\begin{keyword}
stochastic discontinuities; stochastic dynamical systems; uncertainty quantification; long-time integration; stochastic flow map; multi-element flow-driven spectral chaos (ME-FSC).
\end{keyword}
\end{frontmatter}

\section*{Highlights}
\begin{itemize}\setlength\itemsep{0em}
\item A multi-element flow-driven spectral chaos method.
\item ME-FSC handles discontinuities over the probability space adequately.
\item ME-FSC quantifies uncertainties in the long-time response reliably.
\item ME-FSC is highly accurate and computationally efficient.
\end{itemize}

\section{Introduction}

The spectral approach has gained increasing popularity in the past few decades as a powerful tool for solving stochastic problems at low-computational cost.
This assertion is especially true for problems where the dimensionality of the probability space is relatively low.
For problems where the dimensionality of the probability space is relatively high, this assertion may not always hold, especially in cases where the curse of dimensionality at the random-function-space (RFS) level cannot be alleviated noticeably.
This is, for instance, the case of gPC-based methods (e.g.~\cite{xiu2002wiener,gerritsma2010time,heuveline2014hybrid,luchtenburg2014long,ozen2016dynamical}) where the curse-of-dimensionality issue arises naturally whenever a suitable random function space (aka RFS in this paper) is defined over the probability space to represent approximately the stochastic part of the system's response.
In other words, when the random input consists of mutually independent random variables, the stochastic part of the solution space is customarily spanned via a tensor product of vector spaces with each of these vector spaces representing a single space spanned by a finite set of univariate orthogonal polynomials.
As a consequence, the higher the dimensionality of the probability space, the higher the dimensionality of the resulting RFS.
It is worth mentioning that an RFS constructed this way can only ensure exponential convergence to the solution in the early times of the simulation, unless the RFS is kept updated frequently during the simulation.

A brief historical review is given next.
Wan and Karniadakis \cite{wan2005adaptive} developed the multi-element generalized polynomial chaos (ME-gPC) method to handle long-time integration and discontinuities satisfactorily in problems with random input data.
In this approach, the random domain is decomposed into smaller subdomains (aka random elements) whenever the relative error in variance exceeds a predefined threshold value during the simulation.
For computational efficiency, each subdomain uses a relatively low-degree polynomial chaos from the Askey scheme to ensure exponential convergence to the solution locally \cite{xiu2002wiener}.
The effectiveness of ME-gPC to deal with long-time integration was studied in detail in \cite{wan2006long} (see also \cite{wan2009error}), and subsequently the method was further generalized for inputs with arbitrary probability density functions \cite{wan2006beyond}.
To benefit from the computational ease of sampling-based methods, Foo et al.~\cite{foo2008multi} developed the multi-element probabilistic collocation (ME-PCM) method.
This approach uses the decomposition of the random domain and a tensor product grid in each subdomain to deal with long-time integration and stochastic discontinuities more efficiently.
In \cite{foo2010multi}, the rate of convergence was improved by introducing an analysis of variance (ANOVA) into the scheme.
Furthermore, Zheng et al.~\cite{zheng2015adaptive} developed an adaptive multi-element polynomial chaos for discrete probability density functions, and Kawai and Oyama \cite{kawai2020multi} proposed a multi-element stochastic Galerkin method based on edge detection to circumvent the Gibbs phenomenon in the presence of stochastic discontinuities.

Asokan and Zabaras \cite{asokan2006stochastic} developed a variational multiscale framework to solve elliptic problems with heterogeneous random media.
The framework combines ideas from the variational multiscale method, the multiscale finite element method, and the generalized polynomial chaos method to derive a variationally consistent upscaling technique to model diffusion in problems involving heterogeneous random media.
Witteveen et al.~\cite{witteveen2009adaptive} formulated an adaptive stochastic finite elements approach based on the Newton-Cotes quadrature rule.
Here, the stochastic part of the solution is approximated with piecewise polynomial chaos by subdividing the random domain into simplex elements.
Nouy and Cl\'ement \cite{nouy2010extended} combined the extended finite element method with spectral stochastic methods to simulate heterogeneous materials with random material interfaces.
Here, to enhance the convergence rate of the scheme, an extension of the partition of unity method was proposed in the spectral stochastic framework.

Ganapathysubramanian and Zabaras \cite{ganapathysubramanian2007sparse} presented an application of sparse grid collocation schemes based on the Smolyak algorithm to solve high-dimensional stochastic convection problems more effectively.
Ma and Zabaras \cite{ma2009adaptive} proposed an adaptive sparse grid collocation strategy using piecewise multi-linear hierarchical basis functions.
To automatically detect stochastic discontinuities, and thus, adaptively refine the collocation points in the problematic subdomains, hierarchical surpluses were used as error indicators.
Bhaduri and Graham-Brady \cite{bhaduri2018efficient} improved this technique by avoiding unnecessary function evaluations in smoother regions of the random domain by using successive derivative estimations along all random dimensions (see also \cite{bhaduri2018stochastic}).
Agarwal and Aluru \cite{agarwal2009domain} proposed an adaptive sparse grid collocation strategy to investigate the performance of micro-electromechanical systems under uncertainties.

Marzouk et al.~\cite{marzouk2007stochastic} presented a reformulation of the Bayesian approach to inverse problems to accelerate the Bayesian inference via the gPC method.
Mohan et al.~\cite{mohan2008multi} developed the multi-element stochastic reduced basis methods (ME-SRBMs) for solving linear stochastic partial differential equations. 
In their approach, the random domain is decomposed into several subdomains and the solution is approximated on each subdomain using a set of basis vectors spanning a preconditioned stochastic Krylov subspace.
Sarrouy et al.~\cite{sarrouy2013piecewise} proposed a piecewise polynomial chaos expansion to perform a stability analysis of a linear brake system with uncertainties using the complex eigenvalue analysis method and decomposition of the random domain to realize a low degree piecewise polynomial approximation. 
Jakeman et al.~\cite{jakeman2013minimal} developed the minimal multi-element stochastic method to quantify uncertainties of discontinuous functions more efficiently than standard locally-adaptive sparse grid methods.
This approach consists of two major steps.
First, a discontinuity detector is introduced into the scheme to partition the random domain into elements of high regularity.
Then, an adaptive technique based on the least orthogonal interpolant is employed to construct a gPC approximation on each element.

Ma\^itre et al.~\cite{le2004uncertainty} introduced the Wiener-Haar expansion in the context of uncertainty propagation to investigate situations in which the solution can vary rapidly over the probability space.
By combining concepts from the generalized polynomial chaos methods \cite{xiu2002wiener,wan2005adaptive} and stochastic Galerkin methods, this wavelet-based technique was demonstrated to be suitable for resolving stochastic systems with multiple stochastic discontinuities (see also \cite{le2004multi}).
To guarantee optimal convergence with respect to the full stochastic and spatial discretization, adaptive wavelets were used in \cite{gittelson2014adaptive} to construct simultaneously the stochastic and spatial bases and omit the intermediate semidiscrete approximation stage.
Gittelson \cite{gittelson2013adaptive} derived an adaptive solver for elliptic boundary value problems with random coefficients.
Though this technique uses ideas from adaptive wavelet methods, orthonormal polynomials were used instead.

Cho et al.~\cite{cho2013adaptive} developed an adaptive discontinuous Galerkin method for response-excitation PDF equations using a nonconforming adaptive discontinuous Galerkin method for the response space and a probabilistic collocation method for the excitation space.
Heuveline and Schick \cite{heuveline2014hybrid} proposed a hybrid generalized polynomial chaos method to tackle problems featuring strong nonlinear dependencies on the stochastic inputs (with or without stochastic discontinuities).
This was achieved by combining ideas from TD-gPC \cite{gerritsma2010time} and ME-gPC \cite{wan2005adaptive} methods.
Chen et al.~\cite{chen2015local} developed a localized polynomial chaos expansion for high-dimensional stochastic problems to circumvent the difficulty of having a high-dimensional probability space in the mathematical model of the system.
Li and Stinis \cite{li2016unified} presented a unified framework for mesh refinement in random and physical space.
The benefit of the framework is that it does not require any explicit knowledge of a reduced model to perform the mesh refinement in both these spaces.
Finally, B-splines basis functions have also been used for uncertainty propagation to overcome the limitations of classical approaches in the presence of discontinuities (see e.g.~\cite{abdedou2019non,eckert2020polynomial}).

Lastly, the flow-driven spectral chaos (FSC) \cite{esquivel2020flow,esquivel2021flow} is a new numerical method developed by the authors of the present paper to tackle mainly the long-time integration issue found in the gPC method using the spectral approach.
The method uses the concept of \emph{enriched} stochastic flow maps to track the evolution of a finite-dimensional RFS efficiently in time.
In FSC, the dimensionality of the random phase space is deliberately increased to allow both the system's state and its first few time derivatives to be pushed forward over time.
Then, the enriched state of the system is used as a germ to construct a suitable RFS for use within the current time step of the simulation.
It is worth mentioning that for dynamical systems with an $n$-tuple state and driven by a stochastic flow map of order $M$, the cardinality of the random basis is bounded from above by $n+M+1$.
This boundedness from above is what makes the FSC method be curse-of-dimensionality free at the RFS level, even when the probability space is high-dimensional.

In this paper, we present a multi-element version of the FSC method to deal with stochastic discontinuities and long-time integration of stochastic dynamical systems more efficiently.
This new technique is called the \emph{multi-element flow-driven spectral chaos} (ME-FSC) method.
In ME-FSC, the random domain is partitioned into several elements, and then the problem is solved separately on each random element using the FSC method.
Then, the results are aggregated to compute the mean and variance of the response with the law of total probability.
This approach is similar to the multi-element gPC (ME-gPC) method \cite{wan2005adaptive}, with the only difference being that the gPC method is not employed on each random element to fully annihilate the curse of dimensionality at the RFS level.
The benefit of using ME-FSC is threefold.
First, the simulation can be run simultaneously on machines with multiple CPU cores (or if needed on separate machines) to reduce the excessive computational burden associated with the simulation.
Second, suppose the random input is discontinuous over the probability space.
In that case, the random domain can be partitioned into several elements so that the discontinuity only appears on regions of measure zero.
Third, if an adaptive criterion is introduced within the ME-FSC scheme (such that the number of elements gets smaller on-the-fly whenever a threshold value is exceeded), the errors can be kept to a minimum during the simulation.

This paper has been structured in the following way.
The setting and notation used in the manuscript are formally introduced in Section \ref{sec3SetNot}, and then the problem statement we are interested in is presented in Section \ref{sec3ProSta}.
A quick overview of the concept associated with `enriched stochastic flow maps' is provided in Section \ref{sec3BothStoFloMap} to allow the construction of an optimal, low-dimensional RFS during the simulation.
Subsequently, in Section \ref{sec3FSCmet} an overview of the FSC method is imparted to highlight the key idea behind the FSC method.
In Section \ref{sec3MEFSCmet} the ME-FSC method is presented in detail along with an outline of the numerical scheme used in this work.
Finally, in Section \ref{sec3NumRes}, four numerical examples are explored to test the performance of the ME-FSC method (against well-established UQ techniques) and document the findings.

\section{Setting and notation}\label{sec3SetNot}

\subsection{Basic spaces needed in the ME-FSC method}\label{sec3BasSpaMEFSC}

Definitions 1 to 5 presented in two previous works by the authors \cite{esquivel2020flow,esquivel2021flow} are also considered herein.
The \emph{temporal space} is a topological space defined by $\mathfrak{T}=[0,T]$, where $T\in\mathbb{R}^+$ symbolizes the duration of the simulation.
The \emph{probability space} is defined by $(\Omega,\boldsymbol{\Omega},\lambda)$, where $\Omega$ is the sample space, $\boldsymbol{\Omega}\subset 2^\Omega$ is the $\sigma$-algebra on $\Omega$, and $\lambda:\boldsymbol{\Omega}\to[0,1]$ is the probability measure on $\boldsymbol{\Omega}$. 
Because this probability space can be abstract, a $d$-valued random variable $\xi=(\xi^1,\ldots,\xi^d):\Omega\to\mathbb{R}^d$ given by $\xi=\xi(\omega)$ is defined for computational purposes.
Thus, the \emph{random space} is a measure space defined by $(\Xi,\boldsymbol{\Xi},\mu)$, where $\Xi=\xi(\Omega)\subset\mathbb{R}^d$ is a set representing the random domain of the system, $\boldsymbol{\Xi}=\mathcal{B}_{\mathbb{R}^d}\cap\Xi$ is the $\sigma$-algebra on $\Xi$, $\mu:\boldsymbol{\Xi}\to[0,1]$ is the probability measure on $\boldsymbol{\Xi}$ given by $\mu=\xi_*(\lambda):=\lambda\circ\xi^{-1}$ (i.e.~the push-forward of $\lambda$ by $\xi$).
Naturally, $d$ symbolizes the dimensionality of the random space.
The \emph{temporal function space} is a continuous $n$-differentiable function space defined by $\mathscr{T}(n)=C^n(\mathfrak{T};\mathbb{R})$.
In other words, this is the space of all functions $f:\mathfrak{T}\to\mathbb{R}$ that have continuous first $n$ derivatives on $\mathfrak{T}$.
Moreover, the \emph{random function space}, defined by $\mathscr{Z}=(L^2(\Xi,\boldsymbol{\Xi},\mu;\mathbb{R}),\langle\,\cdot\,,\cdot\,\rangle)$, is a Lebesgue square-integrable space equipped with its standard inner product.
This is the space of all (equivalence classes of) measurable functions $f:\Xi\to\mathbb{R}$ that are square $\mu$-integrable on $\Xi$.
The standard inner product $\langle\,\cdot\,,\cdot\,\rangle: L^2(\Xi,\boldsymbol{\Xi},\mu;\mathbb{R})\times L^2(\Xi,\boldsymbol{\Xi},\mu;\mathbb{R})\to\mathbb{R}$ is given by $\langle f,g\rangle=\int fg\,\mathrm{d}\mu$.
For computational purposes, this space is spanned by a complete orthogonal basis, $\{\Psi_j:\Xi\to\mathbb{R}\}_{j=0}^\infty$ such that $\Psi_0(\xi)=1$ for all $\xi\in\Xi$.
Finally, the \emph{solution space} and the \emph{root space} are defined by $\mathscr{U}=\mathscr{T}(n)\otimes\mathscr{Z}$ and $\mathscr{V}=\mathscr{T}(0)\otimes\mathscr{Z}$, respectively, making $\mathscr{U}$ and $\mathscr{V}$ be two tensor products of vector spaces.

Four remarks are in order.
\emph{First}, it is well-known that $\mathscr{Z}$ forms a Hilbert space because it is complete under the metric induced by $\langle\,\cdot\,,\cdot\,\rangle$.
\emph{Second}, if $f\in\mathscr{Z}$, then it can be represented by the Fourier series $f=\sum_{j=0}^\infty f^j\Psi_j$, where $f^j$ denotes the $j$-th coefficient of the series with the superscript not symbolizing an exponentiation.
\emph{Third}, the orthogonality property of the basis in $\mathscr{Z}$ means that $\langle\Psi_i,\Psi_j\rangle:=\int\Psi_i\Psi_j\,\mathrm{d}\mu=\langle\Psi_i,\Psi_i\rangle\,\delta_{ij}$, where $\delta_{ij}$ is the Kronecker delta.
\emph{Fourth}, the dual space of $\mathscr{Z}$, which is denoted by $\mathscr{Z}'$ herein, is the space spanned by the continuous linear functionals $\{\Psi^i:\mathscr{Z}\to\mathbb{R}\}_{i=0}^\infty$ defined by $\Psi^i[f]=\langle\Psi_i,f\rangle/\langle\Psi_i,\Psi_i\rangle\equiv f^i$.

\begin{figure}
\centering
\includegraphics[]{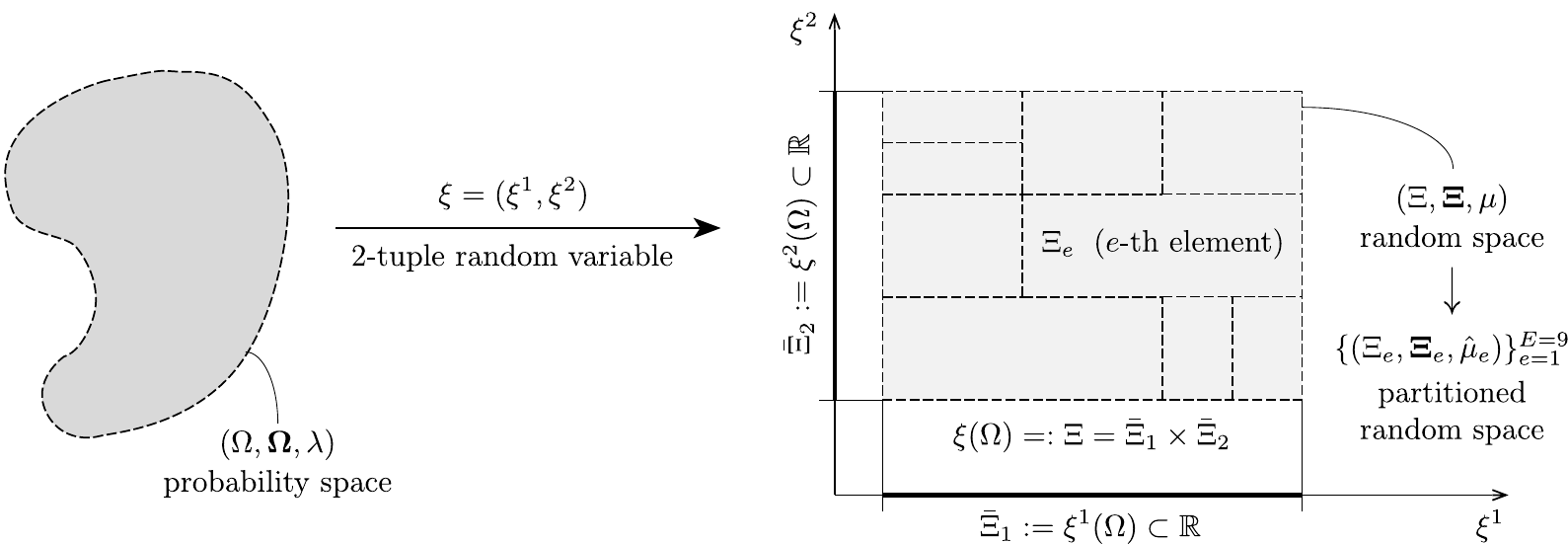}
\caption{Relationship between probability space and random space for $d=2$ (along with a case partition of the random space)}
\label{fig3ProbabilityRandomSpace}
\end{figure}

Fig.~\ref{fig3ProbabilityRandomSpace} depicts the relationship between probability space and random space for $d=2$, and a case scenario of a random space with 9 elements.
As shown, the components of $\xi=(\xi^1,\ldots,\xi^d)$ will always be assumed to be mutually independent, which means that the random domain $\Xi$ is a rectangular hyper-solid of $d$ dimensions constructed by performing a $d$-fold Cartesian product of intervals $\bar{\Xi}_i:=\xi^i(\Omega)$.

\subsection{Additional spaces needed in the ME-FSC method}

The following two definitions are needed in the development of the ME-FSC method.

\begin{definition}[Partitioned random space]\label{sec3SetNotDef6}
Let $\{\Xi_e\}_{e=1}^E$ be a partition of the random domain, where $\Xi_e\neq\emptyset$ represents the $e$-th element of the partition, and $E\in\mathbb{N}_2$ is the number of random elements employed in the partition.
Let $\boldsymbol{\Xi}_e=\boldsymbol{\Xi}\cap\Xi_e$ be the $\sigma$-algebra on $\Xi_e$, and $\mu_e=\mu|_{\boldsymbol{\Xi}_e}$ be the restriction of $\mu$ to $\boldsymbol{\Xi}_e$.
Moreover, let $\hat{\mu}_e:\boldsymbol{\Xi}_e\to[0,1]$ be the probability measure on $\boldsymbol{\Xi}_e$ given by
\begin{equation*}
\hat{\mu}_e=\frac{\mu_e}{\mu(\Xi_e)},\quad\text{or equivalently,}\quad
\mathrm{d}\hat{\mu}_e=\frac{\mathrm{d}\mu_e}{\mu(\Xi_e)}.
\end{equation*}
As a result, in this work the \emph{partitioned random space} is defined as a finite sequence of disjoint random spaces $\{(\Xi_e,\boldsymbol{\Xi}_e,\hat{\mu}_e)\}_{e=1}^E$.
A prototype depiction of such a partition is displayed in Fig.~\ref{fig3ProbabilityRandomSpace}.
\end{definition}

\begin{definition}[Partitioned random function space]
If $\{(\Xi_e,\boldsymbol{\Xi}_e,\hat{\mu}_e)\}_{e=1}^E$ is a partition of the random space, then let $\{\mathscr{Z}_e\}_{e=1}^E$ be its associated \emph{partitioned random function space}, where $\mathscr{Z}_e$ represents the $e$-th RFS for the random element $\Xi_e$.
That is, $\mathscr{Z}_e=(L^2(\Xi_e,\boldsymbol{\Xi}_e,\hat{\mu}_e;\mathbb{R});\langle\,\cdot\,,\cdot\,\rangle_e)$, where
\begin{equation*}
\langle\,\cdot\,,\cdot\,\rangle_e: L^2(\Xi_e,\boldsymbol{\Xi}_e,\hat{\mu}_e;\mathbb{R})\times L^2(\Xi_e,\boldsymbol{\Xi}_e,\hat{\mu}_e;\mathbb{R})\to\mathbb{R}\qquad:\Leftrightarrow\qquad\langle f,g\rangle_e=\int fg\,\mathrm{d}\hat{\mu}_e.
\end{equation*}
This is the space of all (equivalence classes of) measurable functions $f:\Xi_e\to\mathbb{R}$ that are square $\hat{\mu}_e$-integrable on $\Xi_e$.
As in $\mathscr{Z}$, the complete orthogonal basis in $\mathscr{Z}_e$, $\{\Psi_{j.e}:(\Xi_e,\boldsymbol{\Xi}_e)\to(\mathbb{R},\mathcal{B}_\mathbb{R})\}_{j=0}^\infty$, is defined such that $\Psi_{0.e}(\xi)=1$ for all $\xi\in\Xi_e$.
\end{definition}

Note that the four remarks made in Section \ref{sec3BasSpaMEFSC} are also applicable for each of the $\mathscr{Z}_e$'s mentioned above, provided that the following changes are made:
\begin{equation*}
\mathscr{Z}\mapsto\mathscr{Z}_e,\quad\mu\mapsto\hat{\mu}_e,\quad \langle\,\cdot\,,\cdot\,\rangle\mapsto\langle\,\cdot\,,\cdot\,\rangle_e,\quad f\mapsto f_{.e},\quad f^k\mapsto f{^k}_{.e},\quad\Psi_k\mapsto\Psi_{k.e}\quad\text{and}\quad\Psi^k\mapsto\Psi{^k}_{.e}
\end{equation*}
with $k$ symbolizing $i$ or $j$.

\section{Problem statement}\label{sec3ProSta}

For convenience, the same stochastic problem described in \cite{esquivel2020flow} is considered herein.
Namely:

Find the real-valued stochastic process $u:\mathfrak{T}\times\Xi\to\mathbb{R}$ in $\mathscr{U}$, such that ($\mu$-a.e.):
\begin{subequations}\label{eq3ProSta1000}
\begin{align}
\mathcal{L}[u]=f&\qquad\text{on $\mathfrak{T}\times\Xi$}\label{eq3ProSta1000a}\\
\big\{\mathcal{B}_k[u](0,\cdot\,)=b_k\big\}_{k=1}^n&\qquad\text{on $\{0\}\times\Xi$},\label{eq3ProSta1000b}
\end{align}
\end{subequations}
where the partial differential operators $\mathcal{L}:\mathscr{U}\to\mathscr{V}$ and $\mathcal{B}_k:\mathscr{U}\to\mathscr{Z}$ are of order $(n,0)$ and $(n-1,0)$, respectively.
The meaning of $\mathcal{B}_k[u](0,\cdot\,)$ is that once the operator $\mathcal{B}_k$ is applied to $u$ the resulting function is evaluated at $t=0$.
The functions $f:\mathfrak{T}\times\Xi\to\mathbb{R}$ and $b_k:\Xi\to\mathbb{R}$ are elements of $\mathscr{V}$ and $\mathscr{Z}$, respectively, and they are given by $f=f(t,\xi)$ and $b_k=b_k(\xi)$.

However, for illustration purposes, we very often consider the much simpler case\footnote{We highlight the fact that this work has been generalized to deal with problems as complicated as \eqref{eq3ProSta1000} in Section \ref{sec3MEFSCmet}.} where the mathematical model of the system can be represented ($\mu$-a.e.) by an undamped oscillator under the action of an external force:
\addtocounter{equation}{-1}
\begin{subequations}
\makeatletter
\def\@currentlabel{\ref{eq3ProSta1000}*}
\makeatother
\label{eq3ProSta1010}
\begin{align}
\mathcal{L}[u]=f\qquad &:\Leftrightarrow\qquad m\ddot{u}+ku=p\tag{\ref{eq3ProSta1000a}*}\label{eq3ProSta1010a}\\
\big\{\mathcal{B}_k[u](0,\cdot\,)=b_k\big\}_{k=1}^{n=2}\qquad &:\Leftrightarrow\qquad\big\{u(0,\cdot\,)=\mathscr{u},\,\dot{u}(0,\cdot\,)=\mathscr{v}\big\},\tag{\ref{eq3ProSta1000b}*}\label{eq3ProSta1010b}
\end{align}
\end{subequations}
where $m:\Xi\to\mathbb{R}^+$ is the mass of the system, $k:\Xi\to\mathbb{R}^+$ is the stiffness of the system, $p:\mathfrak{T}\times\Xi\to\mathbb{R}$ is the external force applied to the system, and $\mathscr{u},\mathscr{v}:\Xi\to\mathbb{R}$ are the prescribed initial conditions of the system.
Observe that $\dot{u}:=\partial_t u$ and $\ddot{u}:=\partial_t^2 u$ are the first and second partial derivatives of $u$ with respect to time.
It is assumed that $m,k,\mathscr{u},\mathscr{v}\in\mathscr{Z}$, $p\in\mathscr{V}$ and of course that $u\in\mathscr{U}$.

\section{Ordinary and enriched stochastic flow maps}\label{sec3BothStoFloMap}

The attention is now turned to the two flow maps considered in this manuscript: the ordinary stochastic flow map, $\varphi$, and the enriched stochastic flow map, $\hat{\varphi}$.\footnote{Please refer to \cite{esquivel2020flow} for a generalization of the two stochastic flow maps considered in this section.}

\subsection{Ordinary stochastic flow map}\label{sec3StoFloMap}

First of all, the stochastic system given by \eqref{eq3ProSta1010} is assumed to be sufficiently regular over $\mathfrak{T}\times\Xi$. Thus,
\begin{subequations}\label{eq3StoFloMap1000}
\begin{align}
\partial_t^2 u(t,\xi):=\mathscr{f}(t,\xi,s(t,\xi))=\bar{p}(t,\xi)-\bar{k}(\xi)\,u(t,\xi) & \qquad\text{on $\mathfrak{T}\times\Xi$}\label{eq3StoFloMap1000a}\\
\big\{u(0,\xi)=\mathscr{u}(\xi),\,\dot{u}(0,\xi)=\mathscr{v}(\xi)\big\}&\qquad\text{on $\{0\}\times\Xi$},\label{eq3StoFloMap1000b}
\end{align}
\end{subequations}
where $s=(u,\dot{u})\in\prod_{j=1}^2\mathscr{T}(3-j)\otimes\mathscr{Z}$ is the configuration state of the system over $\mathfrak{T}\times\Xi$, $\mathscr{f}:\mathfrak{T}\times\Xi\times\mathbb{R}^2\to\mathbb{R}$ is a noisy, non-autonomous function (which can also be regarded as a function in $\mathscr{V}$) defining the response $\ddot{u}=\partial_t^2 u$, $\bar{p}=p/m:\mathfrak{T}\times\Xi\to\mathbb{R}$ is the external force applied to the system per unit mass, and $\bar{k}=k/m:\Xi\to\mathbb{R}^+$ is the stiffness of the system per unit mass.
In the context of \eqref{eq3ProSta1010}, by `sufficiently regular over $\mathfrak{T}\times\Xi$' we mean that $\bar{p}\in\mathscr{V}$ and $\bar{k}\in\mathscr{Z}$ (and thus, $\mathscr{f}\in\mathscr{V}$).

If the solution is further assumed to be analytic on $\mathfrak{T}$ for all $\xi\in\Xi$, then $u$ can be expanded around $t=t_i$ by the Taylor series:
\begin{equation*}
u(t_i+h,\xi)=\sum_{j=0}^\infty \frac{h^j}{j!}\partial_t^j u(t_i,\xi)=\sum_{j=0}^M \frac{h^j}{j!}\partial_t^j u(t_i,\xi)+O(h^{M+1})(\xi),
\end{equation*}
where $t_i\in\mathfrak{T}$ is the time instant of the simulation, and $h:=t-t_i$ is the associated time-step size (once $t$ is fixed).

Hence, a (local) stochastic flow map of order $M$, $\varphi(M):\mathbb{R}\times\mathscr{Z}^2\to\mathscr{Z}^2$, can be stipulated for \eqref{eq3StoFloMap1000} to be:
\begin{equation}\label{eq3StoFloMap1080}
\varphi(M)(h,s(t_i,\cdot\,)):=\big(u(t_i+h,\cdot\,),\dot{u}(t_i+h,\cdot\,)\big)-O(h^{M+1}),
\end{equation}
where $\varphi^1(M),\varphi^2(M):\mathbb{R}\times\mathscr{Z}^2\to\mathscr{Z}$ are given by:
\begin{gather*}
\varphi^1(M)(h,s(t_i,\cdot\,)):=u(t_i,\cdot\,)+h\,\dot{u}(t_i,\cdot\,)+\tfrac{1}{2}h^2\,\partial_t^2u(t_i,\cdot\,)+\cdots+\tfrac{1}{M!}h^M\,\partial_t^Mu(t_i,\cdot\,)\\
\varphi^2(M)(h,s(t_i,\cdot\,)):=\dot{u}(t_i,\cdot\,)+h\,\partial_t^2u(t_i,\cdot\,)+\tfrac{1}{2}h^2\,\partial_t^3u(t_i,\cdot\,)+\cdots+\tfrac{1}{M!}h^M\,\partial_t^{M+1}u(t_i,\cdot\,).
\end{gather*}

It is worth mentioning that the second time derivative $\partial_t^2u(t_i,\cdot\,)$ is computed using \eqref{eq3StoFloMap1000a}, and that the next time derivatives $\{\partial_t^{j+2}u(t_i,\cdot\,)\}_{j=1}^{M-1}$ are computed using the recursive expression:
\begin{equation}\label{eq3StoFloMap1090}
\partial_t^{j+2}u(t,\xi):=\mathrm{D}_t^j\mathscr{f}(t,\xi,s(t,\xi))=\partial_t^j\bar{p}(t,\xi)-\bar{k}(\xi)\,\partial_t^ju(t,\xi)\qquad\forall j\in\{1,2,\ldots,M-1\}.
\end{equation}

In simple terms, the goal of $\varphi(M)$ is to push the system's state $s(t_i,\cdot\,)=(u(t_i,\cdot\,),\dot{u}(t_i,\cdot\,))$ through time using the local information at $t=t_i$.
We recall that this push is \emph{forward} if $h>0$, it is \emph{backward} if $h<0$, and it is \emph{still} if $h=0$.

\subsection{Enriched stochastic flow map}\label{sec3EnrStoFloMap}

The goal of the associated enriched flow map, $\hat{\varphi}(M):\mathbb{R}\times\mathscr{Z}^{M+2}\to\mathscr{Z}^{M+2}$, is to push the system's state $s(t_i,\cdot\,)=(u(t_i,\cdot\,),\dot{u}(t_i,\cdot\,))$, and the first $M-1$ time derivatives of $\mathscr{f}(t_i,\cdot\,,s(t_i,\cdot\,))$ at $t=t_i$ (including the function itself).
This is the reason why the $k$-th component of $\hat{\varphi}(M)$ is given by:
\begin{equation*}
\hat{\varphi}^k(M)(h,\hat{s}(t_i,\cdot\,))=:\hat{s}^k(t_i+h,\cdot\,)=%
\begin{cases}
\varphi^k(M)(h,s(t_i,\cdot\,))& \text{for $k\in\{1,2\}$}\\
\mathrm{D}_t^{k-3}\mathscr{f}(t_i+h,\cdot\,,s(t_i+h,\cdot\,)) & \text{for $k\in\{3,4,\ldots,M+2\}$}
\end{cases}
\end{equation*}
where $\hat{s}=(u,\partial_tu,\ldots,\partial_t^{M+1}u)\in\prod_{j=1}^{M+2}\mathscr{T}(M-j+2)\otimes\mathscr{Z}$ is the enriched configuration state of the system over $\mathfrak{T}\times\Xi$.
Note that $\mathrm{D}_t^0\mathscr{f}:=\mathscr{f}$ is nothing but the function given by \eqref{eq3StoFloMap1000a}, and that $\{\mathrm{D}_t^{k-3}\mathscr{f}\}_{k=4}^{M+2}$ is the set of functions given by \eqref{eq3StoFloMap1090} with $j=k-3$.

\section{Overview of the FSC method}\label{sec3FSCmet}

The FSC method provides a new alternative for constructing an evolving finite-dimensional RFS efficiently in time.
The construction is said to be efficient because the dimensionality of the RFS is bounded from below by $n+2$ and from above by $n+M+1$, regardless of the dimensionality of the probability space, and where $n$ denotes the order of the governing ODE with respect to time, and $M$ the order of the stochastic flow map.

To illustrate the key idea behind the FSC method\footnote{For more details, please refer to \cite{esquivel2020flow,esquivel2021flow}.}, consider the stochastic dynamical system defined by \eqref{eq3ProSta1010}.
From \eqref{eq3StoFloMap1080}, we learn that the state of the system around $t=t_i$ can be expanded as
\begin{equation}\label{eq3FSCmet1010}
u(t,\xi)=\sum_{j=0}^M \frac{(t-t_i)^j}{j!}\partial_t^j u(t_i,\xi)\qquad\text{and}\qquad
\dot{u}(t,\xi)=\sum_{j=0}^M \frac{(t-t_i)^j}{j!}\partial_t^{j+1} u(t_i,\xi),
\end{equation}
where $h=t-t_i$, but with $t$ not necessarily fixed.

From \eqref{eq3FSCmet1010} we discover that the solution actually consists of two parts: the deterministic part $(t-t_i)^j/j!\in\mathscr{T}$, and the non-deterministic part $\partial_t^j u(t_i,\xi)\in\mathscr{Z}$.
If, for efficiency reasons, $\{\partial_t^j u(t_i,\cdot\,)\}_{j=0}^{M+1}$ is orthogonalized with respect to the measure defined in $\mathscr{Z}$, then $M+2$ orthogonal functions are obtained, namely: $\{\Psi_j\}_{j=1}^{M+2}$.
Hence, a more convenient way of writing \eqref{eq3FSCmet1010} would be
\begin{equation}\label{eq3FSCmet1050}
u(t,\xi)=\sum_{j=1}^{M+2} u^j(t)\,\Psi_j(\xi)\qquad\text{and}\qquad
\dot{u}(t,\xi)=\sum_{j=1}^{M+2} \dot{u}^j(t)\,\Psi_j(\xi).
\end{equation}
However, since the identically-equal-to-one function may not always be an element of $\{\Psi_j\}_{j=1}^{M+1}$, this calls to rewriting \eqref{eq3FSCmet1050} in the following manner:
\begin{equation}\label{eq3FSCmet1100}
u(t,\xi)=\sum_{j=0}^{M+2} u^j(t)\,\Psi_j(\xi)\qquad\text{and}\qquad
\dot{u}(t,\xi)=\sum_{j=0}^{M+2} \dot{u}^j(t)\,\Psi_j(\xi)
\end{equation}
with $\Psi_0\equiv1$ symbolizing the identically-equal-to-one function, as prescribed in Section \ref{sec3BasSpaMEFSC}.

Put differently, in FSC the orthogonalization process is always started with $\Psi_0\equiv1$ and then the germs $\{\partial_t^j u(t_i,\cdot\,)\}_{j=0}^{M+1}$ are orthogonalized with respect to each other to produce $M+3$ (orthogonal) basis vectors.
The resulting set of basis vectors is then used to construct a \emph{complete} RFS for the system's state at $t=t_i$:
\begin{equation*}
\mathscr{Z}^{[M+2]}=\mathrm{span}\{\Psi_j\}_{j=0}^{M+2}.
\end{equation*}

\begin{remark}
For the more general case, i.e.~for a dynamical system with an $n$-tuple state and driven by a stochastic flow map of order $M$, the aforementioned RFS is constructed as $\mathscr{Z}^{[n+M]}=\mathrm{span}\{\Psi_j\}_{j=0}^{n+M}$.
\end{remark}

However, since an exact representation of the system's state is not always needed to obtain accurate results, we highlight the fact that the FSC method is also capable of considering as few as $n+2$ basis vectors, but at the cost of losing some precision along the way.
Naturally, the longer the duration of the simulation, the greater this loss of precision will be.
In fact, this loss of precision is noticeable in the numerical problems explored in Section \ref{sec3NumRes}, because the duration of the simulations was chosen so that this loss would show up in the numerical results.
Therefore, a more effective way to write the RFS required for the simulation is to simply define it as
\begin{equation*}
\mathscr{Z}^{[P]}=\mathrm{span}\{\Psi_j\}_{j=0}^{P}\quad\text{with}\quad P\in\{n+1,n+2,\ldots, M+2\}.
\end{equation*}

\section{Multi-element flow-driven spectral chaos (ME-FSC) method}\label{sec3MEFSCmet}

This section is devoted to presenting the multi-element version of the FSC method in detail.
First, we review the multi-element concept behind the ME-FSC method, and then we present the computation of the (global) probability moments using the local information available on each random element.
The section is concluded with an outline of the ME-FSC scheme used in this work.

\subsection{Overview of the ME-FSC method}

Because the solution is in $\mathscr{U}$ by assumption, then $u$ can be represented with the Fourier series:
\begin{equation}\label{eq3MEFSCmet1000}
u(t,\xi)=\sum_{j=0}^\infty u^j(t)\,\Psi_j(\xi)\quad\text{on $\mathfrak{T}\times\Xi$},
\end{equation}
where $u^j$ is a temporal function in $\mathscr{T}(2)$ symbolizing the $j$-th random mode of $u$.

Now, let $\{(\Xi_e,\boldsymbol{\Xi}_e,\hat{\mu}_e)\}_{e=1}^E$ be a partitioned random space, and let $\{\mathscr{Z}_e\}_{e=1}^E$ be its associated partitioned RFS.
Then, a $p$-discretization for each of these $\mathscr{Z}_e$'s can be stipulated by letting $\mathscr{Z}_e^{[P_e]}=\mathrm{span}\{\Psi_{j.e}\}_{j=0}^{P_e}$ be a finite subspace of $\mathscr{Z}_e$ with $P_e+1\in\mathbb{N}_1$ denoting the dimensionality of $\mathscr{Z}_{e}^{[P_e]}$.

Moreover, if $u_{.e}(t,\cdot\,)=u(t,\cdot\,)|_{\Xi_e}$ is defined to be the restriction of $u(t,\cdot\,)$ to $\Xi_e$ for all $t\in\mathfrak{T}$, then an approximation $u_{.e}^{[P_e]}(t,\cdot\,)$ of $u_{.e}(t,\cdot\,)$ can be represented in $\mathscr{Z}_e^{[P_e]}$ by
\begin{equation}\label{eq3MEFSCmet1010}
u_{.e}(t,\xi)\approx u_{.e}^{[P_e]}(t,\xi)=\sum_{j=0}^{P_e}u{^j}_{\!.e}(t)\,\Psi_{j.e}(\xi)\equiv u{^j}_{\!.e}(t)\,\Psi_{j.e}(\xi)\quad\text{on $\mathfrak{T}\times\Xi_e$},
\end{equation}
where the summation sign has been omitted in the last equality for notational simplicity, and
the summation index is taken over $j\in\{0,1,\ldots,P_e\}$ unless indicated otherwise.
Similarly as before, $u{^j}_{\!.e}$ represents a temporal function in $\mathscr{T}(2)$ symbolizing the $j$-th random mode of $u_{.e}$.

\begin{remark}
A dot symbol is introduced in \eqref{eq3MEFSCmet1010} to distinguish between tensor indices and identification indices (aka ID indices).
Consequently, indices $j$ and $e$ are to be interpreted herein as tensor and identification indices, respectively.
Moreover, no hidden summation signs should ever be expected for identification indices.
\end{remark}

\begin{remark}
In writing \eqref{eq3MEFSCmet1010}, $\{\Psi_{j.e}\}_{j=0}^\infty$ was assumed to be a well-graded orthogonal basis in order to ensure that the approximation of $u_{.e}$ could be carried out that way.
This assumption does not represent an issue in this work, since all random bases utilized in ME-FSC are well-graded by construction.
\end{remark}

The problem now reduces to find $E$ independent solutions to \eqref{eq3ProSta1000} by looping across the random domain from $e=1$ to $e=E$ and using the following procedure.

Substituting \eqref{eq3MEFSCmet1010} into \eqref{eq3ProSta1000} gives
\begin{subequations}\label{eq3MEFSCmet2000}
\begin{align}
\mathcal{L}[u{^j}_{\!.e}\Psi_{j.e}]=f&\qquad\text{on $\mathfrak{T}\times\Xi_e$}\label{eq3MEFSCmet2000a}\\
\big\{\mathcal{B}_k[u{^j}_{\!.e}\Psi_{j.e}](0,\cdot\,)=b_k\big\}_{k=1}^n&\qquad\text{on $\{0\}\times\Xi_e$}.\label{eq3MEFSCmet2000b}
\end{align}
\end{subequations}

Projecting \eqref{eq3MEFSCmet2000} onto $\mathscr{Z}_e^{[P_e]}$ yields a system of $P_e+1$ ordinary differential equations of second order in the variable $t$, where the unknowns are the random modes $u\indices{^j_{.e}}=u\indices{^j_{.e}}(t)$ and their first $n-1$ time derivatives:
\begin{subequations}\label{eq3MEFSCmet3000}
\begin{align}
\Psi{^i}_{.e}\big[\mathcal{L}[u{^j}_{\!.e}\Psi_{j.e}]\big]=\Psi{^i}_{.e}[f]&\qquad\text{on $\mathfrak{T}$}\label{eq3MEFSCmet3000a}\\
\big\{\Psi{^i}_{.e}\big[\mathcal{B}_k[u{^j}_{\!.e}\Psi_{j.e}](0,\cdot\,)\big]=\Psi{^i}_{.e}[b_k]\big\}_{k=1}^n&\qquad\text{on $\{0\}$} \label{eq3MEFSCmet3000b}
\end{align}
\end{subequations}
with $i,j\in\{0,1,\ldots,P_e\}$.
Notice that this projection is done here by applying on both sides of each equation the linear functionals $\{\Psi{^i}_{.e}\in\mathscr{Z}_e'\}_{i=0}^{P_e}$ one by one.

For the specific case of a system given by \eqref{eq3ProSta1010}, it is clear that \eqref{eq3MEFSCmet3000} would reduce to
\begin{subequations}\label{eq3MEFSCmet3010}
\begin{align}
m\indices{^i_{j.e}}\ddot{u}\indices{^j_{.e}}+k\indices{^i_{j.e}}u\indices{^j_{.e}}=p\indices{^i_{.e}}&\qquad\text{on $\mathfrak{T}$}\tag{\ref{eq3MEFSCmet3000a}*}\\
\big\{u\indices{^i_{.e}}(0)=\mathscr{u}\indices{^i_{.e}},\,\dot{u}\indices{^i_{.e}}(0)=\mathscr{v}\indices{^i_{.e}}\big\}&\qquad\text{on $\{0\}$},\tag{\ref{eq3MEFSCmet3000b}*}
\end{align}
\end{subequations}
where $i,j\in\{0,1,\ldots,P_e\}$ (summation sign implied only over repeated index $j$), and:
\begin{gather*}
m\indices{^i_{j.e}}=\langle\Psi_{i.e},m\Psi_{j.e}\rangle_e/\langle\Psi_{i.e},\Psi_{i.e}\rangle_e,\quad
k\indices{^i_{j.e}}=\langle\Psi_{i.e},k\Psi_{j.e}\rangle_e/\langle\Psi_{i.e},\Psi_{i.e}\rangle_e,\\
p\indices{^i_{.e}}(t)=\langle\Psi_{i.e},p(t,\cdot\,)\rangle_e/\langle\Psi_{i.e},\Psi_{i.e}\rangle_e,\\
\mathscr{u}\indices{^i_{.e}}=\langle\Psi_{i.e},\mathscr{u}\rangle_e/\langle\Psi_{i.e},\Psi_{i.e}\rangle_e\quad\text{and}\quad
\mathscr{v}\indices{^i_{.e}}=\langle\Psi_{i.e},\mathscr{v}\rangle_e/\langle\Psi_{i.e},\Psi_{i.e}\rangle_e.
\end{gather*}

\subsection{Computation of probability moments}\label{sec3ComProMom}

In modeling notation, the system specified by \eqref{eq3ProSta1000} can be written as:
\begin{equation*}
y=\boldsymbol{\mathcal{M}}[u][x]\quad\text{subject to initial condition}\quad \boldsymbol{\mathcal{I}}[u],\\
\end{equation*}
where $\boldsymbol{\mathcal{M}}[u]:\mathscr{V}^r\to\mathscr{V}^s$ is as in \eqref{eq3ProSta1000a} and denotes the mathematical model of the system, $x=(x_1,\ldots,x_r):\mathfrak{T}\times\Xi\to\mathbb{R}^r$ is the $r$-tuple input of $\boldsymbol{\mathcal{M}}[u]$, $y=(y_1,\ldots,y_s):\mathfrak{T}\times\Xi\to\mathbb{R}^s$ is the $s$-tuple output of $\boldsymbol{\mathcal{M}}[u]$, and $\boldsymbol{\mathcal{I}}[u]$ is as in \eqref{eq3ProSta1000b} and denotes the initial condition of the system.

Now, define $z=y_k$ to be the $k$-th component of output $y=\boldsymbol{\mathcal{M}}[u][x]$.
If $z\in\mathscr{V}$, then it is clear that such a function can be represented approximately over the $e$-th random element using the following expansion:
\begin{equation*}
z_{.e}(t,\xi)\approx z_{.e}^{[P_e]}(t,\xi)=\sum_{j=0}^{P_e} z{^j}_{\!.e}(t)\,\Psi_{j.e}(\xi)\equiv z{^j}_{\!.e}(t)\,\Psi_{j.e}(\xi)\quad\text{on $\mathfrak{T}\times\Xi_e$},
\end{equation*}
where $z_{.e}(t,\cdot\,)=z(t,\cdot\,)|_{\Xi_e}$ denotes the restriction of $z(t,\cdot\,)$ to $\Xi_e$ for all $t\in\mathfrak{T}$.

In this section, we are interested in computing the mean and variance of $z$ over $\mathfrak{T}\times\Xi$ using the local information available on each random element.
To do so, we first make the following observation regarding the (local) mean and (local) variance of $z$ given that $\xi\in\Xi_e$.
As it should be easy to verify, the (local) mean of $z$ given that $\xi\in\Xi_e$, $\mathbf{E}_e[z]:\mathfrak{T}\to\mathbb{R}$, is nothing but
\begin{equation}\label{eq3ComProMom500}
\mathbf{E}_e[z](t)\equiv\mathbf{E}[z\mid\Xi_e](t):=\int z(t,\cdot\,)\,\mathrm{d}\hat{\mu}_e=
\int z_{.e}(t,\cdot\,)\,\mathrm{d}\hat{\mu}_e=
z{^0}_{.e}(t),
\end{equation}
and that the (local) variance of $z$ given that $\xi\in\Xi_e$, $\mathrm{Var}_e[z]:\mathfrak{T}\to\mathbb{R}_0^+$, is
\begin{multline}\label{eq3ComProMom550}
\mathrm{Var}_e[z](t)\equiv\mathrm{Var}[z\mid\Xi_e](t)
:=\int\big(z(t,\cdot\,)-\mathbf{E}_e[z](t)\big)^2\,\mathrm{d}\hat{\mu}_e
=\int z^2(t,\cdot\,)\,\mathrm{d}\hat{\mu}_e-\mathbf{E}_e[z]^2(t)\\
=\int z{_{.e}}^2(t,\cdot\,)\,\mathrm{d}\hat{\mu}_e-\mathbf{E}_e[z]^2(t)
=\sum_{j=1}^{P_e}\langle\Psi_{j.e},\Psi_{j.e}\rangle_e\,z{^j}_{\!.e}(t)\,z{^j}_{\!.e}(t).
\end{multline}

Therefore, the (global) mean of $z$, $\mathbf{E}[z]:\mathfrak{T}\to\mathbb{R}$, which is given by
\begin{equation*}
\mathbf{E}[z](t):=\int z(t,\cdot\,)\,\mathrm{d}\mu
=\sum_{e=1}^E\int_{\Xi_e} z(t,\cdot\,)\,\mathrm{d}\mu
=\sum_{e=1}^E\int z(t,\cdot\,)\,\mathrm{d}\mu_e,
\end{equation*}
simplifies to:
\begin{equation}\label{eq3ComProMom1000}
\mathbf{E}[z](t)
=\sum_{e=1}^E \mu(\Xi_e)\int z(t,\cdot\,)\,\mathrm{d}\hat{\mu}_e
=\sum_{e=1}^E\mu(\Xi_e)\,\mathbf{E}_e[z](t).
\end{equation}
This final form of $\mathbf{E}[z]$ is also known as the \emph{law of total expectation}.

Moreover, the (global) variance of $z$, $\mathrm{Var}[z]:\mathfrak{T}\to\mathbb{R}_0^+$, is given by
\begin{equation}\label{eq3ComProMom2000}
\mathrm{Var}[z](t)
:=\int\big(z(t,\cdot\,)-\mathbf{E}[z](t)\big)^2\,\mathrm{d}\mu
=\int z^2(t,\cdot\,)\,\mathrm{d}\mu-\mathbf{E}[z]^2(t)
=\sum_{e=1}^E\int_{\Xi_e} z^2(t,\cdot\,)\,\mathrm{d}\mu-\mathbf{E}[z]^2(t),
\end{equation}
The resulting expression can be further simplified if we recognize that \emph{first}:
\begin{equation*}
\int_{\Xi_e} z^2(t,\cdot\,)\,\mathrm{d}\mu
=\mu(\Xi_e)\int z^2(t,\cdot\,)\,\mathrm{d}\hat{\mu}_e
=\mu(\Xi_e)\int z{_{.e}}^2(t,\cdot\,)\,\mathrm{d}\hat{\mu}_e
=\mu(\Xi_e)\,\big(\mathrm{Var}_e[z](t)+\mathbf{E}_e[z]^2(t)\big)
\end{equation*}
and that \emph{second}:
\begin{multline*}
\mathbf{E}[z]^2(t)=\sum_{e_1=1}^E\sum_{e_2=1}^E\mu(\Xi_{e_1})\,\mu(\Xi_{e_2})\,\mathbf{E}_{e_1}[z](t)\,\mathbf{E}_{e_2}[z](t)\\
=\sum_{e=1}^E\mu(\Xi_e)^2\,\mathbf{E}_e[z]^2(t)
+2\sum_{e_1=2}^E\sum_{e_2=1}^{e_1-1}\mu(\Xi_{e_1})\,\mu(\Xi_{e_2})\,\mathbf{E}_{e_1}[z](t)\,\mathbf{E}_{e_2}[z](t).
\end{multline*}
Hence, by replacing these two expressions into \eqref{eq3ComProMom2000} and simplifying, we get
\begin{multline}\label{eq3ComProMom2050}
\mathrm{Var}[z](t)
=\sum_{e=1}^E\mu(\Xi_e)\,\mathrm{Var}_e[z](t)\\
+\sum_{e=1}^E\mu(\Xi_e)\,(1-\mu(\Xi_e))\,\mathbf{E}_e[z]^2(t)
-2\sum_{e_1=2}^E\sum_{e_2=1}^{e_1-1}\mu(\Xi_{e_1})\,\mu(\Xi_{e_2})\,\mathbf{E}_{e_1}[z](t)\,\mathbf{E}_{e_2}[z](t).
\end{multline}
This final form of $\mathrm{Var}[z]$ is also known as the \emph{law of total variance}.

\subsection{ME-FSC scheme}\label{sec3MEFSCsch}

\begin{figure}
\centering
\includegraphics[]{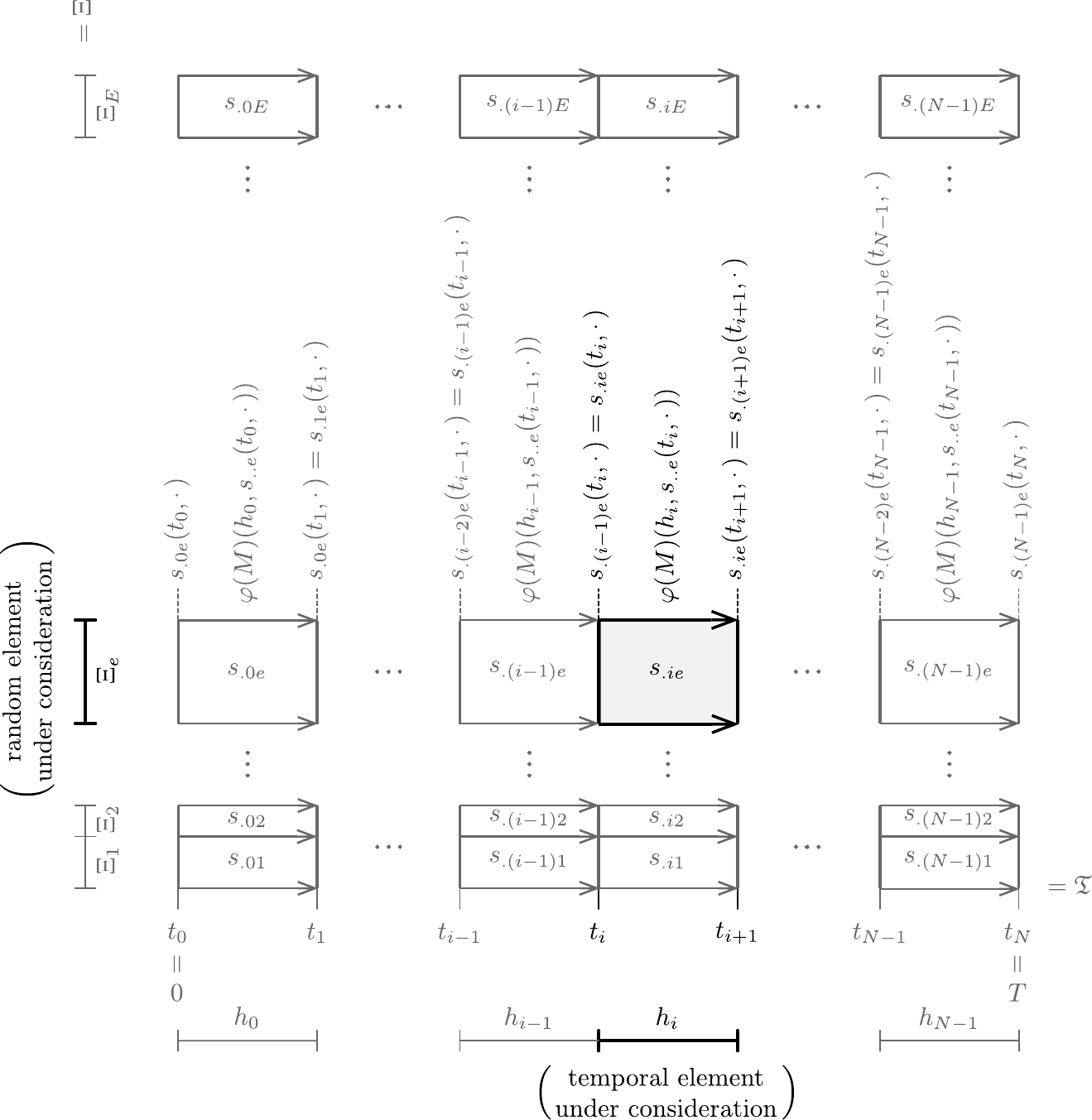}
\caption{Evolution of a dynamical system via a stochastic flow map of order $M$ (with $h_i>0$) over a partitioned random domain}
\label{fig3StochasticFlowMap}
\end{figure}

Consider the stochastic dynamical system given by \eqref{eq3ProSta1000}.
Let $\{\mathfrak{T}_i\times\Xi_e\}_{i=0,e=1}^{N-1,E}$ be a partition of the temporal-random domain $\mathfrak{T}\times\Xi$ (aka system's domain), where $\mathfrak{T}_i\neq\emptyset$ is the $i$-th element of the temporal domain, $\Xi_e\neq\emptyset$ is the $e$-th element of the random domain, and $\mathfrak{T}_i\times\Xi_e$ is the $(i,e)$-th element of the system's domain.
It is worth mentioning that this partition gives rise to the partitioned random space $\{(\Xi_e,\boldsymbol{\Xi}_e,\hat{\mu}_e)\}_{e=1}^E$ and associated partitioned RFS $\{\mathscr{Z}_e\}_{e=1}^E$, as defined in Section \ref{sec3SetNot}.
Moreover, for notational convenience, we define $s_{.ie}=s|_{\mathrm{cl}(\mathfrak{T}_i)\times\Xi_e}$ to be the restriction of $s$ to $\mathfrak{R}_{ie}:=\mathrm{cl}(\mathfrak{T}_i)\times\Xi_e$, and $s_{.\cdot e}=s|_{\mathfrak{T}\times\Xi_e}$ to be the restriction of $s$ to $\mathfrak{T}\times\Xi_e$.
For illustration purposes, Fig.~\ref{fig3StochasticFlowMap} depicts the evolution of a dynamical system by applying successively a stochastic flow map of order $M$ over a partitioned random domain with $E$ elements.

In this sense, if we assume that the system is driven by a stochastic flow map of order $M$, we can proceed as follows.
(Recall that $n$ denotes the order of the governing ODE with respect to time, as displayed in \eqref{eq3ProSta1000}.)

\begin{enumerate}
\item Loop across the temporal domain from $i=0$ to $i=N-1$.
\begin{enumerate}[i.]
\item Loop across the random domain from $e=1$ to $e=E$.
Note that this loop can be parallelized from a computational standpoint since each iteration is independent of the others.
\begin{enumerate}[(a)]
\item Define a solution representation for the configuration state $s_{.ie}$ in the following way.
\begin{itemize}
\item Take $\Phi_{0.ie}\equiv1$ and $\{\Phi_{j.ie}:=\hat{\varphi}^j(M)(0,\hat{s}_{.\cdot e}(t_i,\cdot\,))\}_{j=1}^{P_e}$ to be an ordered set of linearly independent functions in $\mathscr{Z}_e$ with $n+1\leq P\leq n+M$.
Note that $\hat{\varphi}(M)(0,\hat{s}_{.\cdot e}(t_i,\cdot\,))\equiv\hat{s}_{.ie}(t_i,\cdot\,)=\hat{s}_{.(i-1)e}(t_i,\cdot\,)$ for $i\geq 1$.
However, if $i=0$, then $\hat{\varphi}(M)(0,\hat{s}_{.\cdot e}(t_0,\cdot\,))\equiv\hat{s}_{.\cdot e}(0,\cdot\,)$.
\item Orthogonalize the set $\{\Phi_{j.ie}\}_{j=0}^{P_e}$ using the Gram-Schmidt process \cite{cheney2010linear}, so that the resulting set $\{\Psi_{j.ie}\}_{j=0}^{P_e}$ is an orthogonal basis in $\mathscr{Z}_e$.
Specifically, for $j\in\{0,1,\ldots,P_e\}$:
\begin{equation*}
\Psi_{j.ie}:=\Phi_{j.ie}-\sum_{k=0}^{j-1}\frac{\langle\Phi_{j.ie},\Psi_{k.ie}\rangle_e}{\langle\Psi_{k.ie},\Psi_{k.ie}\rangle_e}\Psi_{k.ie}.
\end{equation*}
\item Define $\mathscr{Z}_{ie}^{[P_e]}=\mathrm{span}\{\Psi_{j.ie}\}_{j=0}^{P_e}$ to be an $(h,p)$-discretization of $\mathscr{Z}$ over the region $\mathfrak{R}_{ie}$.
Then, because $\mathscr{Z}_{ie}^{[P_e]}$ is an evolving function space over $\Xi_e$, expansion \eqref{eq3MEFSCmet1010} is now to be read as:
\begin{equation}\label{eq3MEFSCsch3000}
u_{.ie}(t,\xi)\approx u^{[P_e]}_{.ie}(t,\xi)=\sum_{j=0}^{P_e}u{^j}_{\!.ie}(t)\,\Psi_{j.ie}(\xi)\equiv u{^j}_{\!.ie}(t)\,\Psi_{j.ie}(\xi).\tag{\ref{eq3MEFSCmet1010}'}
\end{equation}
Thus, to compute the $l$-th component of the configuration state, $s{^l}_{\!.ie}$, we simply need to take the $(l-1)$-th time derivative of \eqref{eq3MEFSCsch3000}.
Here $l\in\{1,2,\ldots,n\}$.
\end{itemize}

\item Transfer at $t=t_i$ the random modes from the old definition of the system's configuration state
\begin{equation*}
s_{.(i-1)e}(t_i,\cdot\,)=(u_{.(i-1)e}(t_i,\cdot\,),\partial_tu_{.(i-1)e}(t_i,\cdot\,),\ldots,\partial_t^{n-1}u_{.(i-1)e}(t_i,\cdot\,))
\end{equation*}
to the new definition of the system's configuration state
\begin{equation*}
s_{.ie}(t_i,\cdot\,)=(u_{.ie}(t_i,\cdot\,),\partial_tu_{.ie}(t_i,\cdot\,),\ldots,\partial_t^{n-1}u_{.ie}(t_i,\cdot\,)),
\end{equation*}
given that $i\geq1$.
To do so, we implement the FSC-2 approach presented in \cite{esquivel2020flow} to obtain the random modes of each of the components of $s_{.ie}$ at $t=t_i$.
The reason why the FSC-2 is chosen over the FSC-1 is that the former can transfer the probability information exactly at any given time.
Thus, by resorting to Theorem 1 of \cite{esquivel2020flow}, one can show that the $j$-th random mode of the $l$-th component of $s_{.ie}(t_i,\cdot\,)$ is given by:
\begin{equation}\label{eq3MEFSCsch3050}
(s^l){^j}_{\!.ie}(t_i)=
\begin{cases}
\mathbf{E}[\Phi_{l.ie}] & \text{for $j=0$}\\[1.25ex]
\displaystyle\frac{\det\triangle_j(l)}{\det\square_j} & \text{for $0<j<l$}\\[1.25ex]
1 & \text{for $j=l$}\\
0 & \text{otherwise,}
\end{cases}\tag{\ref{eq3MEFSCmet3000b}'}
\end{equation}
from where we have taken $\{\Phi_{l.ie}:=\varphi^l(M)(0,s_{.\cdot e}(t_i,\cdot\,))\equiv s{^l}_{\!.ie}(t_i,\cdot\,)=s{^l}_{\!.(i-1)e}(t_i,\cdot\,)\}_{l=1}^n$, and
\begin{align*}
\square_j=&
\begin{bmatrix}
\mathrm{Cov}[\Phi_{1.ie},\Phi_{1.ie}] & \cdots & \mathrm{Cov}[\Phi_{1.ie},\Phi_{j.ie}]\\
\vdots & \ddots & \vdots\\
\mathrm{Cov}[\Phi_{j.ie},\Phi_{1.ie}] & \cdots & \mathrm{Cov}[\Phi_{j.ie},\Phi_{j.ie}]
\end{bmatrix}\in\mathscr{M}(j\times j,\mathbb{R}),\\[1ex]
\triangle_j(l)=&
\begin{bmatrix}
\mathrm{Cov}[\Phi_{1.ie},\Phi_{1.ie}] & \cdots & \mathrm{Cov}[\Phi_{1.ie},\Phi_{j.ie}]\\
\vdots & \ddots & \vdots\\
\mathrm{Cov}[\Phi_{(j-1).ie},\Phi_1] & \cdots & \mathrm{Cov}[\Phi_{(j-1).ie},\Phi_{j.ie}]\\
\mathrm{Cov}[\Phi_{l.ie},\Phi_{1.ie}] & \cdots & \mathrm{Cov}[\Phi_{l.ie},\Phi_{j.ie}]
\end{bmatrix}\in\mathscr{M}(j\times j,\mathbb{R})
\end{align*}
with $\triangle_1(l)=\mathrm{Cov}[\Phi_{l.ie},\Phi_{1.ie}]\in\mathbb{R}$ and
\begin{equation*}
\triangle_2(l)=
\begin{bmatrix}
\mathrm{Cov}[\Phi_{1.ie},\Phi_{1.ie}] & \mathrm{Cov}[\Phi_{1.ie},\Phi_{2.ie}]\\
\mathrm{Cov}[\Phi_{l.ie},\Phi_{1.ie}] & \mathrm{Cov}[\Phi_{l.ie},\Phi_{2.ie}]
\end{bmatrix}\in\mathscr{M}(2\times2,\mathbb{R})
\end{equation*}
as special cases of $\triangle_j$ by definition.

\item Substitute \eqref{eq3MEFSCsch3000} into \eqref{eq3ProSta1000} to obtain \eqref{eq3MEFSCmet2000}.

\item Project \eqref{eq3MEFSCmet2000a} onto $\mathscr{Z}_{ie}^{[P_e]}$ to obtain \eqref{eq3MEFSCmet3000a} subject to \eqref{eq3MEFSCsch3050}.
Note that if $i=0$, \eqref{eq3MEFSCmet3000a} is subject to \eqref{eq3MEFSCmet3000b}.

\item Integrate \eqref{eq3MEFSCmet3000} over time, provided that a suitable time integration method has been chosen for the system of equations in hand.
This step requires finding the random modes of each of the components of the configuration state $s_{.ie}$ at $t=t_{i+1}$; that is, $\{(s^l){^j}_{.ie}(t_{i+1})\}_{l=1,j=0}^{n,P_e}$

\item Compute both the (local) mean and the (local) variance of each of the components of output $y=\boldsymbol{\mathcal{M}}[x][u]$ over $\mathfrak{R}_{ie}$, by resorting to the formulas prescribed by \eqref{eq3ComProMom500} and \eqref{eq3ComProMom550}.
\end{enumerate}

\item Aggregate results to compute over $\mathfrak{R}_i:=\mathfrak{T}_i\times\Xi$ the (global) mean and (global) variance of $y=\boldsymbol{\mathcal{M}}[x][u]$ using the formulas stipulated by \eqref{eq3ComProMom1000} and \eqref{eq3ComProMom2050}.
\end{enumerate}

\item Post-process results.
\end{enumerate}

\section{Numerical results}\label{sec3NumRes}

The local and global errors, $\epsilon:\mathscr{T}\to\mathscr{T}$ and $\epsilon_G:\mathscr{T}\to\mathbb{R}$, are defined using the following expressions:
\begin{gather*}
\epsilon[f](t)=|f(t)-f_\mathrm{exact}(t)|\\
\epsilon_G[f]=\frac{1}{T}\int_\mathfrak{T}|f(t)-f_\mathrm{exact}(t)|\,\mathrm{d}t\approx\frac{\Delta t}{T}\sum_{i=0}^N|f(t_i)-f_\mathrm{exact}(t_i)|,
\end{gather*}
where $t_i\in\mathfrak{T}$ is the time instant of the simulation, $\Delta t$ is the associated time-step size, and $N$ is the number of time steps utilized in the simulation (such that $t_0=0$ and $t_N=N\,\Delta t=T$).

In an effort to reduce significantly the source of errors coming from the discretization of $\mathscr{T}$, the time-step size used in the simulations is taken as $\Delta t=0.001$ s for Problems 1 and 2 and $\Delta t=0.005$ s for Problems 3 and 4.
To integrate \eqref{eq3MEFSCmet3000} over time, we use the RK4 method over each random element, and in order to obtain accurate results, the random function space is updated at every time step.
For simplicity, the partition of the random domain is such that all elements are the same size.
The temporal domain employed in Problems 1, 2 and 3 is $\mathfrak{T}=[0,150]$ s, and in Problem 4 is $\mathfrak{T}=[0,50]$ s.
To prevent ill-conditioned matrices in the early times of a simulation with deterministic initial conditions, the gPC method is used till the first second of the simulation.
This is to ensure that the system's state is well-developed for the analysis with ME-FSC.
Specifically, for Problems 1 and 2 we employ the gPC method with $P=7$, and for Problem 3 we employ the gPC method with $P=9$.

The inner products are computed using a set of Legendre quadrature rules on each random axis defined by 
\begin{equation*}
\mathrm{Uniform}\sim\text{Gauss-Legendre (10 points/element)}.
\end{equation*}
This means that for distributions other than uniform, the probability distribution function must be included in the integrand in order to correctly obtain the numerical value of the inner product.

All problems are run using Apple's Foundation and Accelerate frameworks \cite{apple2021frameworks} on a 2020 MacBook Air with Apple M1 chip (\emph{8-Core} CPU at \emph{3.20 GHz}, \emph{8-Core} GPU, \emph{16-Core} Neural Engine, and \emph{16 GB} unified memory) and \emph{1 TB Apple-Fabric SSD} storage (APFS-formatted), running macOS Big Sur (version 11.2).
The code is written entirely in the Swift 5.3 language \cite{swift2021language}.

\subsection[Problem 1 (linear system)]{Problem 1: A linear system governed by a 2nd-order stochastic ODE}

We first consider the problem of an undamped single-degree-of-freedom system under free vibration.
The law of motion for this system is defined by
\begin{equation*}
m\ddot{u}+ku=0,
\end{equation*}
where the mass of the system is $m=100$ kg, and the stiffness of the system, $k:\Xi\to\mathbb{R}^+$, is stochastic and given by $k(\xi)=\xi$.
The system has an initial displacement of $u(0,\cdot\,)=\mathscr{u}\equiv 0.05$ m, and an initial velocity of $\dot{u}(0,\cdot\,)=\mathscr{v}\equiv 0.20$ m/s.
Formally, one can express this problem in the following way:

Find the displacement of the system $u:\mathfrak{T}\times\Xi\to\mathbb{R}$ in $\mathscr{U}$, such that ($\mu$-a.e.):
\begin{subequations}
\label{eq3NumRes1000}
\begin{align}
m\ddot{u}+ku=0 &\qquad\text{on $\mathfrak{T}\times\Xi$}\label{eq3NumRes1000a}\\
\big\{u(0,\cdot\,)=\mathscr{u},\,\dot{u}(0,\cdot\,)=\mathscr{v}\big\} &\qquad\text{on $\{0\}\times\Xi$}.\label{eq3NumRes1000b}
\end{align}
\end{subequations}
This problem statement is similar to \eqref{eq3ProSta1010}, with the only difference being that $p\equiv0$ and $m,\mathscr{u},\mathscr{v}$ are real numbers.

Three different probability distributions are considered for $\xi$.
The first distribution is a \emph{uniform distribution}, $\mathrm{Uniform}\sim\xi\in\Xi=[a,b]$, the second distribution is a \emph{beta distribution}, $\mathrm{Beta}(\alpha_1,\beta_1)\sim\xi\in\Xi=[a,b]$, and the third distribution is a \emph{gamma distribution}, $\mathrm{Gamma}(\alpha_2,\beta_2)\sim\xi\in\Xi=[a,\infty)$.
Because gamma is a distribution with non-compact support, the simulations are run with \emph{truncated gamma distribution}, $\mathrm{Gamma}(\alpha_2,\beta_2,b_2)\sim\xi\in\Xi=[a,b_2]$.
But still, the exact solution is obtained with the gamma distribution.
The parameters for each of the distributions mentioned above are taken as: $(a,b)=(340,460)$ N/m, $(\alpha_1,\beta_1)=(2,5)$ and $(\alpha_2,\beta_2,b_2)=(10,0.1,920\text{ N/m})$.

\begin{figure}
\centering
\begin{subfigure}[b]{0.495\textwidth}
\includegraphics[width=\textwidth]{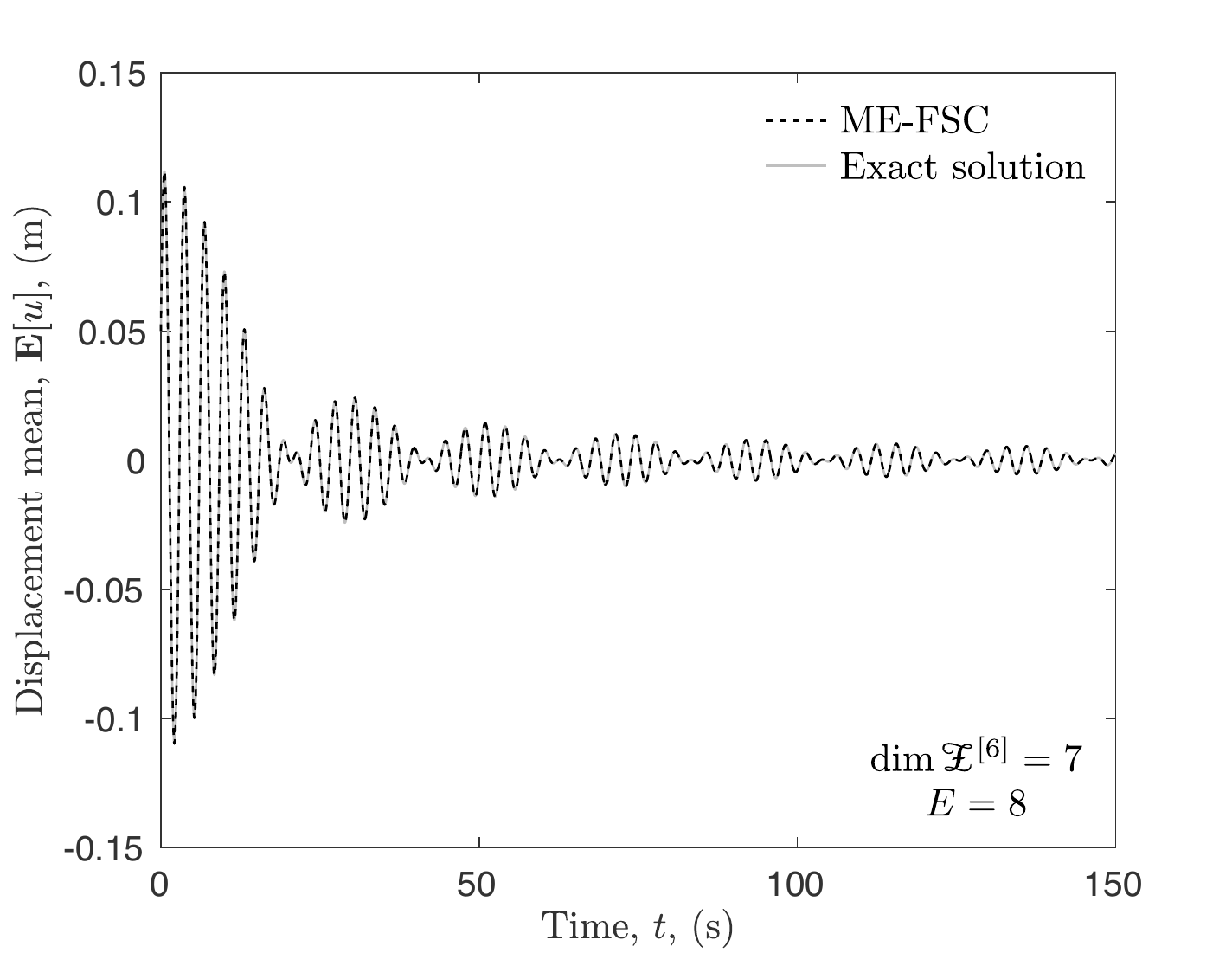}
\caption{Mean}
\label{fig3Problem1_Uniform_MEFSC_7BV_8E_Disp_Mean}
\end{subfigure}\hfill
\begin{subfigure}[b]{0.495\textwidth}
\includegraphics[width=\textwidth]{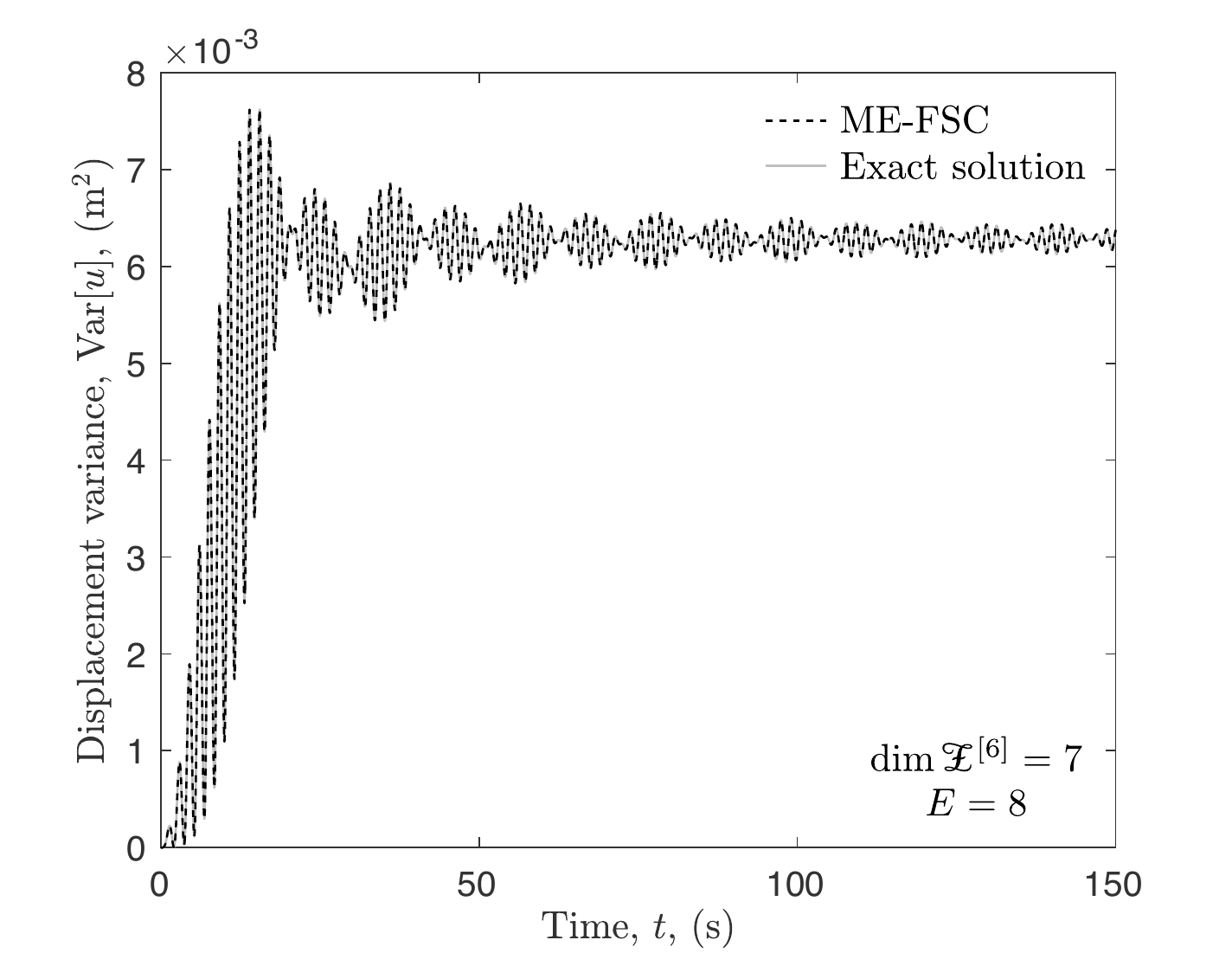}
\caption{Variance}
\label{fig3Problem1_Uniform_MEFSC_7BV_8E_Disp_Var}
\end{subfigure}
\caption{\emph{Problem 1} --- Evolution of $\mathbf{E}[u]$ and $\mathrm{Var}[u]$ for the case when the $(h,p)$-discretization level of RFS is $(P,E)=(6,8)$ and $\mu\sim\mathrm{Uniform}$}
\label{fig3Problem1_Uniform_MEFSC_7BV_8E_Disp}
\end{figure}

In Fig.~\ref{fig3Problem1_Uniform_MEFSC_7BV_8E_Disp} we show the evolution of the mean and variance of the system's displacement.
These results are obtained using 8 elements in the random domain and 7 basis vectors per element.
We see that when this particular discretization is used, the ME-FSC method is able to reproduce the exact response with high fidelity, explaining why the two plots look indistinguishable from each other.

% Problem 1 [1/3]:

\begin{figure}
\centering
\begin{subfigure}[b]{0.495\textwidth}
\includegraphics[width=\textwidth]{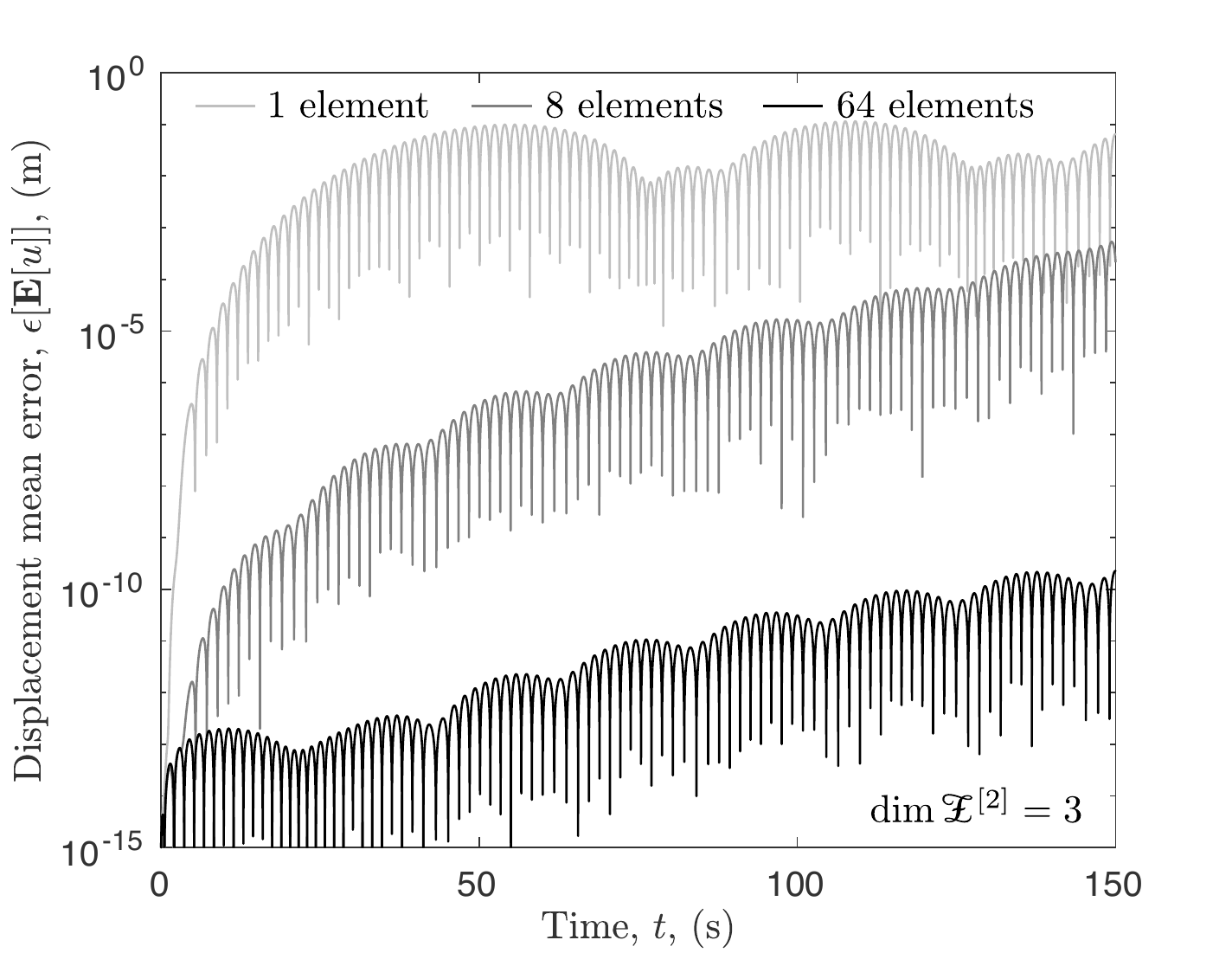}
\caption{Mean error for $\mathscr{Z}^{[2]}$}
\label{fig3Problem1_Uniform_MEFSC_3BV_Disp_Mean_LocalError}
\end{subfigure}\hfill
\begin{subfigure}[b]{0.495\textwidth}
\includegraphics[width=\textwidth]{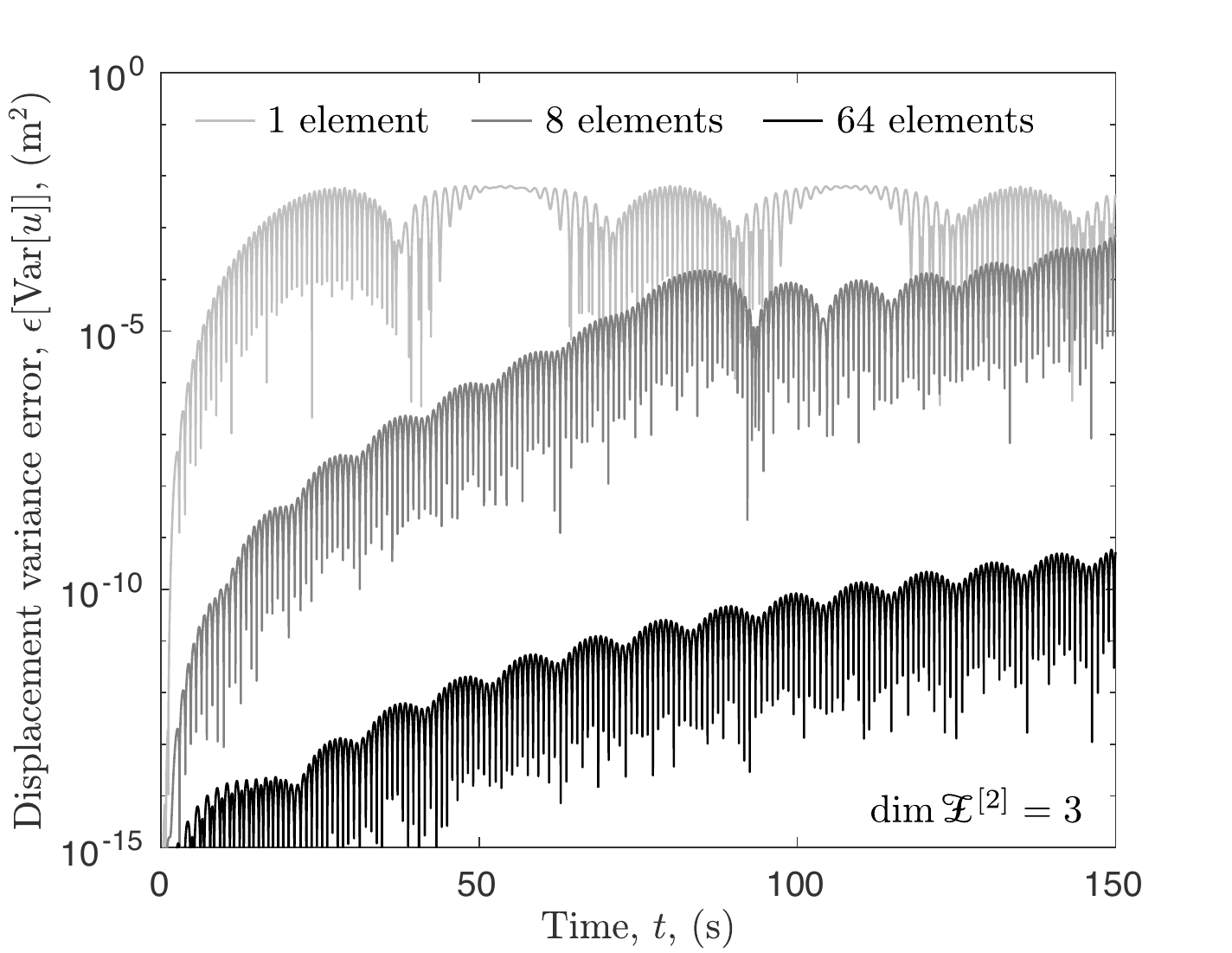}
\caption{Variance error for $\mathscr{Z}^{[2]}$}
\label{fig3Problem1_Uniform_MEFSC_3BV_Disp_Var_LocalError}
\end{subfigure}\quad
\begin{subfigure}[b]{0.495\textwidth}
\includegraphics[width=\textwidth]{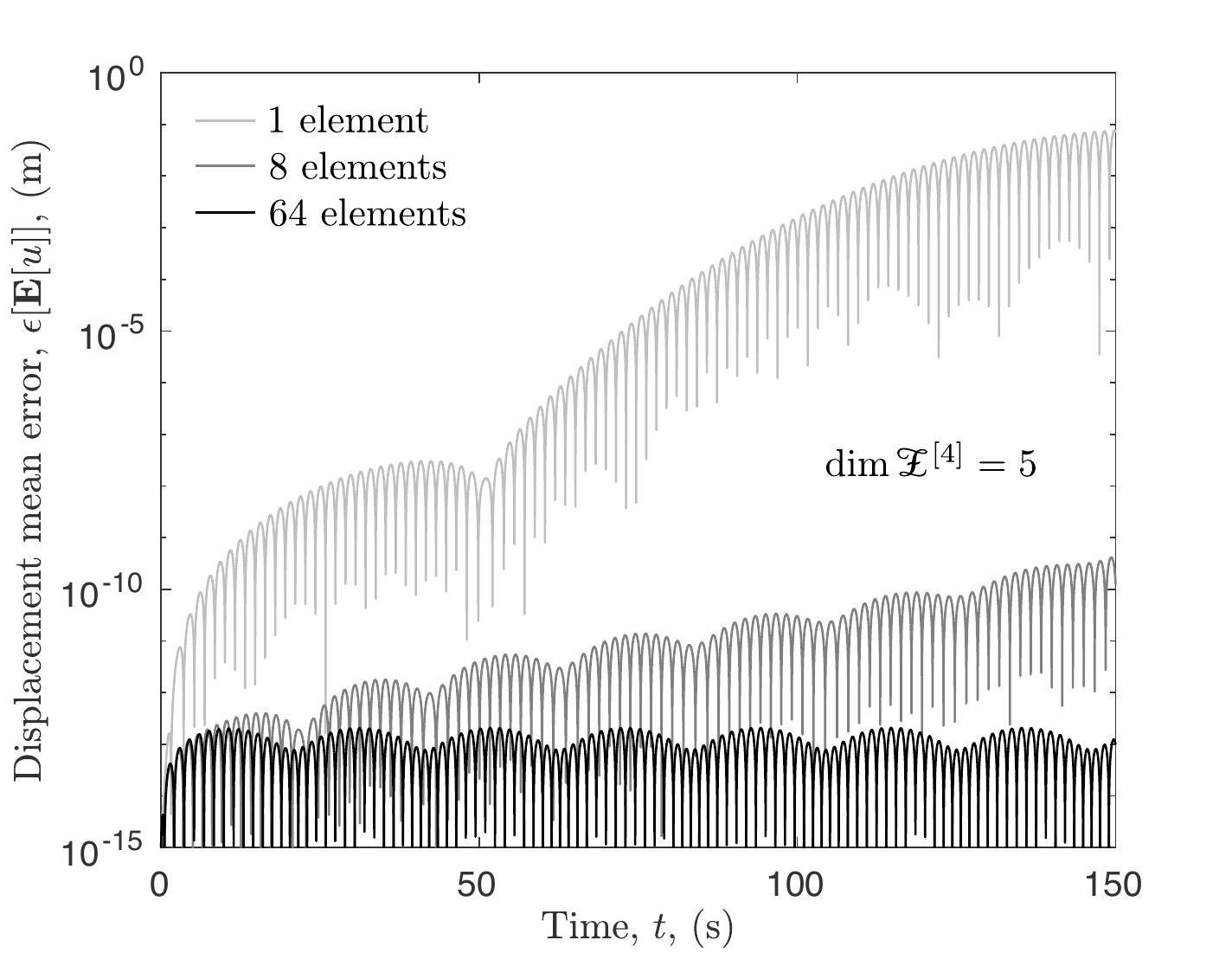}
\caption{Mean error for $\mathscr{Z}^{[4]}$}
\label{fig3Problem1_Uniform_MEFSC_5BV_Disp_Mean_LocalError}
\end{subfigure}\hfill
\begin{subfigure}[b]{0.495\textwidth}
\includegraphics[width=\textwidth]{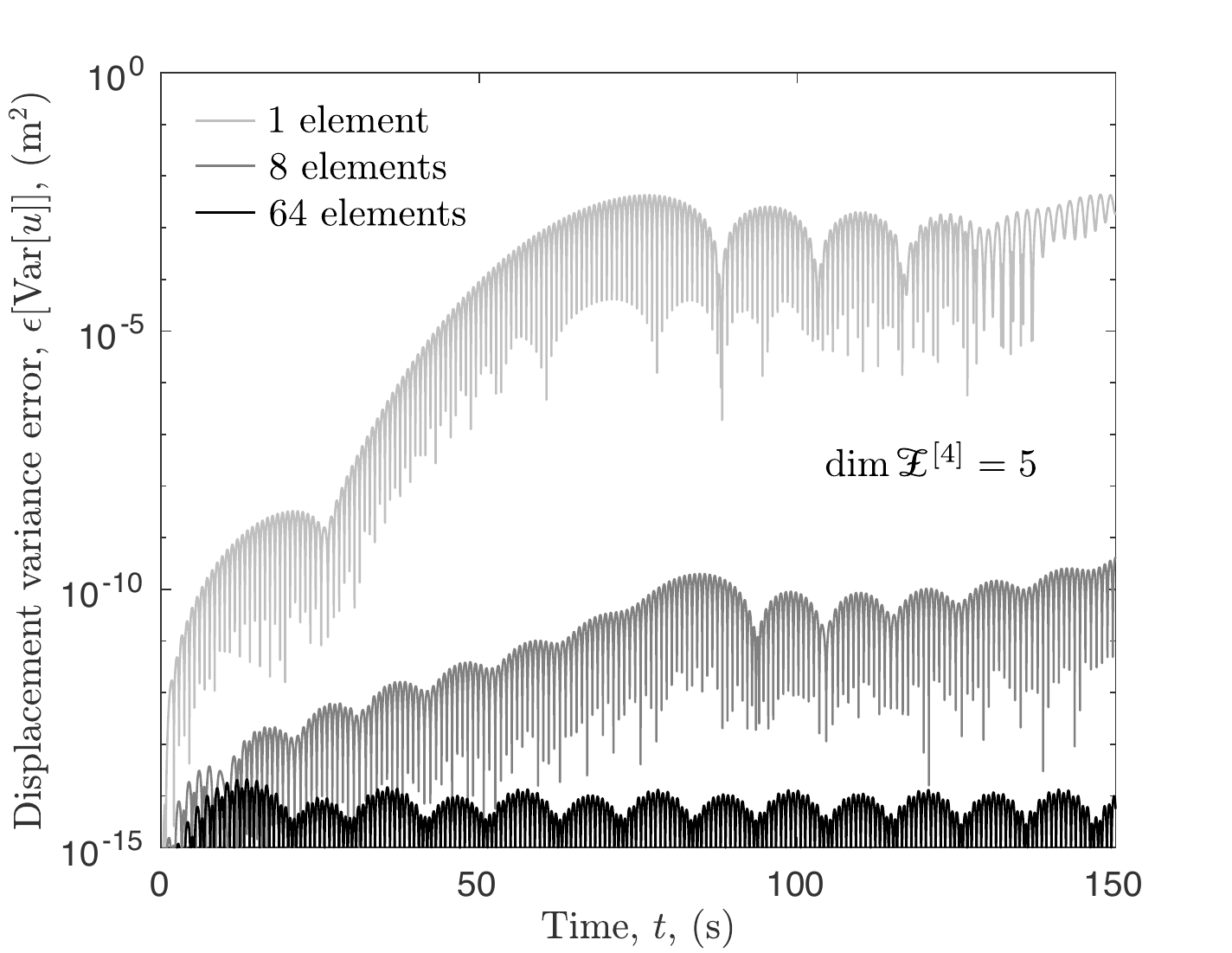}
\caption{Variance error for $\mathscr{Z}^{[4]}$}
\label{fig3Problem1_Uniform_MEFSC_5BV_Disp_Var_LocalError}
\end{subfigure}\quad
\begin{subfigure}[b]{0.495\textwidth}
\includegraphics[width=\textwidth]{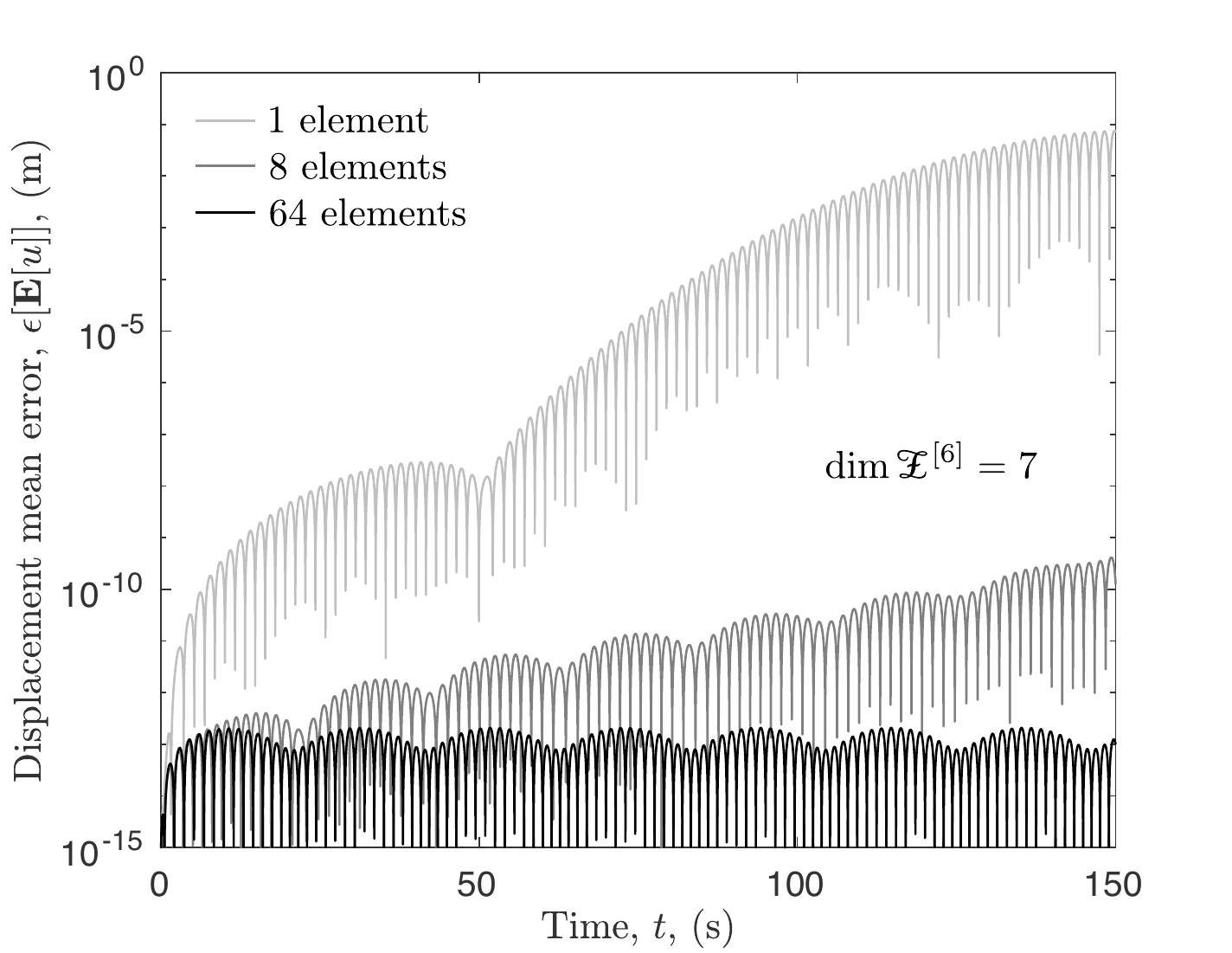}
\caption{Mean error for $\mathscr{Z}^{[6]}$}
\label{fig3Problem1_Uniform_MEFSC_7BV_Disp_Mean_LocalError}
\end{subfigure}\hfill
\begin{subfigure}[b]{0.495\textwidth}
\includegraphics[width=\textwidth]{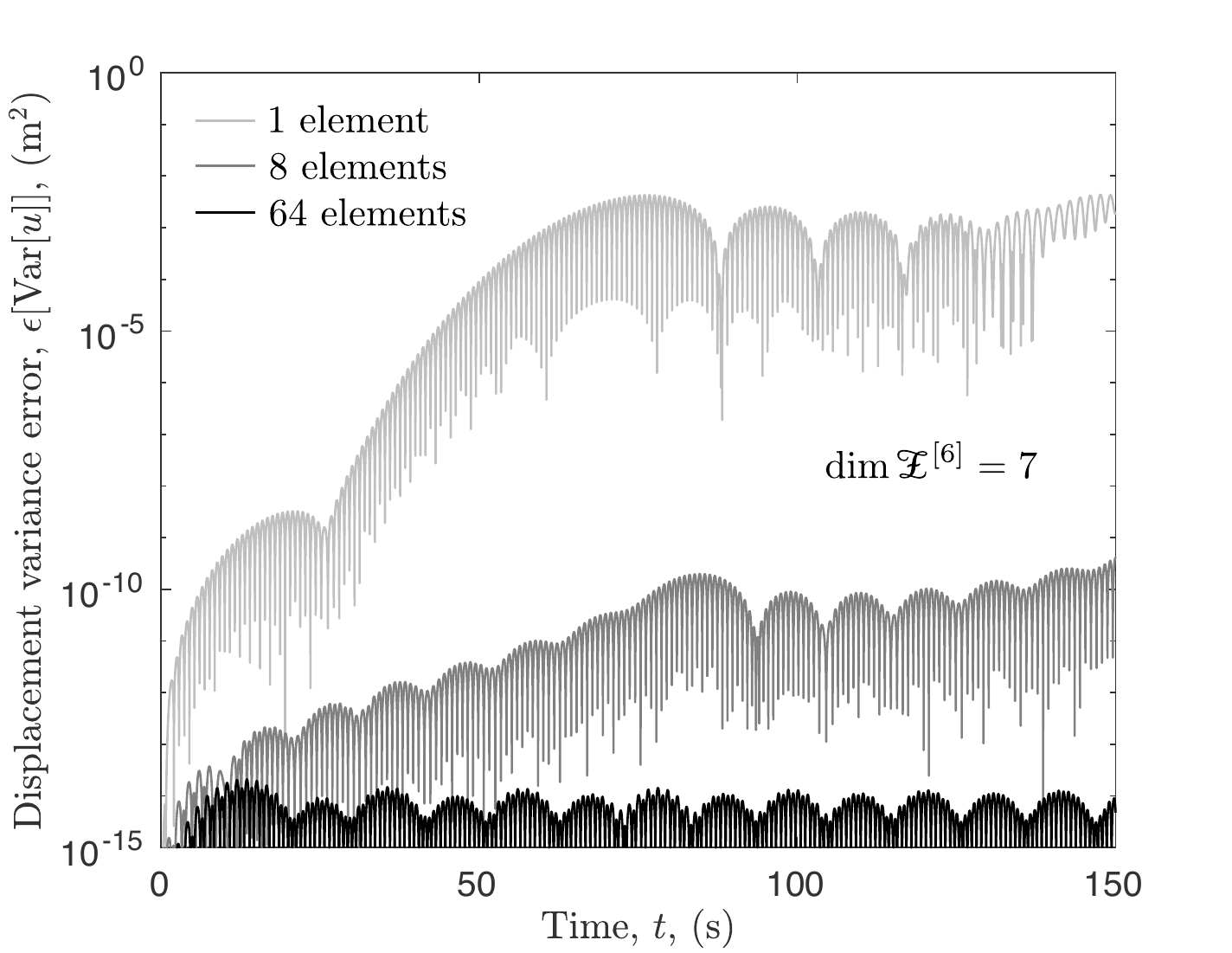}
\caption{Variance error for $\mathscr{Z}^{[6]}$}
\label{fig3Problem1_Uniform_MEFSC_7BV_Disp_Var_LocalError}
\end{subfigure}
\caption{\emph{Problem 1} --- Local error evolution of $\mathbf{E}[u]$ and $\mathrm{Var}[u]$ for different $(h,p)$-discretization levels of RFS and for $\mu\sim\mathrm{Uniform}$}
\label{fig3Problem1_Uniform_MEFSC_Disp_LocalError}
\end{figure}

Fig.~\ref{fig3Problem1_Uniform_MEFSC_Disp_LocalError} presents the local errors in mean and variance of the system's displacement using different numbers of elements and basis vectors.
For brevity, we only present the case when the probability measure is uniform, although similar convergence trends are also achievable when the probability measure is beta or gamma.
For sake of comparison, we also include the case when $P=2$, even though this is not allowed by the FSC method.\footnote{This is because in FSC the lower bound for $P$ is always $n+1$, where $n$ is the order of the governing ODE with respect to time.}
These plots show that the results get better if the number of elements is increased.
However, this is not always the case whenever the number of basis vectors is increased.
For instance, when the number of basis vectors is increased from 3 to 5, we observe that the results improve noticeably, but when the number of basis vectors is increased from 5 to 7 they do not.
This is chiefly because the quadrature rule employed in this work is not optimal---recall that Gaussian-based quadrature rules only ensure exponential convergence to the sought integral if the selected quadrature points are the roots of an orthogonal polynomial that is concordant with the measure defined on the integral's domain (which is manifestly not the case here).
The above assertion even holds for the case when 1 element is used in the simulations (i.e.~when there is no partition in the random domain), since distributing only 10 Legendre quadrature points in the random domain is not sufficient enough to allow the higher basis vectors to play a role in the global response of the system.
However, if more quadrature points were to be employed in the simulations, an improvement between the plots with $P=4$ and $P=6$ could have been discerned for the case of using 1 element, but virtually no improvement for the case of using more than 1 element.
It is worth pointing out that the accuracy of these results is limited by machine precision.
As a result, better results than those depicted in Fig.~\ref{fig3Problem1_Uniform_MEFSC_Disp_LocalError} are difficult to obtain for other values of $P$ and $E$.

% Problem 1 [2/3]:

\begin{figure}
\centering
\begin{subfigure}[b]{0.495\textwidth}
\includegraphics[width=\textwidth]{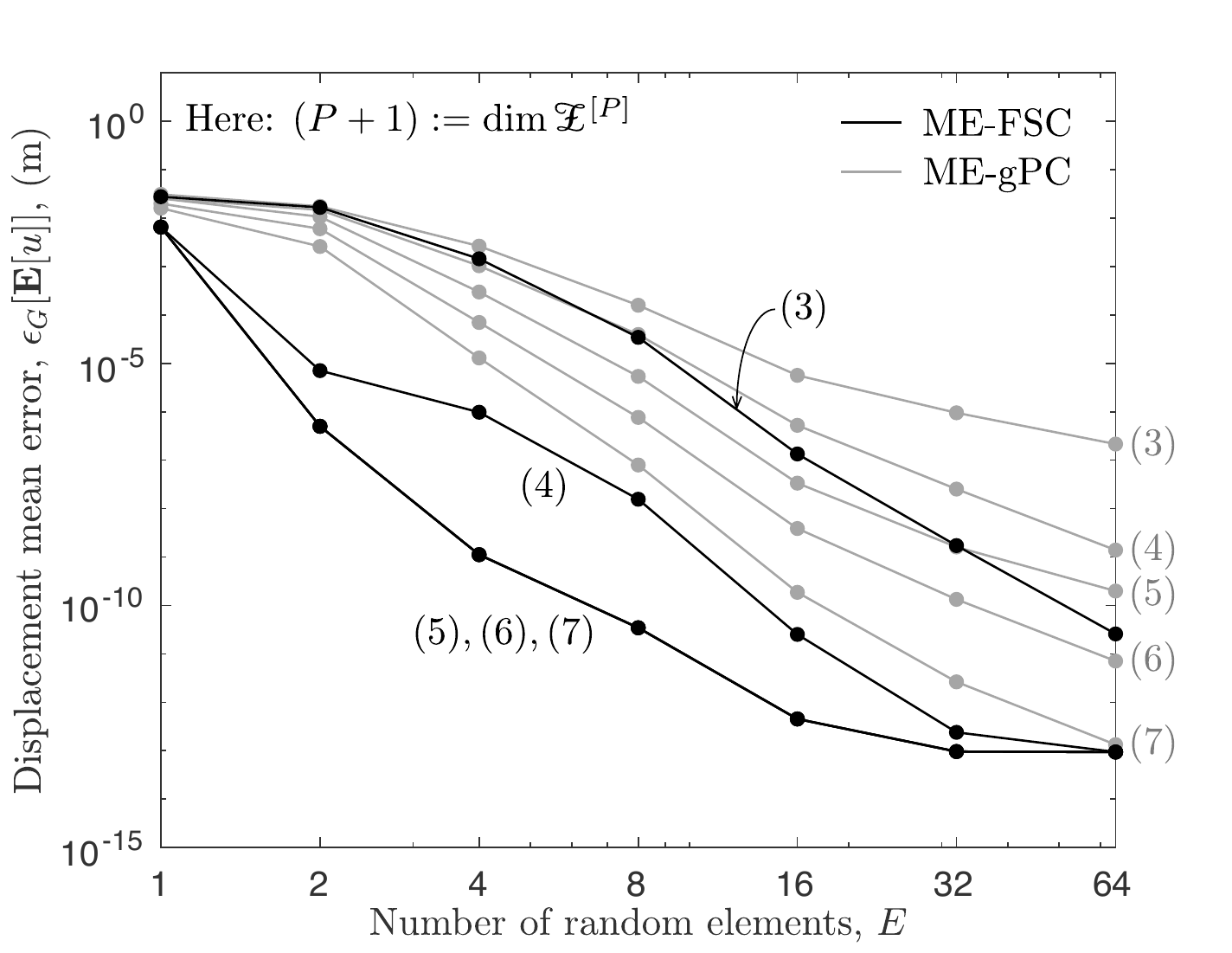}
\caption{Mean error for $\mu\sim\mathrm{Uniform}$}
\label{fig3Problem1_Uniform_MEFSC_Disp_Mean_GlobalError}
\end{subfigure}\hfill
\begin{subfigure}[b]{0.495\textwidth}
\includegraphics[width=\textwidth]{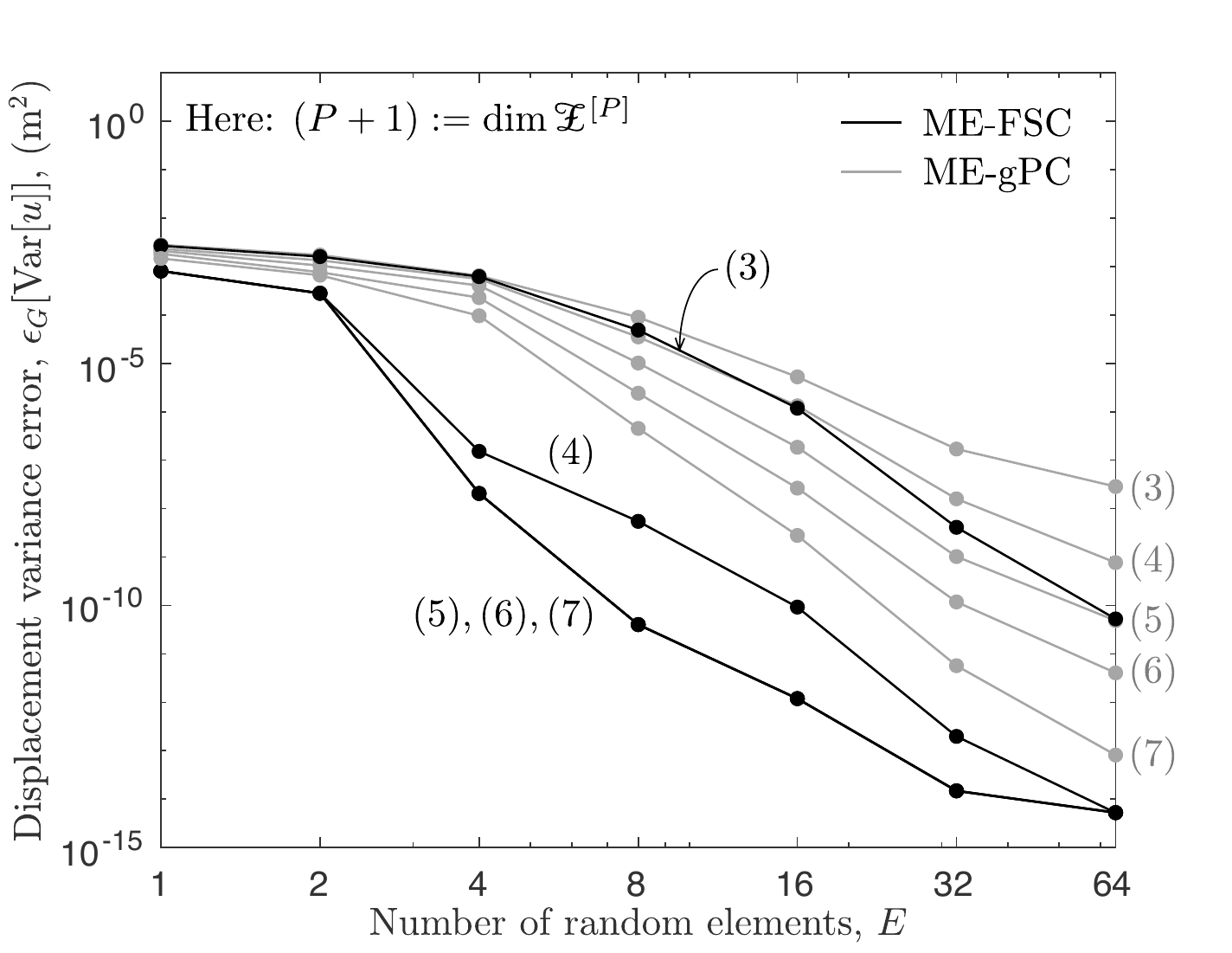}
\caption{Variance error for $\mu\sim\mathrm{Uniform}$}
\label{fig3Problem1_Uniform_MEFSC_Disp_Var_GlobalError}
\end{subfigure}\quad
\begin{subfigure}[b]{0.495\textwidth}
\includegraphics[width=\textwidth]{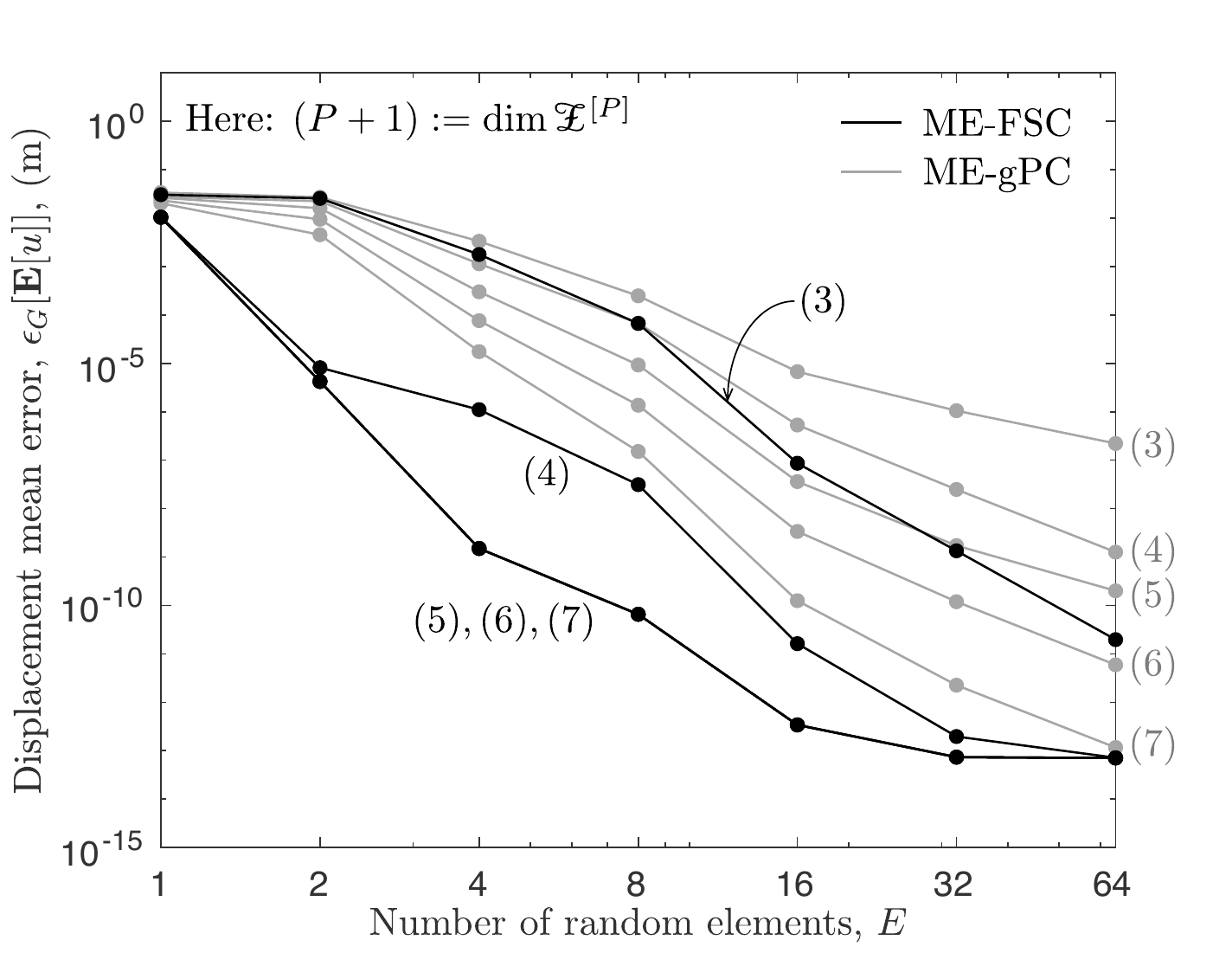}
\caption{Mean error for $\mu\sim\mathrm{Beta}$}
\label{fig3Problem1_Beta_MEFSC_Disp_Mean_GlobalError}
\end{subfigure}\hfill
\begin{subfigure}[b]{0.495\textwidth}
\includegraphics[width=\textwidth]{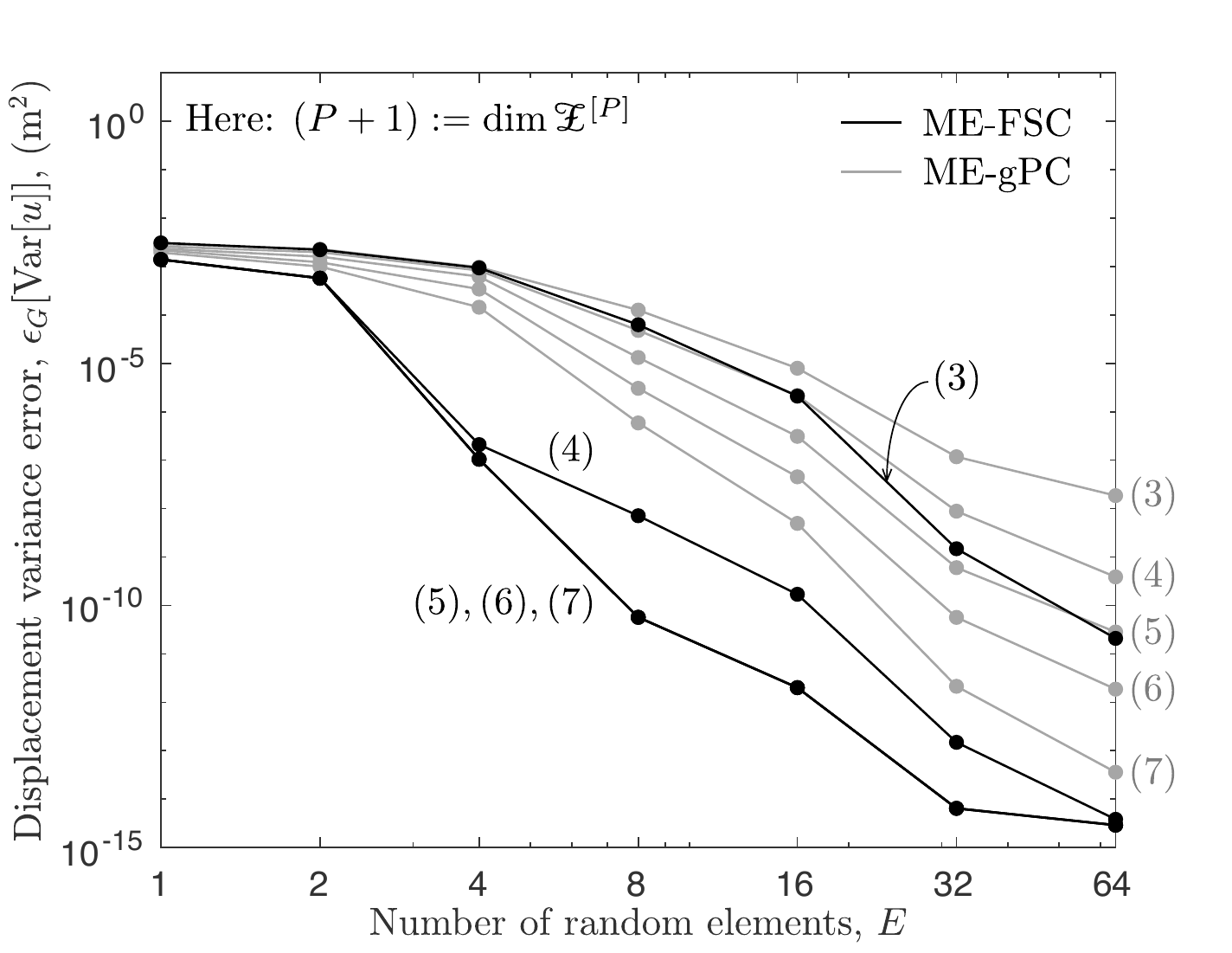}
\caption{Variance error for $\mu\sim\mathrm{Beta}$}
\label{fig3Problem1_Beta_MEFSC_Disp_Var_GlobalError}
\end{subfigure}\quad
\begin{subfigure}[b]{0.495\textwidth}
\includegraphics[width=\textwidth]{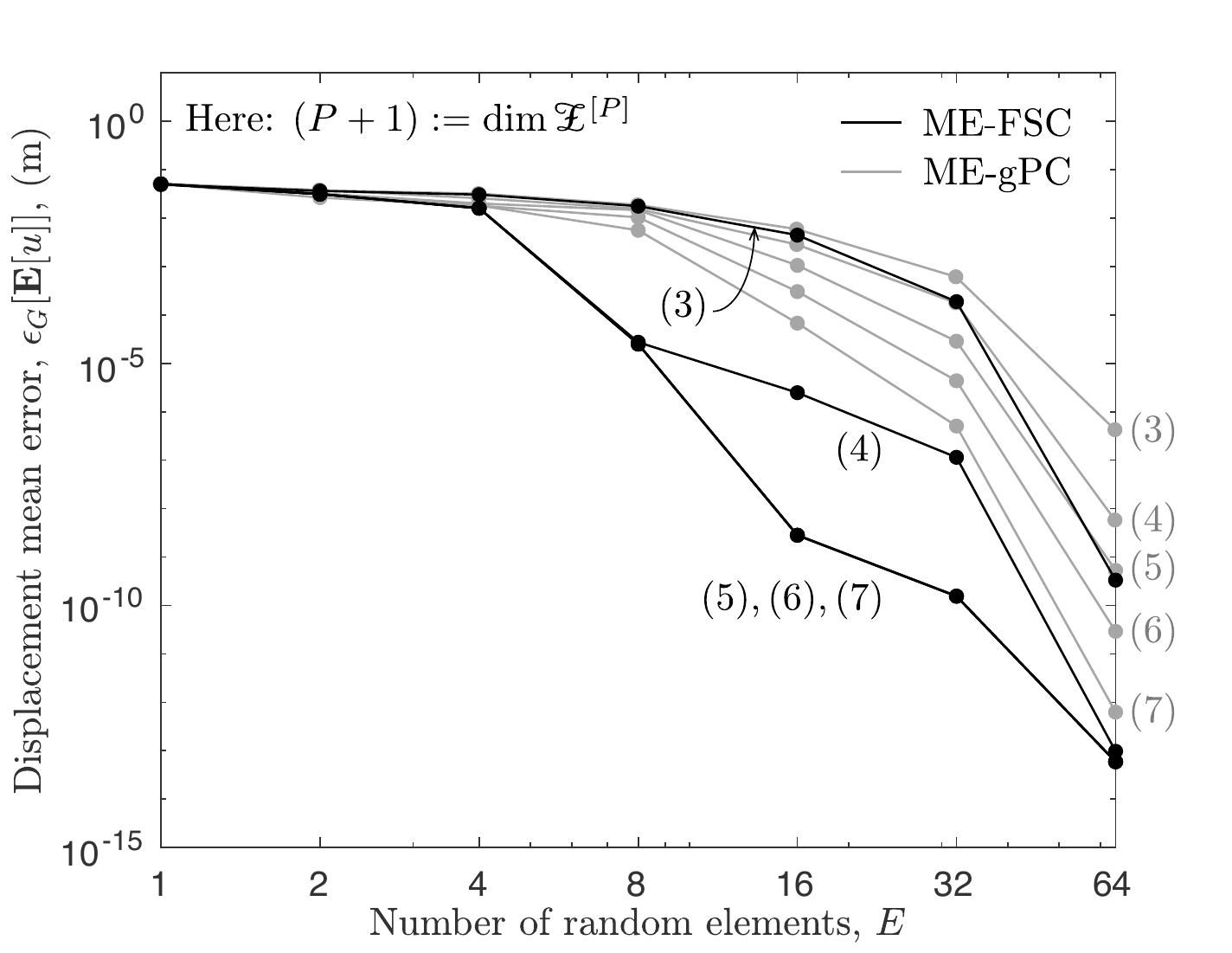}
\caption{Mean error for $\mu\sim\mathrm{Gamma}$}
\label{fig3Problem1_Gamma_MEFSC_Disp_Mean_GlobalError}
\end{subfigure}\hfill
\begin{subfigure}[b]{0.495\textwidth}
\includegraphics[width=\textwidth]{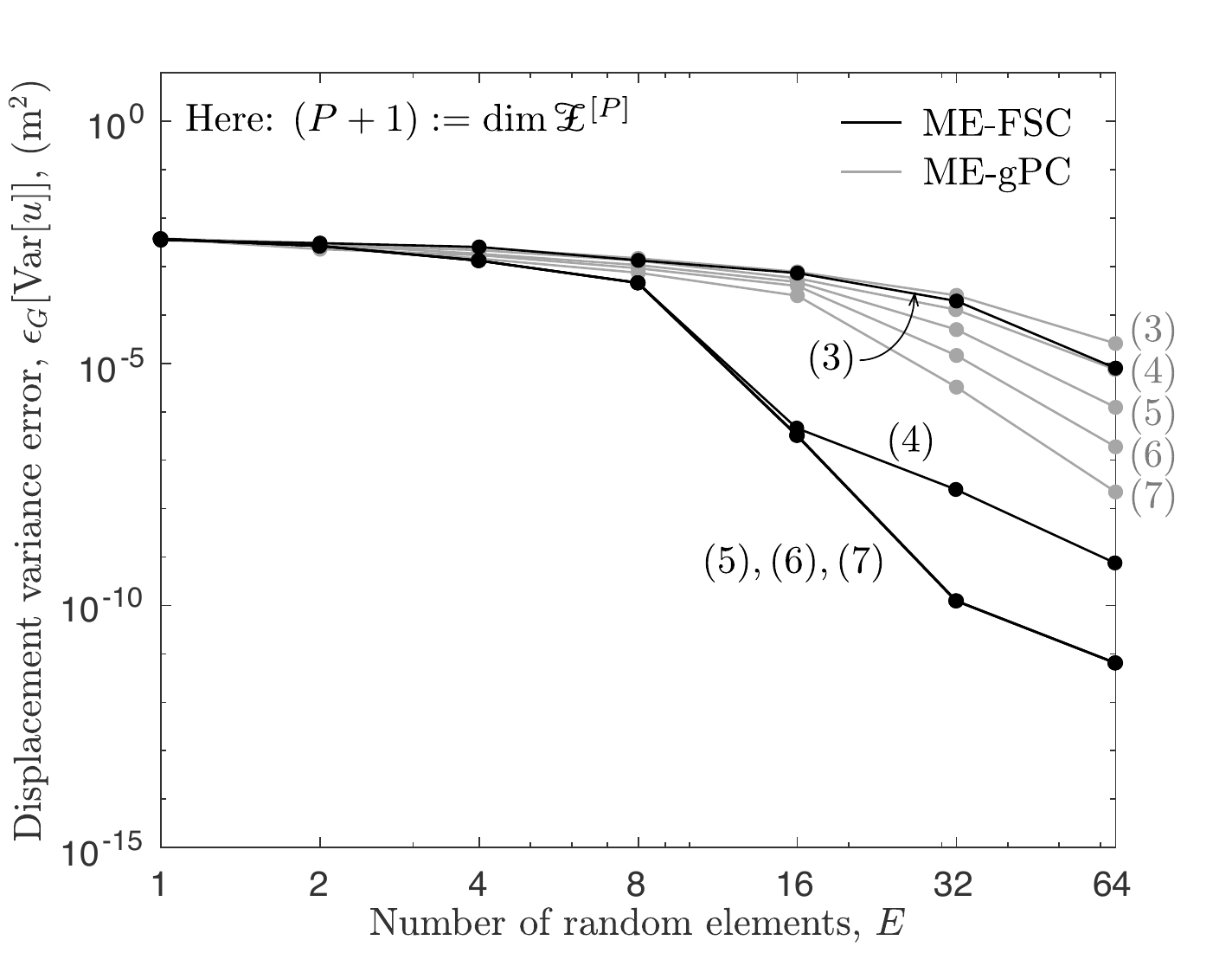}
\caption{Variance error for $\mu\sim\mathrm{Gamma}$}
\label{fig3Problem1_Gamma_MEFSC_Disp_Var_GlobalError}
\end{subfigure}
\caption{\emph{Problem 1} --- Global error of $\mathbf{E}[u]$ and $\mathrm{Var}[u]$ for different $(h,p)$-discretization levels of RFS}
\label{fig3Problem1_MEFSC_Disp_GlobalError}
\end{figure}

Fig.~\ref{fig3Problem1_MEFSC_Disp_GlobalError} depicts the convergence of global errors as a function of the number of elements and the number of basis vectors used.
Included in this figure are the cases where the probability measure is uniform, beta or gamma.
Overall, algebraic convergence can be attained if the number of elements is increased.
In particular, when the number of basis vectors is 5, the convergence is observed to be much steeper than when it is, say, 3 or 4.
However, when the number of basis vectors is greater than 5, no improvement in the accuracy of the results can be achieved.
This same outcome occurs consistently with the three distributions, which means that when a Legendre quadrature rule with 10 points per element and 5 basis vectors are used in the simulations, the maximum accuracy allowed by the machine and the ME-FSC method is ultimately reached.
Plots for the convergence of global errors using ME-gPC (with no adaptation) are also provided in the same figure.
From these, it is observed that the convergence rate is slower than ME-FSC and that for the case of uniform and beta distributions (distributions with compact support), the global errors are capable of reaching the smallest value achievable by the machine.
This, however, is not the case for the gamma distribution (a distribution with non-compact support), from which it can be observed that the variance error obtained with ME-FSC is more than 3 orders of magnitude more accurate.
This shows the great benefit of using ME-FSC versus ME-gPC to solve stochastic problems that have strong nonlinear dependencies over the probability space.

% Problem 1 [3/3]:

\subsection[Problem 2 (linear system)]{Problem 2: A linear system governed by a 3rd-order stochastic ODE}

We next consider the problem of a linear mechanical system governed by a third-order stochastic ODE.
The governing differential equation for this system is defined by
\begin{equation*}
\partial_t^3u+\tfrac{1}{2}\partial_t^2u+k\,\partial_tu+u=0,
\end{equation*}
where $k:\Xi\to\mathbb{R}$ is a stochastic mechanical parameter given by $k(\xi)=\xi$, and $u:\mathfrak{T}\times\Xi\to\mathbb{R}$ denotes the displacement of the system with $\partial_tu,\partial_t^2u,\partial_t^3u$ representing the velocity, acceleration and jerk of the system, respectively.
The initial conditions of the system are deterministic and given by: $u(0,\cdot\,)\equiv1$ m, $\partial_tu(0,\cdot\,)\equiv-1$ m/s, and $\partial_t^2{u}(0,\cdot\,)\equiv2$ m/s$^2$.

As with Problem 1, three different probability distributions are explored for $\xi$.
Namely, a \emph{uniform distribution} defined by $\mathrm{Uniform}\sim\xi\in\Xi=[a,b]$, a \emph{beta distribution} defined by $\mathrm{Beta}(\alpha,\beta)\sim\xi\in\Xi=[a,b]$, and a \emph{normal distribution} defined by $\mathrm{Normal}(\mu,\sigma^2)\sim\xi\in\Xi=\mathbb{R}$, whence $(a,b)=(2,3)$ N/m, $(\alpha,\beta)=(2,5)$ and $(\mu,\sigma)=(2.5,0.125)$ N/m.
However, since the normal distribution has non-compact support, the simulations are run with \emph{truncated normal distribution}, $\mathrm{Normal}(\mu,\sigma^2,c,d)\sim\xi\in\Xi=[c,d]$, where $(c,d)=(1.4,3.6)$ N/m.
The exact solution, nonetheless, is obtained with the normal distribution.

\begin{figure}
\centering
\begin{subfigure}[b]{0.495\textwidth}
\includegraphics[width=\textwidth]{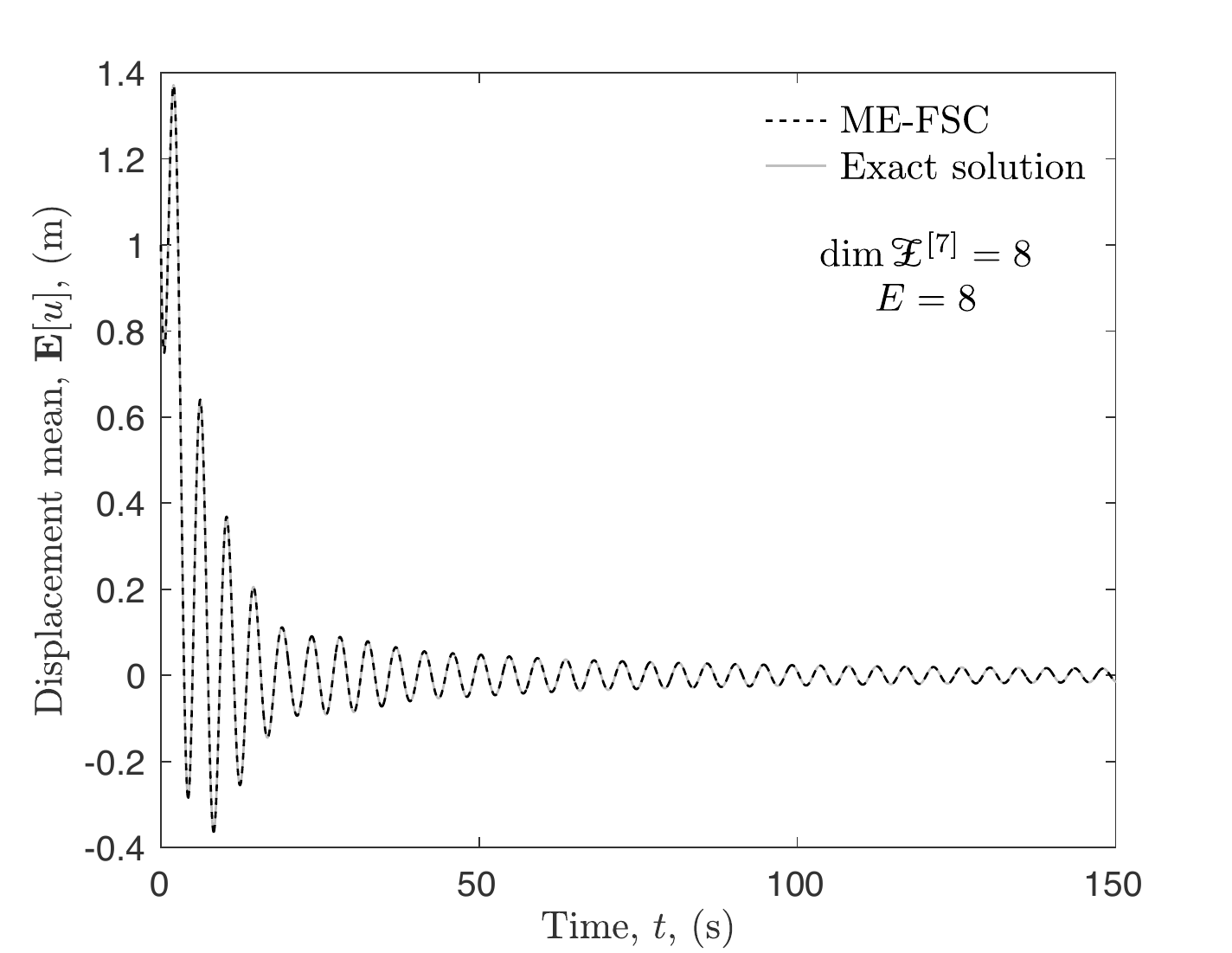}
\caption{Mean}
\label{fig3Problem2_Uniform_MEFSC_8BV_8E_Disp_Mean}
\end{subfigure}\hfill
\begin{subfigure}[b]{0.495\textwidth}
\includegraphics[width=\textwidth]{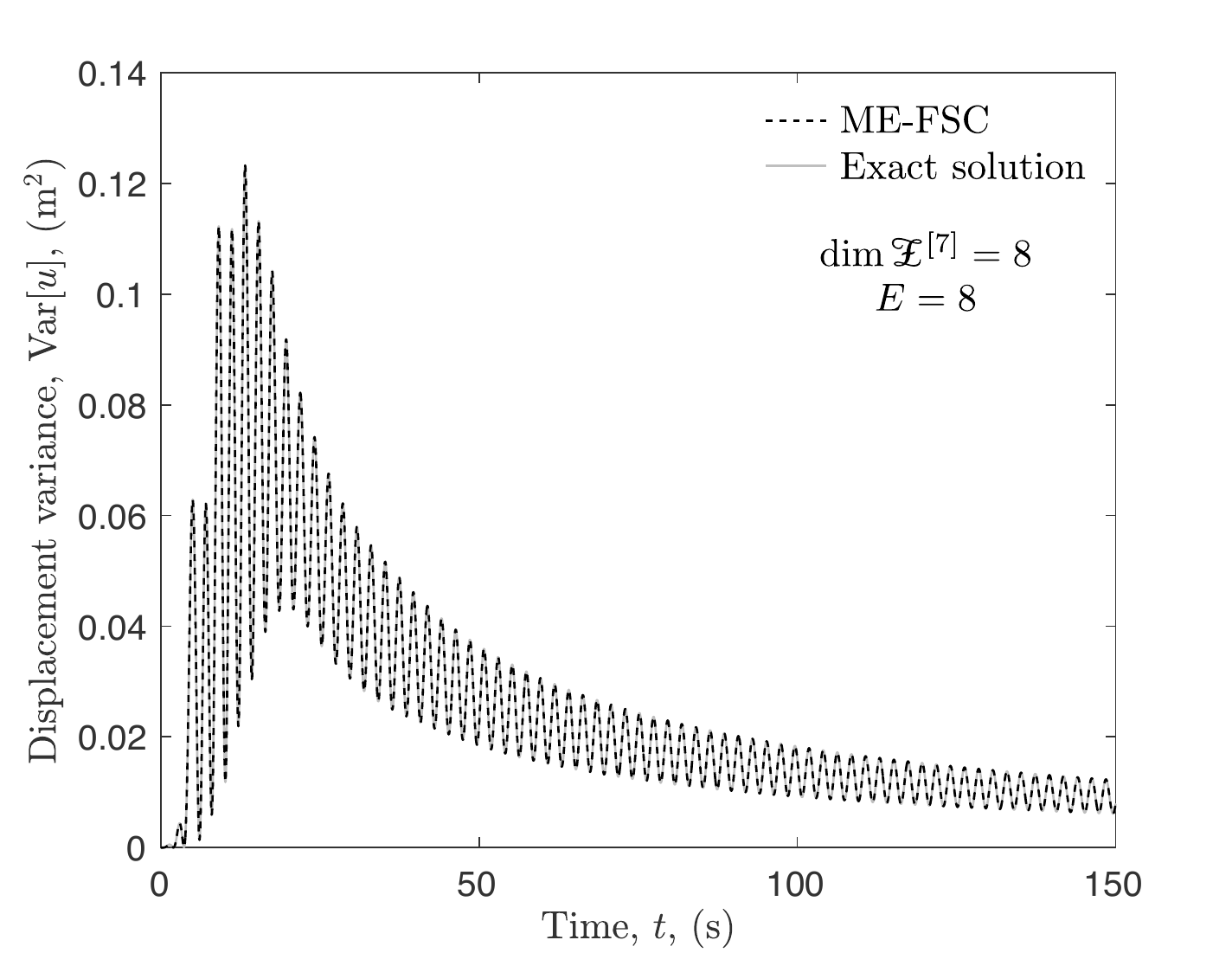}
\caption{Variance}
\label{fig3Problem2_Uniform_MEFSC_8BV_8E_Disp_Var}
\end{subfigure}
\caption{\emph{Problem 2} --- Evolution of $\mathbf{E}[u]$ and $\mathrm{Var}[u]$ for the case when the $(h,p)$-discretization level of RFS is $(P,E)=(7,8)$ and $\mu\sim\mathrm{Uniform}$}
\label{fig3Problem2_Uniform_MEFSC_8BV_8E_Disp}
\end{figure}

Fig.~\ref{fig3Problem2_Uniform_MEFSC_8BV_8E_Disp} shows the evolution of the mean and variance of the system's displacement for the particular case where the probability measure is uniform and 8 elements and 8 basis vectors are used.
Like before, the results obtained with ME-FSC are again indistinguishable from the exact response.

% Problem 2 [1/3]:

\begin{figure}
\centering
\begin{subfigure}[b]{0.495\textwidth}
\includegraphics[width=\textwidth]{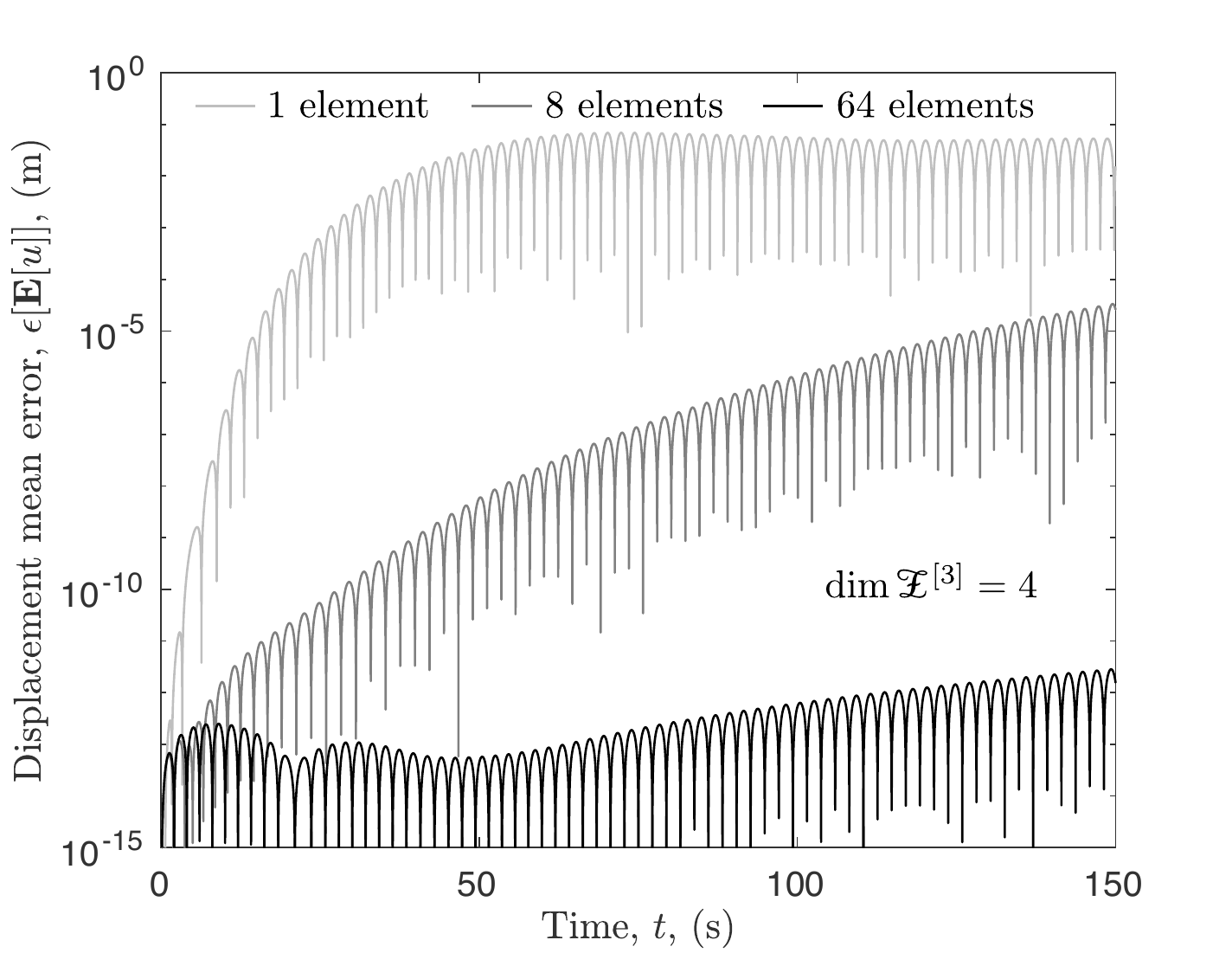}
\caption{Mean error for $\mathscr{Z}^{[3]}$}
\label{fig3Problem2_Uniform_MEFSC_4BV_Disp_Mean_LocalError}
\end{subfigure}\hfill
\begin{subfigure}[b]{0.495\textwidth}
\includegraphics[width=\textwidth]{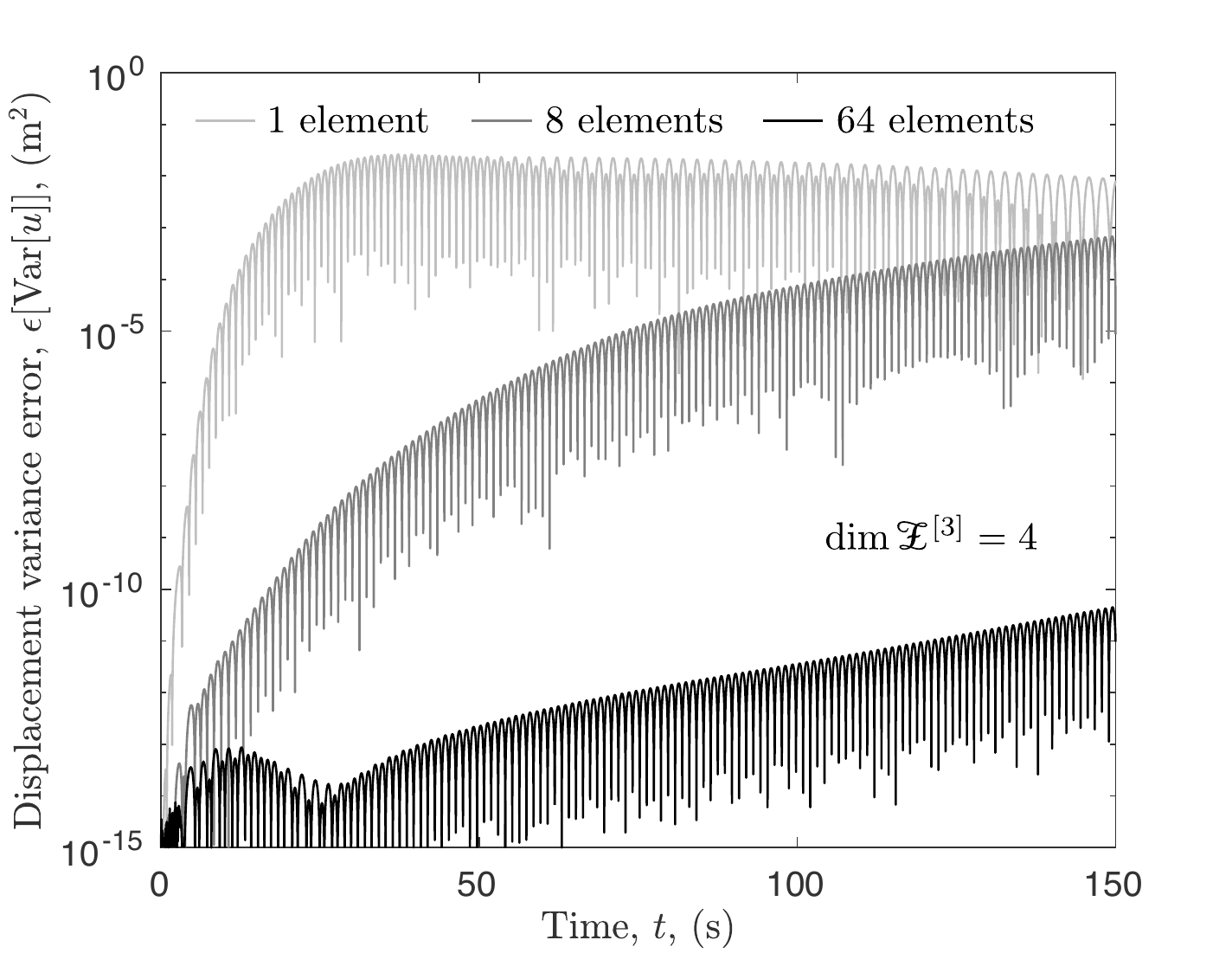}
\caption{Variance error for $\mathscr{Z}^{[3]}$}
\label{fig3Problem2_Uniform_MEFSC_4BV_Disp_Var_LocalError}
\end{subfigure}\quad
\begin{subfigure}[b]{0.495\textwidth}
\includegraphics[width=\textwidth]{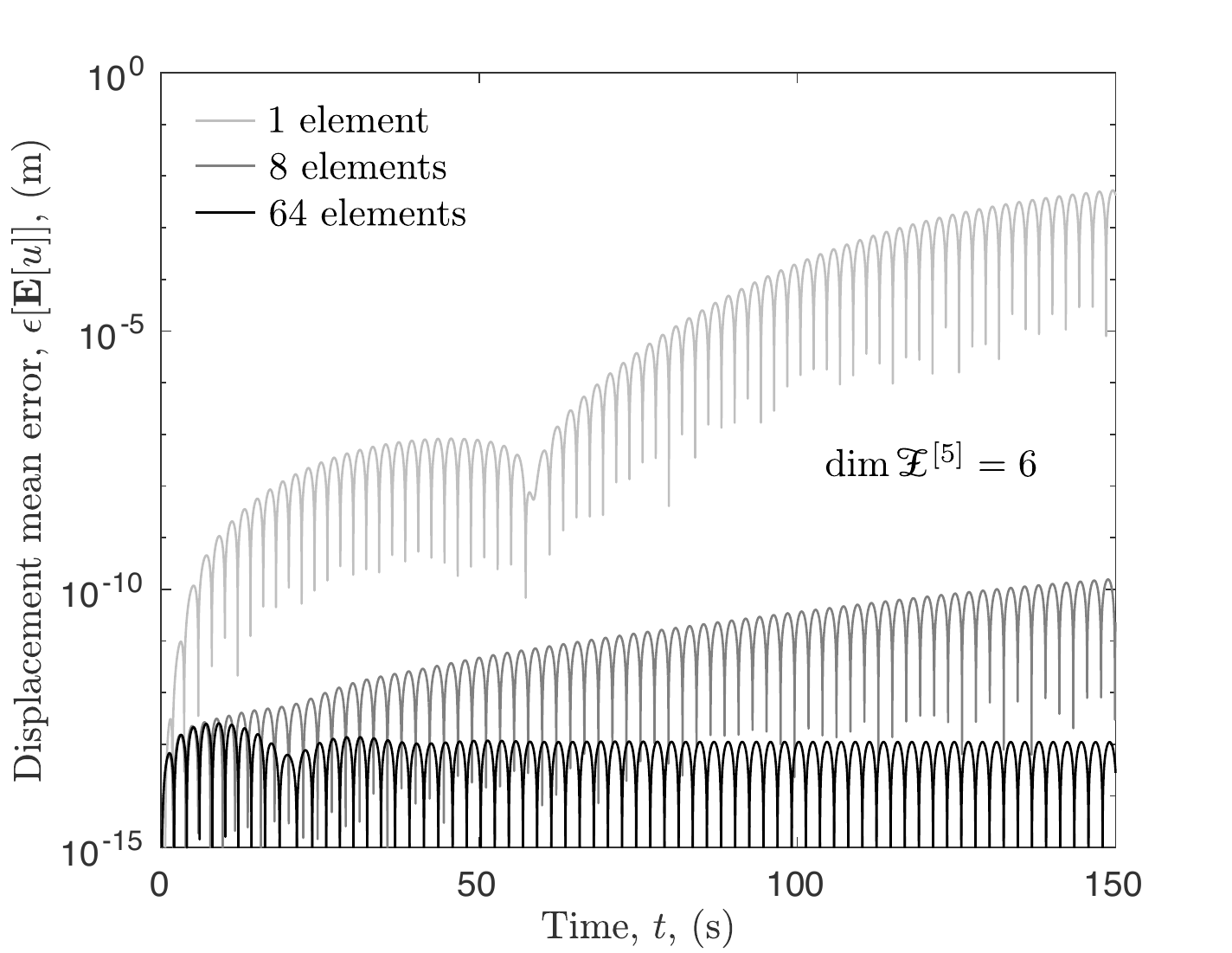}
\caption{Mean error for $\mathscr{Z}^{[5]}$}
\label{fig3Problem2_Uniform_MEFSC_6BV_Disp_Mean_LocalError}
\end{subfigure}\hfill
\begin{subfigure}[b]{0.495\textwidth}
\includegraphics[width=\textwidth]{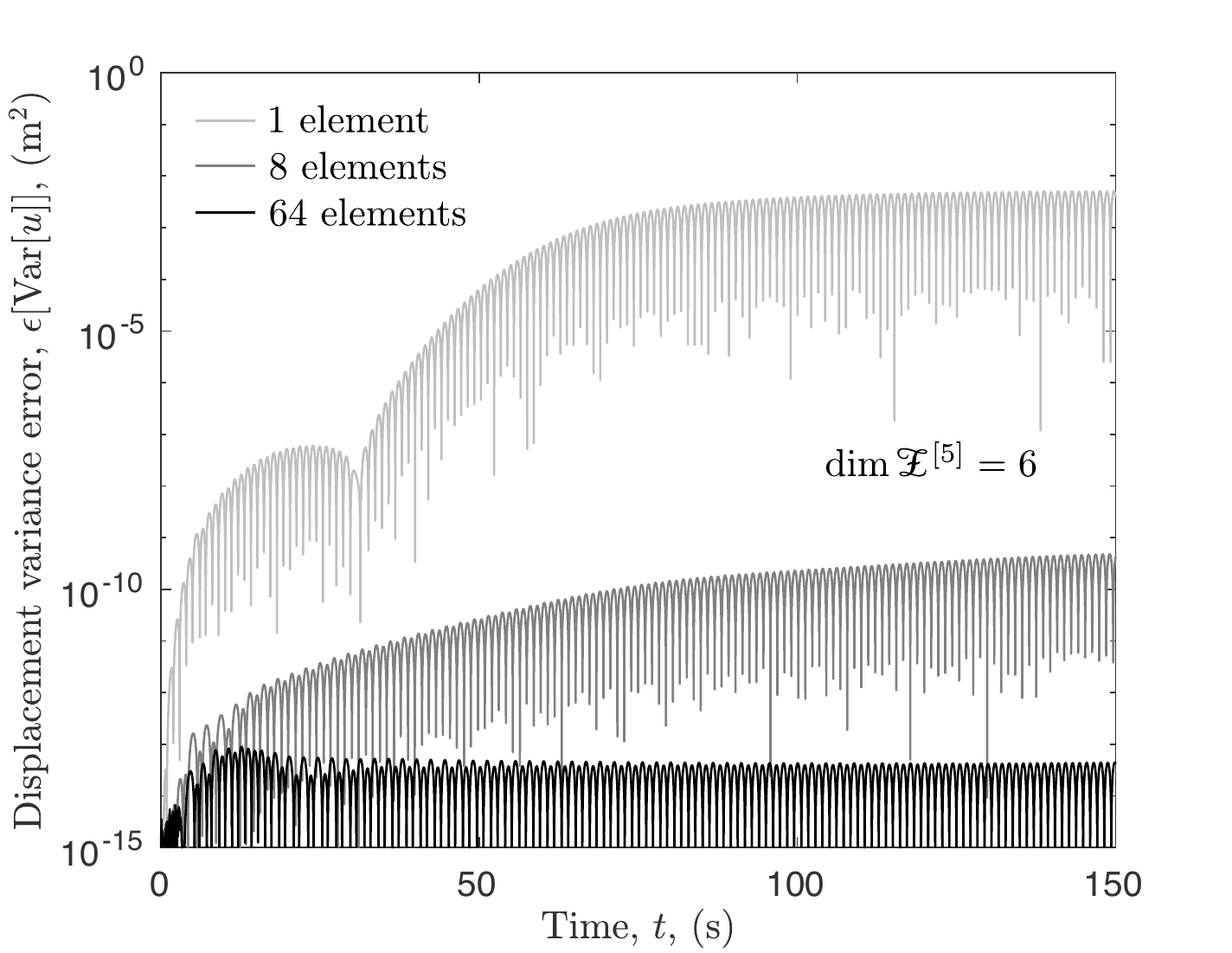}
\caption{Variance error for $\mathscr{Z}^{[5]}$}
\label{fig3Problem2_Uniform_MEFSC_6BV_Disp_Var_LocalError}
\end{subfigure}\quad
\begin{subfigure}[b]{0.495\textwidth}
\includegraphics[width=\textwidth]{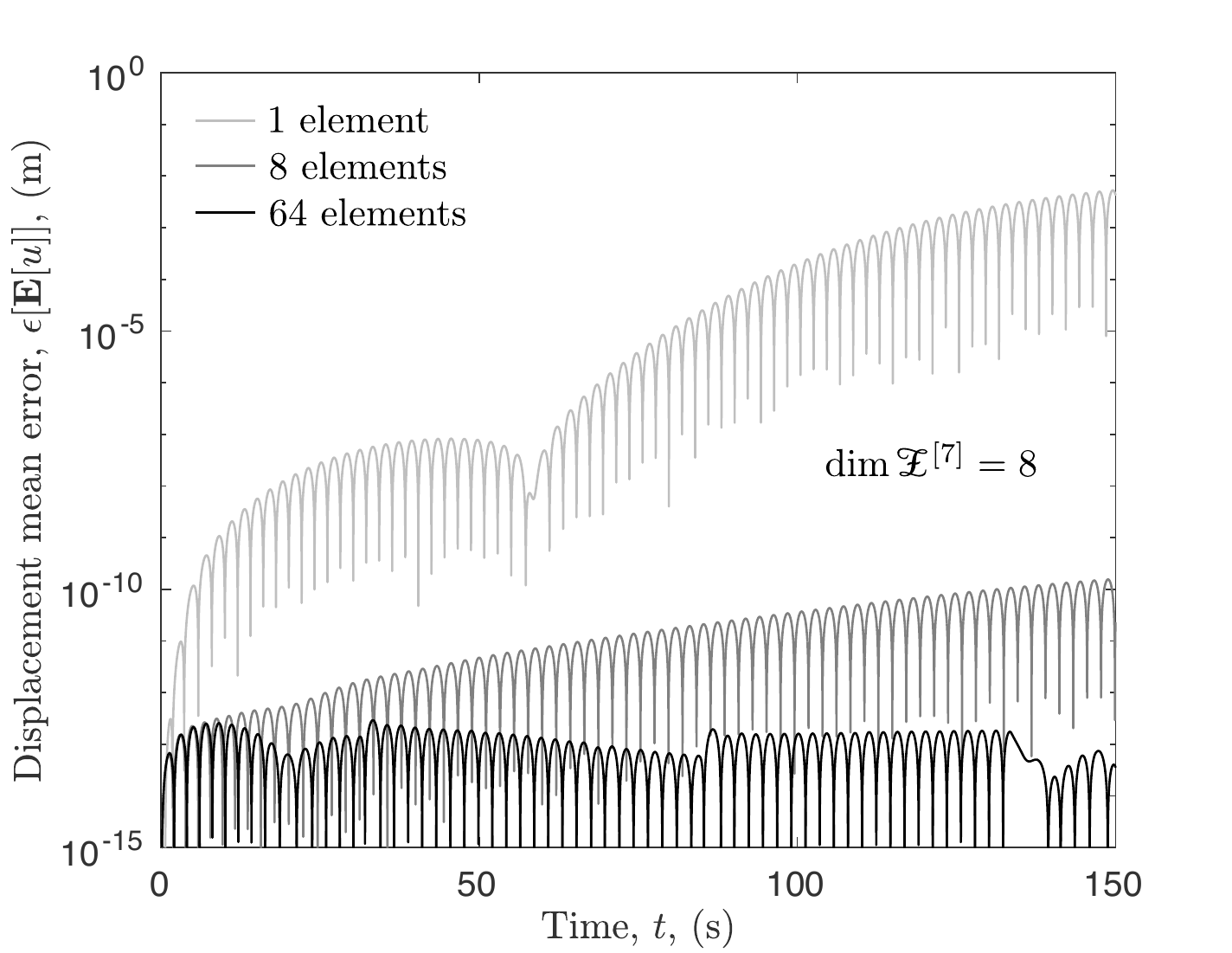}
\caption{Mean error for $\mathscr{Z}^{[7]}$}
\label{fig3Problem2_Uniform_MEFSC_8BV_Disp_Mean_LocalError}
\end{subfigure}\hfill
\begin{subfigure}[b]{0.495\textwidth}
\includegraphics[width=\textwidth]{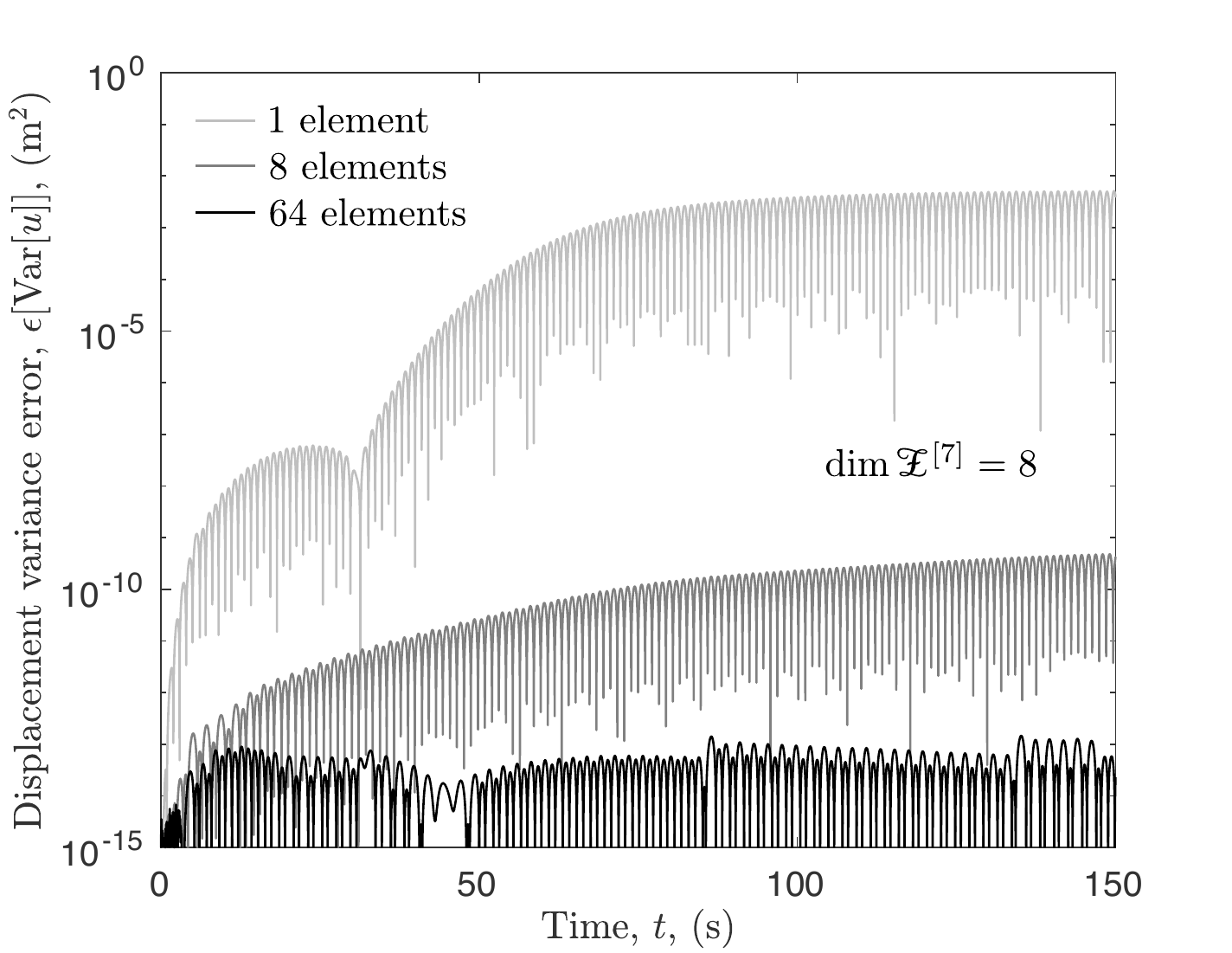}
\caption{Variance error for $\mathscr{Z}^{[7]}$}
\label{fig3Problem2_Uniform_MEFSC_8BV_Disp_Var_LocalError}
\end{subfigure}
\caption{\emph{Problem 2} --- Local error evolution of $\mathbf{E}[u]$ and $\mathrm{Var}[u]$ for different $(h,p)$-discretization levels of RFS and for $\mu\sim\mathrm{Uniform}$}
\label{fig3Problem2_Uniform_MEFSC_Disp_LocalError}
\end{figure}

In Fig.~\ref{fig3Problem2_Uniform_MEFSC_Disp_LocalError} we present the local errors in mean and variance of the system's displacement using different numbers of elements and basis vectors.
For brevity, the results are only presented for $\mu\sim\mathrm{Uniform}$.
Once again, it is observed that the accuracy of the results improves as the number of elements increases, but it necessarily does not as the number of basis vectors increases.
This is exemplified in Figs.~\ref{fig3Problem2_Uniform_MEFSC_6BV_Disp_Mean_LocalError} to \ref{fig3Problem2_Uniform_MEFSC_8BV_Disp_Var_LocalError} from where it is deduced that the results do not improve if the number of basis vectors is increased from 6 to 8.
This once again is due to the quadrature rule used to compute the inner products (i.e.~10 Legendre quadrature points per element) and the limited precision of the machine.

% Problem 2 [2/3]:

\begin{figure}
\centering
\begin{subfigure}[b]{0.495\textwidth}
\includegraphics[width=\textwidth]{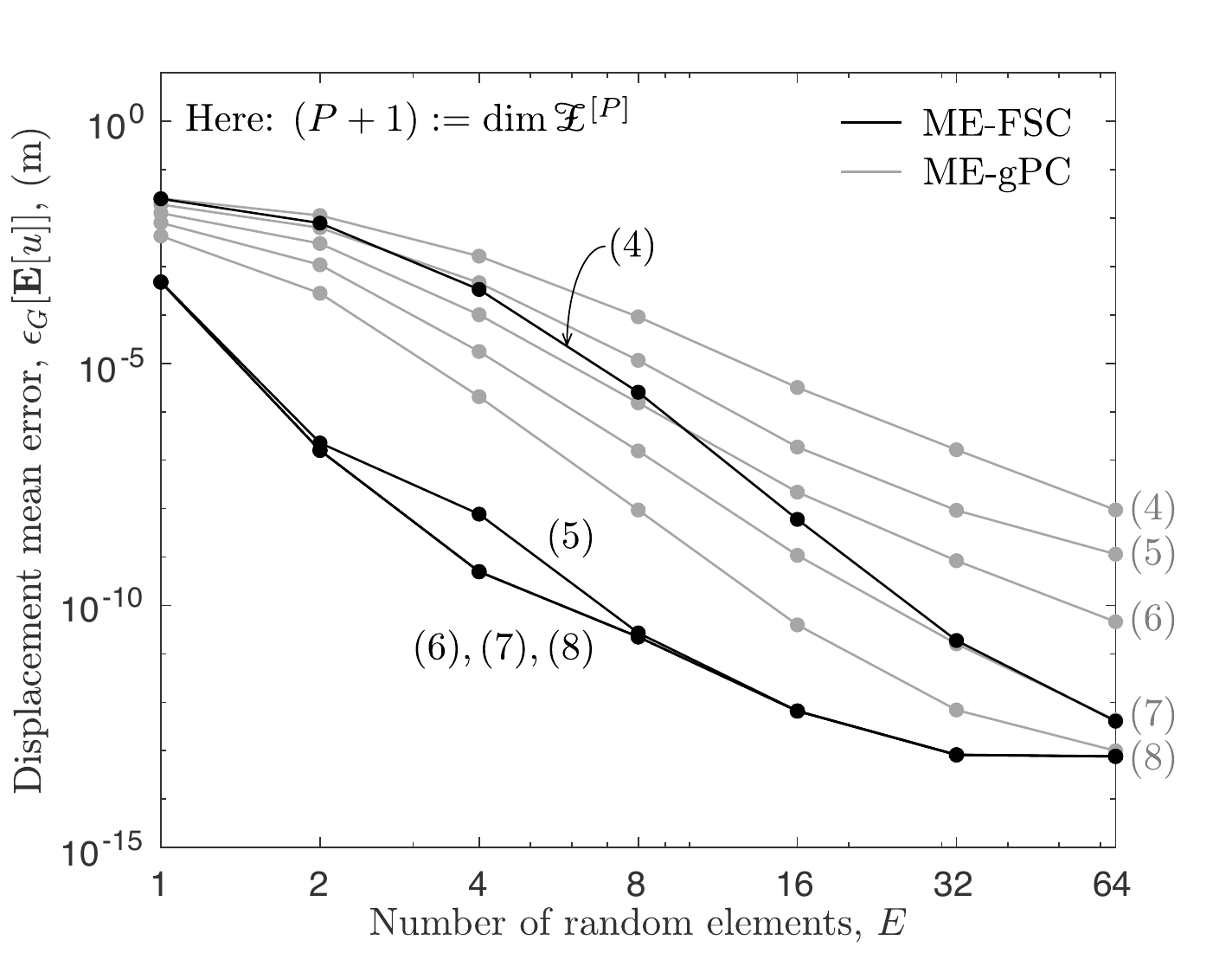}
\caption{Mean error for $\mu\sim\mathrm{Uniform}$}
\label{fig3Problem2_Uniform_MEFSC_Disp_Mean_GlobalError}
\end{subfigure}\hfill
\begin{subfigure}[b]{0.495\textwidth}
\includegraphics[width=\textwidth]{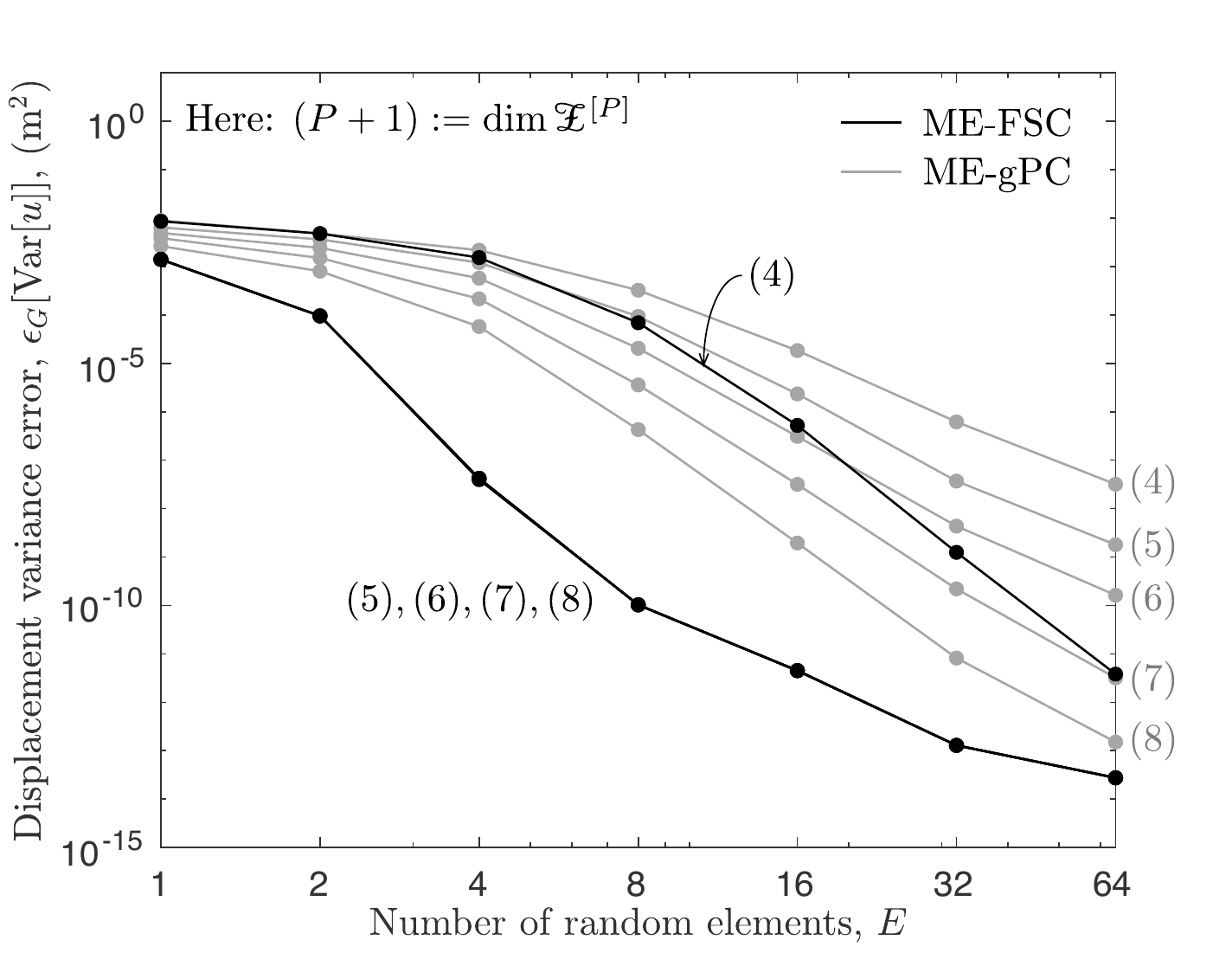}
\caption{Variance error for $\mu\sim\mathrm{Uniform}$}
\label{fig3Problem2_Uniform_MEFSC_Disp_Var_GlobalError}
\end{subfigure}\quad
\begin{subfigure}[b]{0.495\textwidth}
\includegraphics[width=\textwidth]{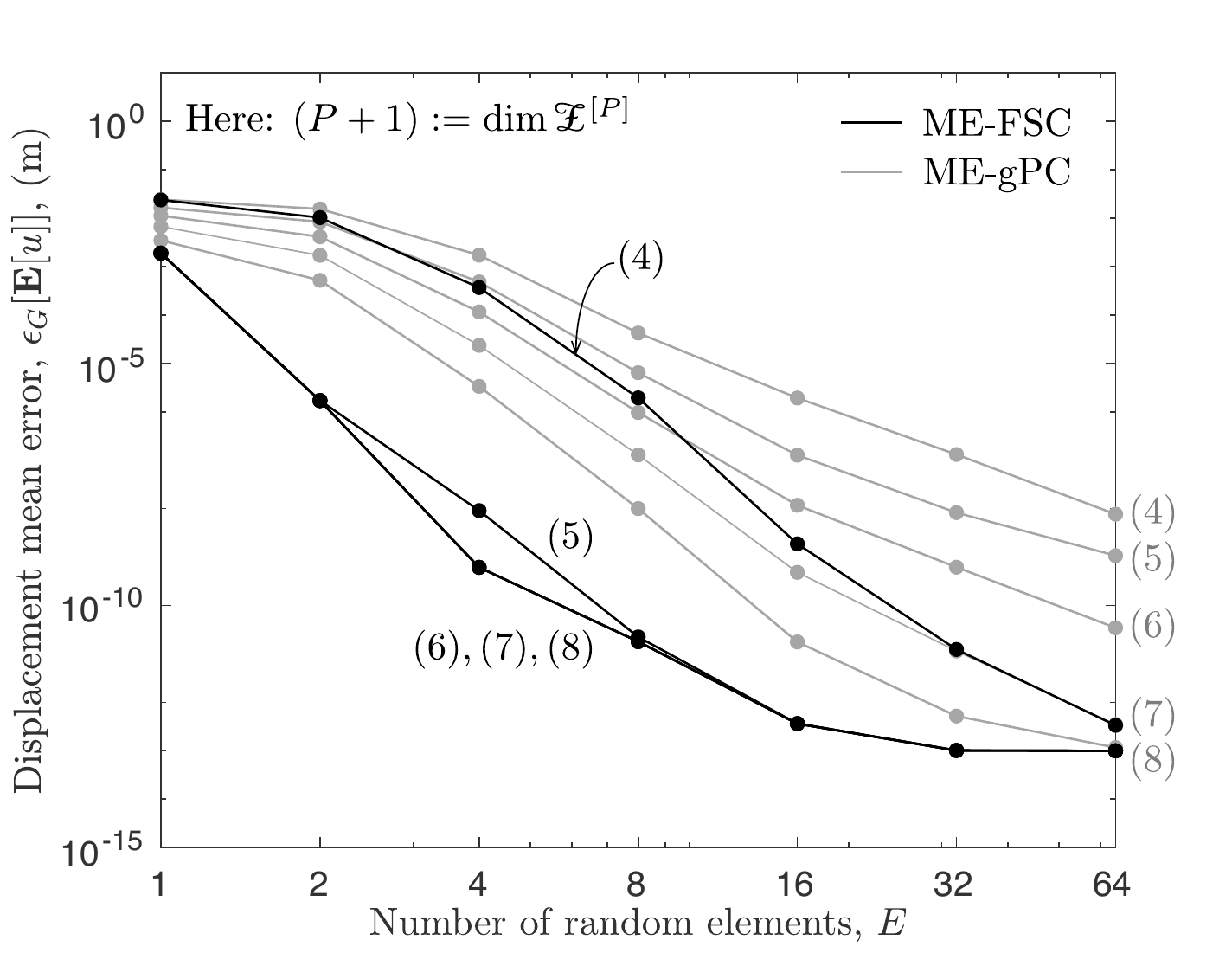}
\caption{Mean error for $\mu\sim\mathrm{Beta}$}
\label{fig3Problem2_Beta_MEFSC_Disp_Mean_GlobalError}
\end{subfigure}\hfill
\begin{subfigure}[b]{0.495\textwidth}
\includegraphics[width=\textwidth]{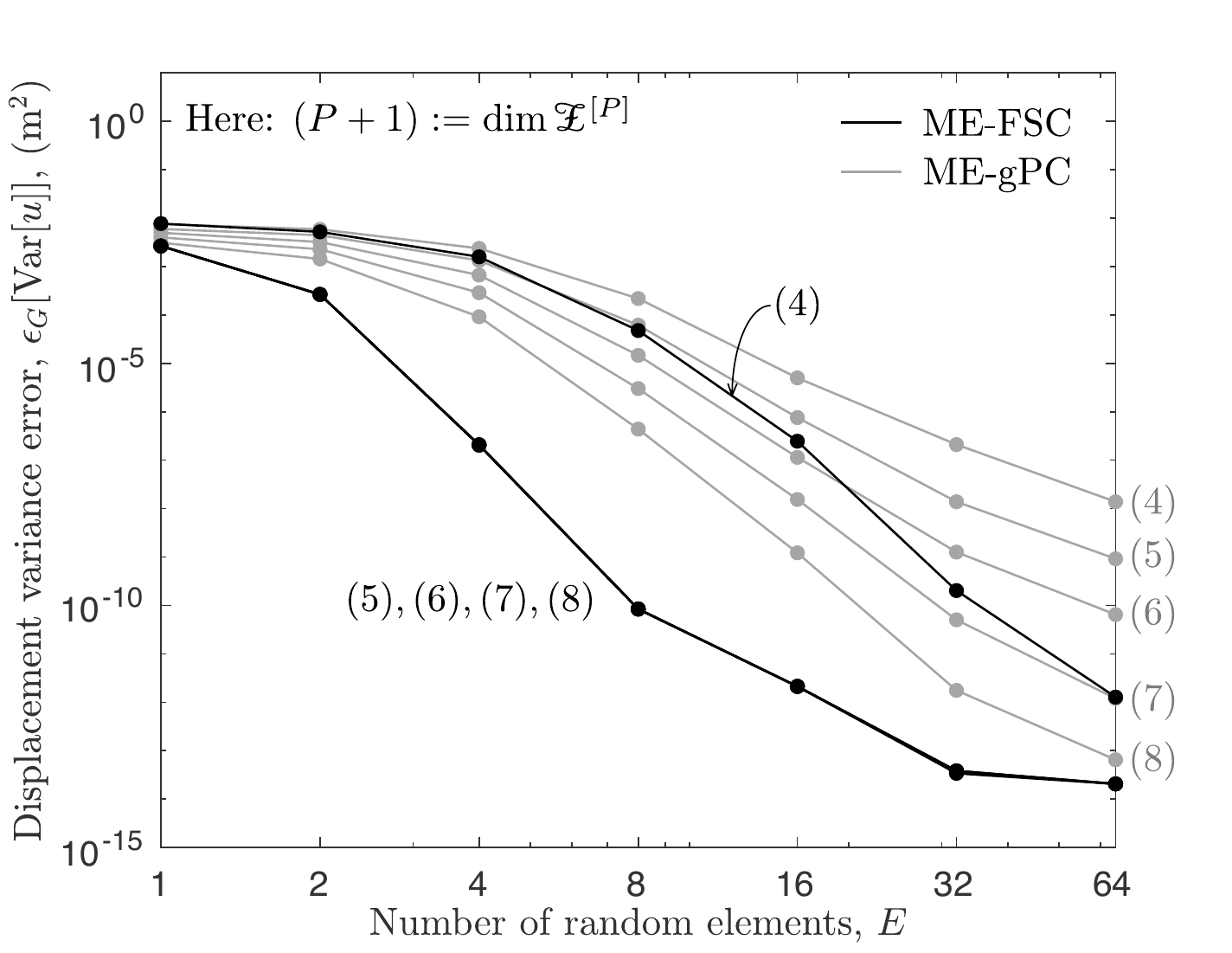}
\caption{Variance error for $\mu\sim\mathrm{Beta}$}
\label{fig3Problem2_Beta_MEFSC_Disp_Var_GlobalError}
\end{subfigure}\quad
\begin{subfigure}[b]{0.495\textwidth}
\includegraphics[width=\textwidth]{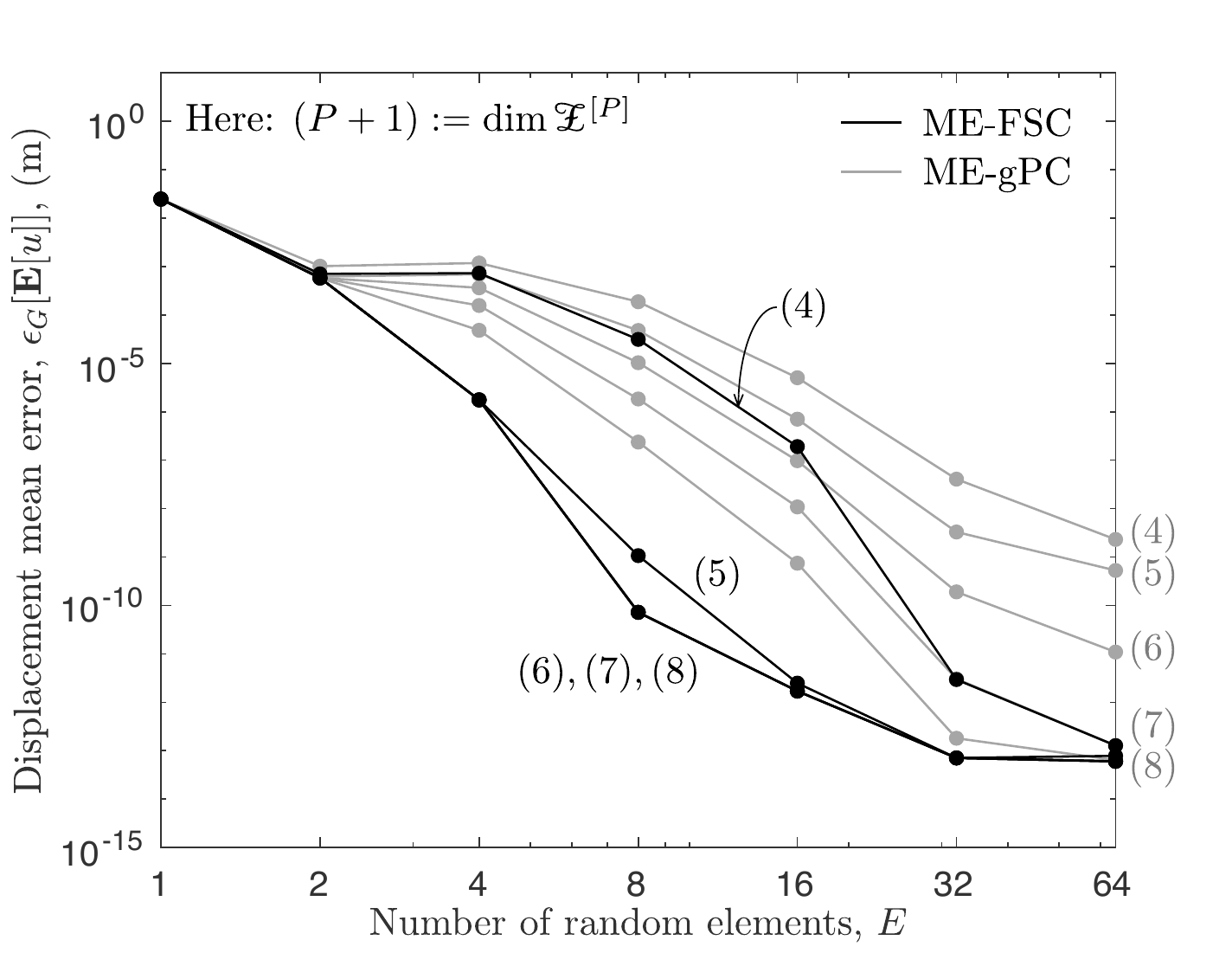}
\caption{Mean error for $\mu\sim\mathrm{Normal}$}
\label{fig3Problem2_Normal_MEFSC_Disp_Mean_GlobalError}
\end{subfigure}\hfill
\begin{subfigure}[b]{0.495\textwidth}
\includegraphics[width=\textwidth]{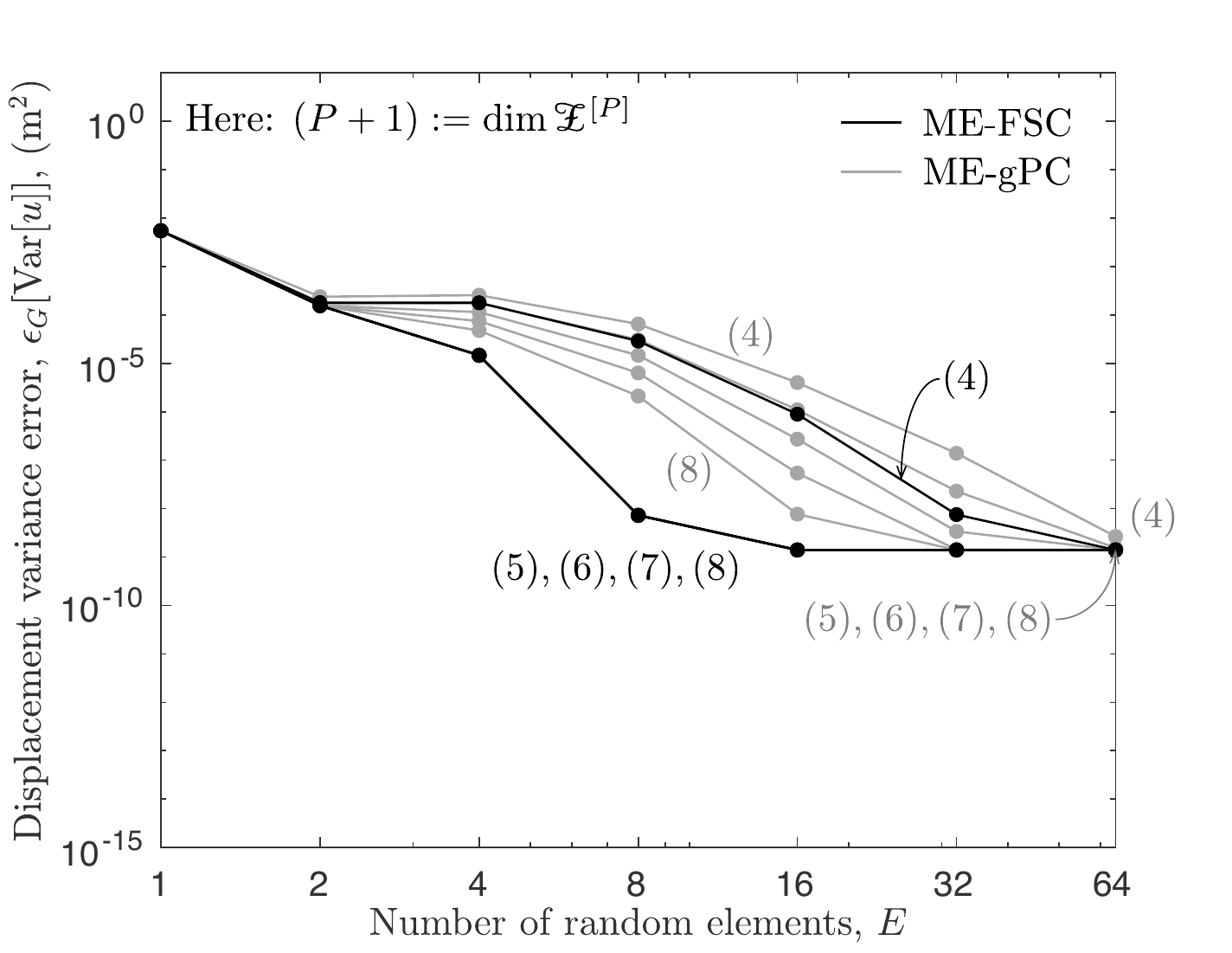}
\caption{Variance error for $\mu\sim\mathrm{Normal}$}
\label{fig3Problem2_Normal_MEFSC_Disp_Var_GlobalError}
\end{subfigure}
\caption{\emph{Problem 2} --- Global error of $\mathbf{E}[u]$ and $\mathrm{Var}[u]$ for different $(h,p)$-discretization levels of RFS}
\label{fig3Problem2_MEFSC_Disp_GlobalError}
\end{figure}

Finally, Fig.~\ref{fig3Problem2_MEFSC_Disp_GlobalError} depicts the global errors in mean and variance of the system's displacement.
The results are presented as a function of the number of elements and the number of basis vectors used for each of the probability distributions explored for $\xi$.
Again, algebraic convergence is overall attainable for ME-FSC if the number of elements is increased.
However, as shown, using more than 4 basis vectors in the simulations does not help improve the overall accuracy of the results, except of course when 4 elements and 5 basis vectors are used for the mean response.
This figure therefore suggests that in some situations, it might be beneficial to refine the partition of the random domain (instead of making the random function space bigger) to obtain more accurate results.
As for the simulations with ME-gPC, the convergence trend shown in Problem 1 is also showcased in this problem through the same type of plot.
For those distributions with compact support (i.e.~uniform and beta distributions), the global errors can eventually reach the smallest value attainable by the machine.
In fact, this same conclusion applies to the distribution with non-compact support (i.e.~normal distribution), which, as observed, is also capable of achieving global errors as small as those produced by ME-FSC.

% Problem 2 [3/3]:

\subsection[Problem 3 (nonlinear system)]{Problem 3: A nonlinear system governed by a 2nd-order stochastic ODE (the Van-der-Pol oscillator)}

In this section, we investigate the nonlinear behavior of a single-degree-of-freedom system with Van-der-Pol damping.
The governing differential equation for this system is given by
\begin{equation*}
m\ddot{u}-(1-\rho u^2)c\dot{u}+ku=0,
\end{equation*}
where $m=100$ kg is the mass of the system, $\rho=150$ m$^{-2}$ is the contributing factor to the nonlinearity of the system, $c:\Xi\to\mathbb{R}^+$ is a coefficient representing the strength of the damping in the system, and $k:\Xi\to\mathbb{R}^+$ is the stiffness of the system.
The functions $c$ and $k$ are assumed to be given by $c(\xi)=\xi^1$ and $k(\xi)=\xi^2$.
The initial conditions of the system are: $u(0,\cdot\,)\equiv 0.20$ m and $\dot{u}(0,\cdot\,)\equiv 0.30$ m/s.
Notice here that $u:\mathfrak{T}\times\Xi\to\mathbb{R}$ denotes the displacement of the system, and that $\dot{u}:=\partial_t u$ and $\ddot{u}:=\partial_t^2 u$ represent, respectively, the velocity and acceleration of the system.

For concreteness, we take $\xi^1$ to be \emph{uniformly distributed} in $\bar{\Xi}_1=[150,450]$ kg/s, and $\xi^2$ to be \emph{beta distributed} with parameters $(\alpha,\beta)=(2,5)$ in $\bar{\Xi}_2=[340,460]$ N/m.
Hence, the random space is two-dimensional and defined by $\Xi=\bar{\Xi}_1\times\bar{\Xi}_2$ with $\mu\sim\mathrm{Uniform}\otimes\mathrm{Beta}$.

\begin{figure}
\centering
\begin{subfigure}[b]{0.495\textwidth}
\includegraphics[width=\textwidth]{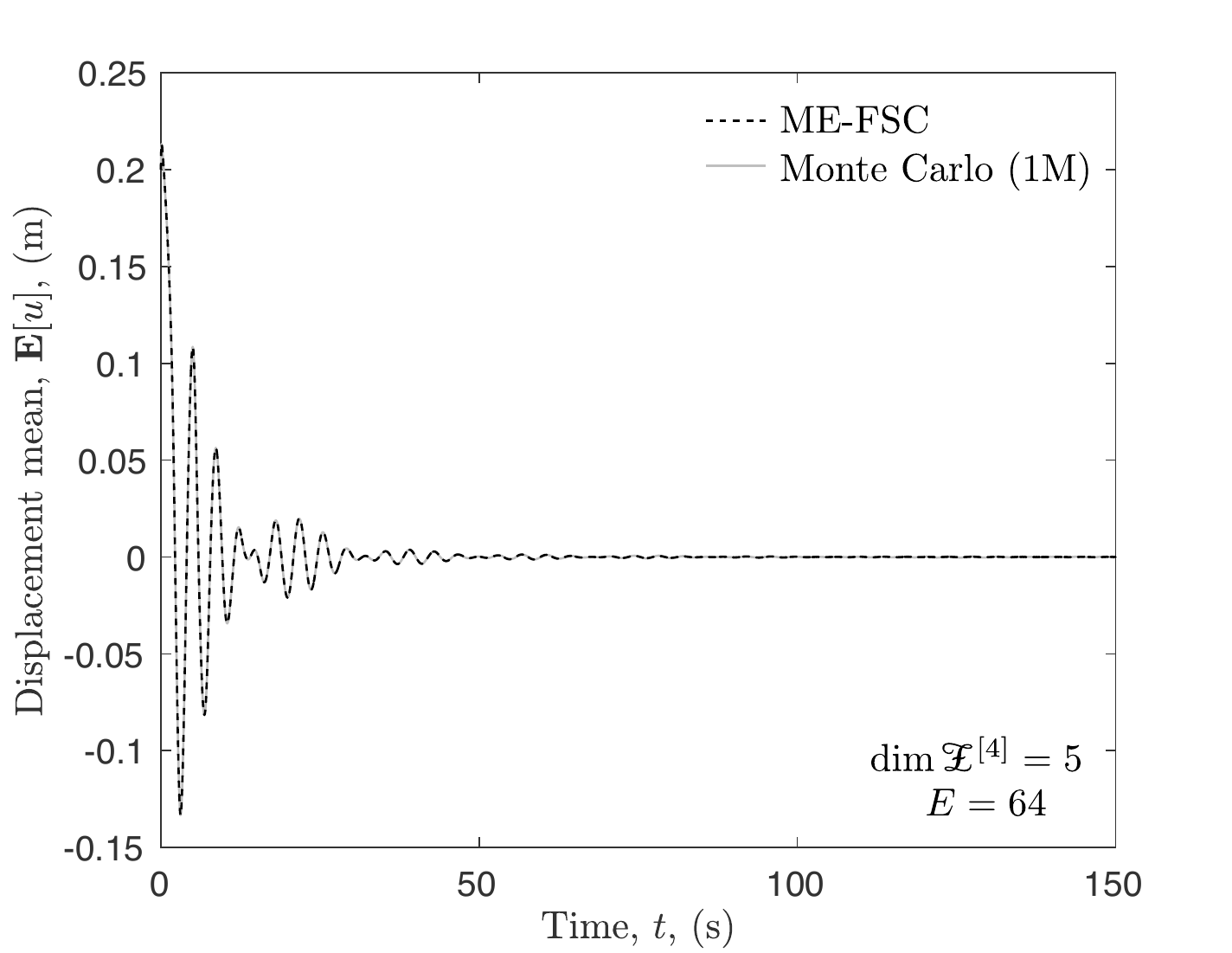}
\caption{Mean}
\label{fig3Problem3_UniformBeta_MEFSC_5BV_64E_Disp_Mean}
\end{subfigure}\hfill
\begin{subfigure}[b]{0.495\textwidth}
\includegraphics[width=\textwidth]{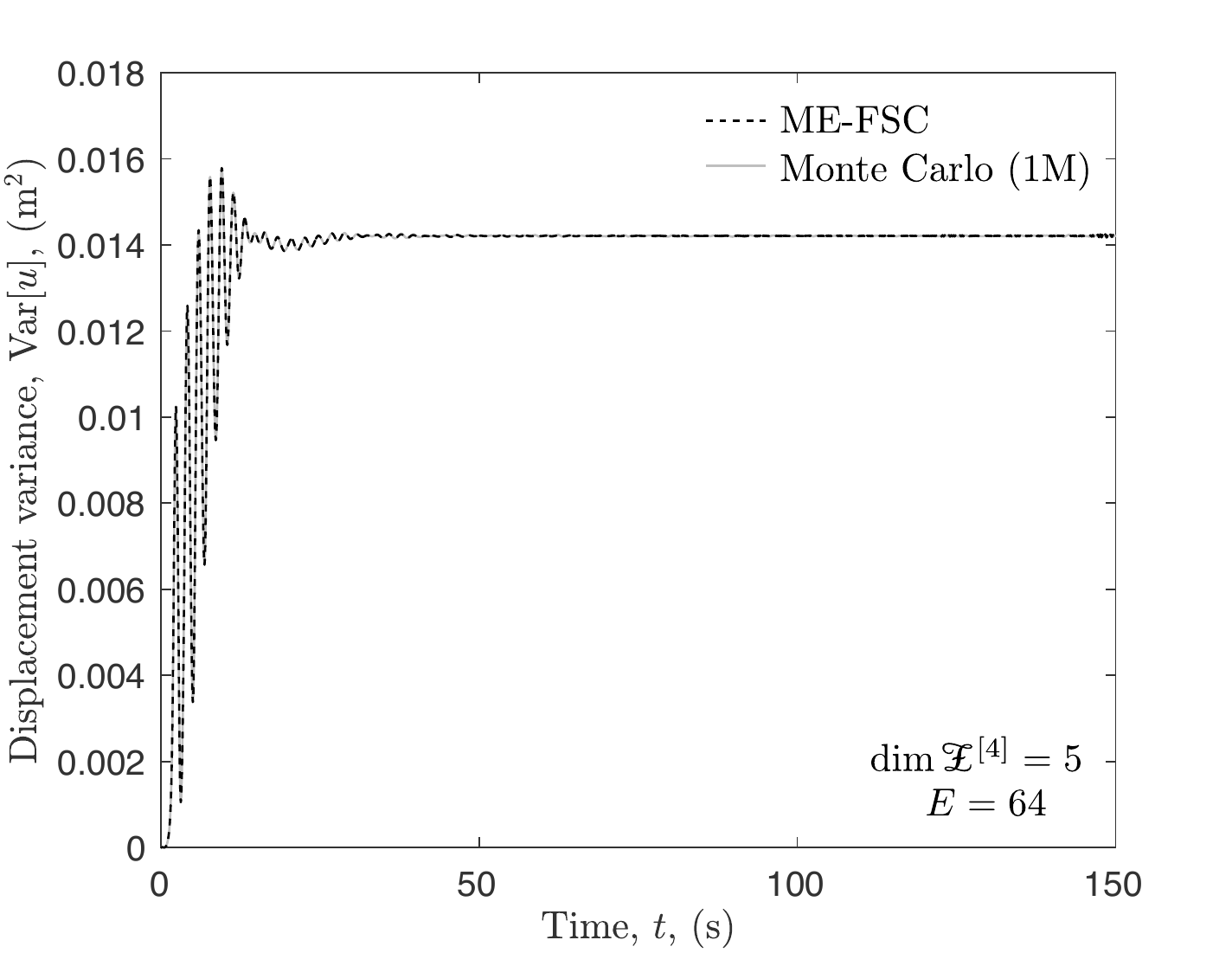}
\caption{Variance}
\label{fig3Problem3_UniformBeta_MEFSC_5BV_64E_Disp_Var}
\end{subfigure}
\caption{\emph{Problem 3} --- Evolution of $\mathbf{E}[u]$ and $\mathrm{Var}[u]$ for the case when the $(h,p)$-discretization level of RFS is $(P,E)=(4,64)$ and $\mu\sim\mathrm{Uniform}\otimes\mathrm{Beta}$}
\label{fig3Problem3_UniformBeta_MEFSC_5BV_64E_Disp}
\end{figure}

\begin{figure}
\centering
\begin{subfigure}[b]{0.495\textwidth}
\includegraphics[width=\textwidth]{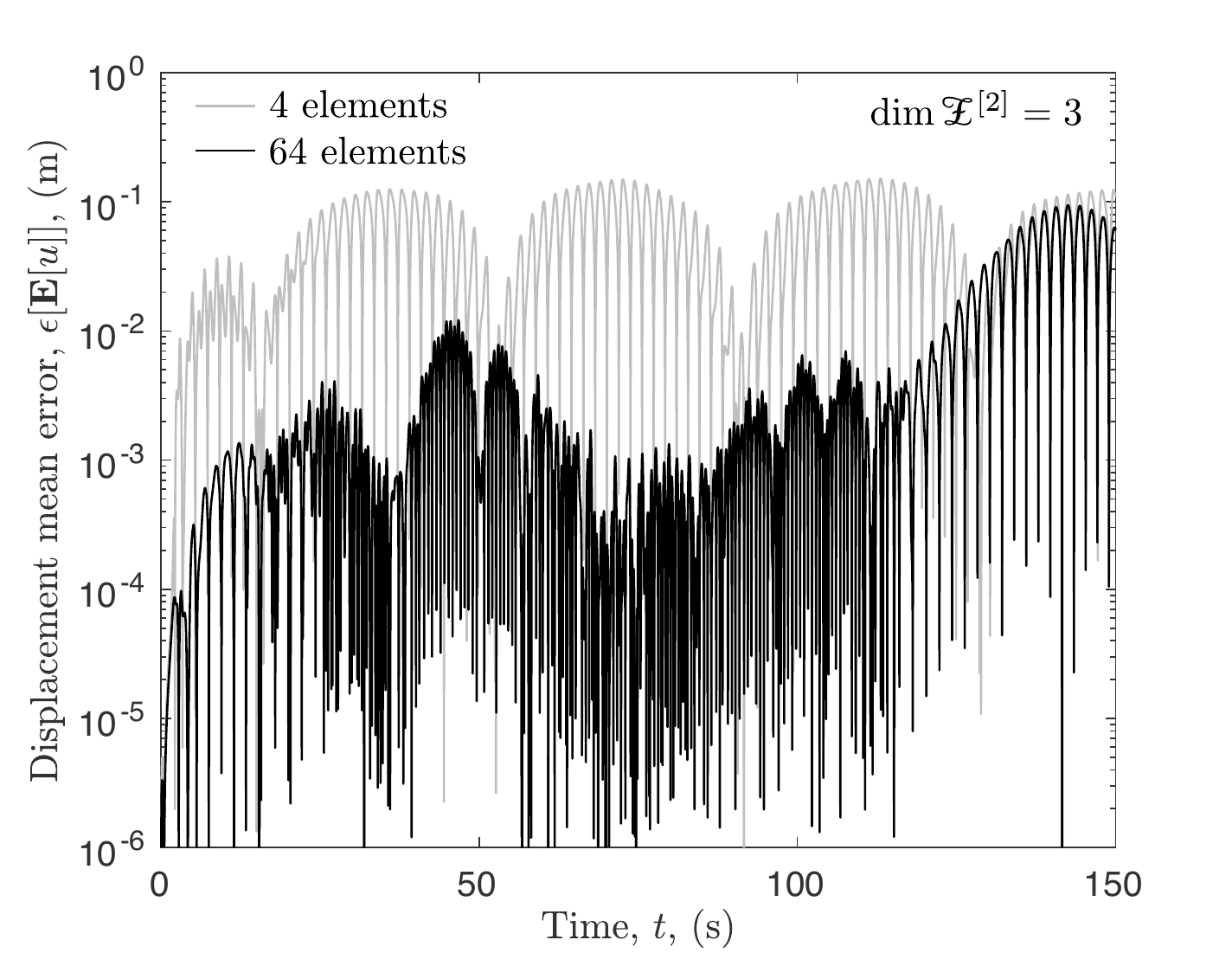}
\caption{Mean error for $\mathscr{Z}^{[2]}$}
\label{fig3Problem3_UniformBeta_MEFSC_3BV_Disp_Mean_LocalError}
\end{subfigure}\hfill
\begin{subfigure}[b]{0.495\textwidth}
\includegraphics[width=\textwidth]{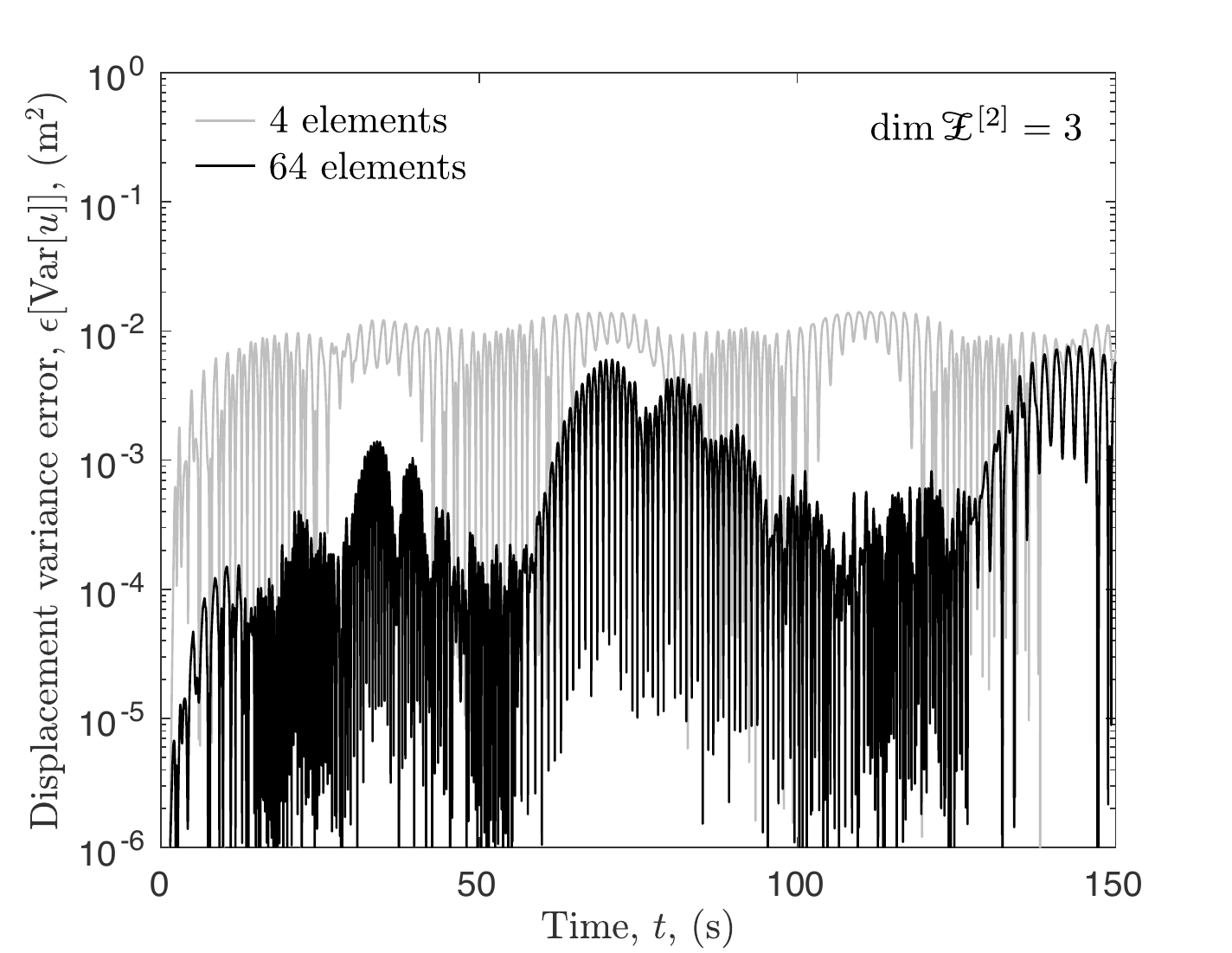}
\caption{Variance error for $\mathscr{Z}^{[2]}$}
\label{fig3Problem3_UniformBeta_MEFSC_3BV_Disp_Var_LocalError}
\end{subfigure}\quad
\begin{subfigure}[b]{0.495\textwidth}
\includegraphics[width=\textwidth]{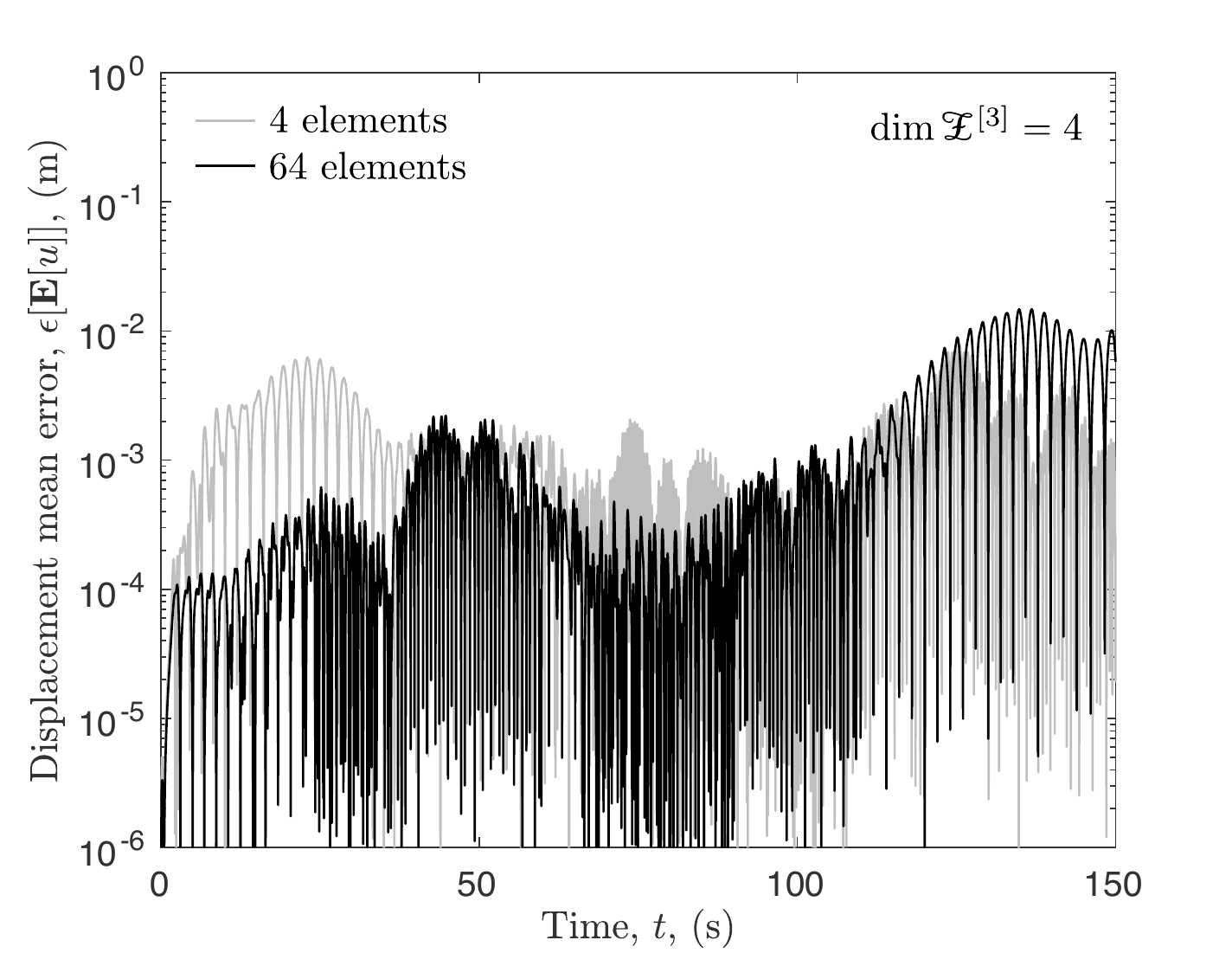}
\caption{Mean error for $\mathscr{Z}^{[3]}$}
\label{fig3Problem3_UniformBeta_MEFSC_4BV_Disp_Mean_LocalError}
\end{subfigure}\hfill
\begin{subfigure}[b]{0.495\textwidth}
\includegraphics[width=\textwidth]{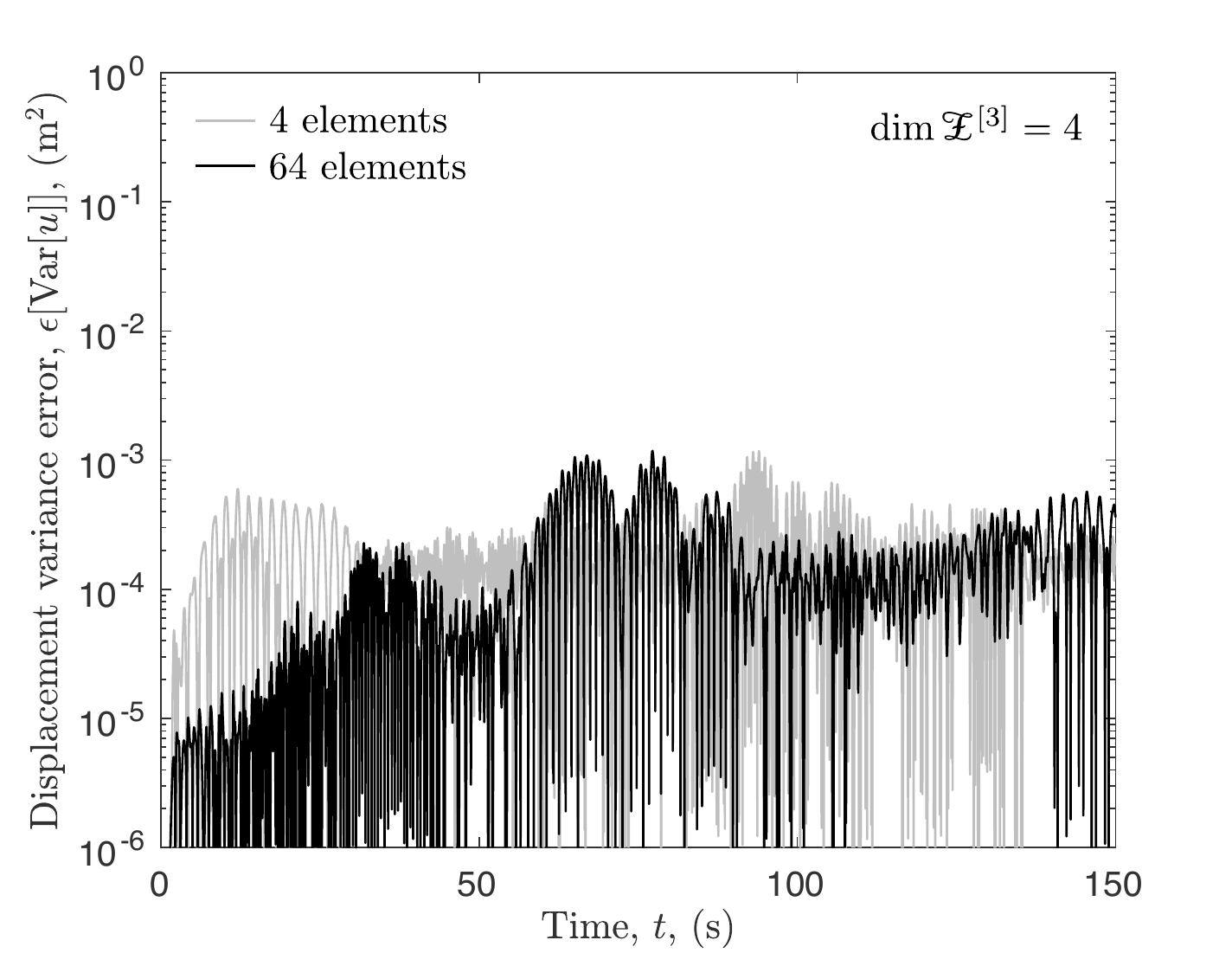}
\caption{Variance error for $\mathscr{Z}^{[3]}$}
\label{fig3Problem3_UniformBeta_MEFSC_4BV_Disp_Var_LocalError}
\end{subfigure}\quad
\begin{subfigure}[b]{0.495\textwidth}
\includegraphics[width=\textwidth]{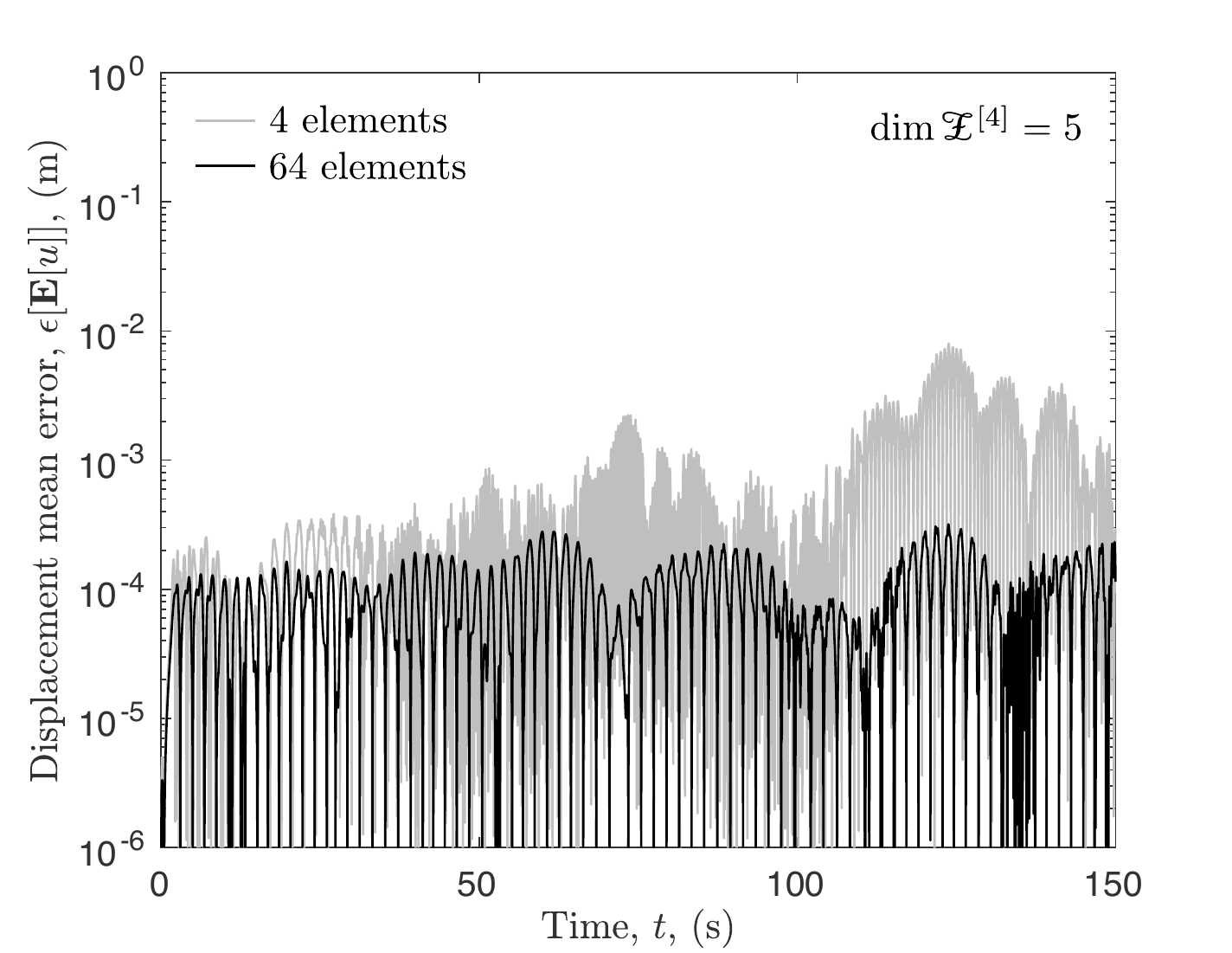}
\caption{Mean error for $\mathscr{Z}^{[4]}$}
\label{fig3Problem3_UniformBeta_MEFSC_5BV_Disp_Mean_LocalError}
\end{subfigure}\hfill
\begin{subfigure}[b]{0.495\textwidth}
\includegraphics[width=\textwidth]{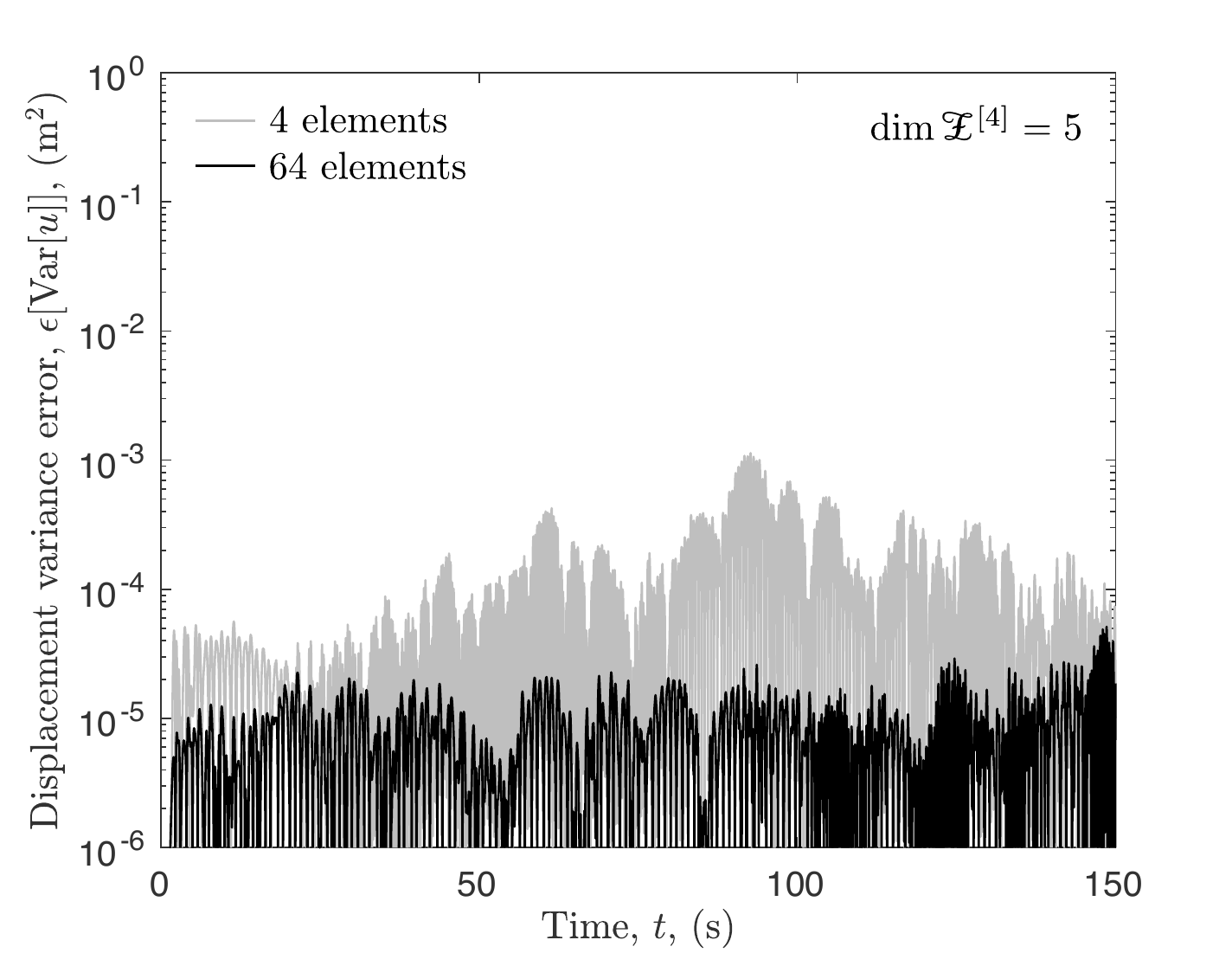}
\caption{Variance error for $\mathscr{Z}^{[4]}$}
\label{fig3Problem3_UniformBeta_MEFSC_5BV_Disp_Var_LocalError}
\end{subfigure}
\caption{\emph{Problem 3} --- Local error evolution of $\mathbf{E}[u]$ and $\mathrm{Var}[u]$ for different $(h,p)$-discretization levels of RFS and for $\mu\sim\mathrm{Uniform}\otimes\mathrm{Beta}$}
\label{fig3Problem3_UniformBeta_MEFSC_Disp_LocalError}
\end{figure}

\begin{figure}
\centering
\begin{subfigure}[b]{0.495\textwidth}
\includegraphics[width=\textwidth]{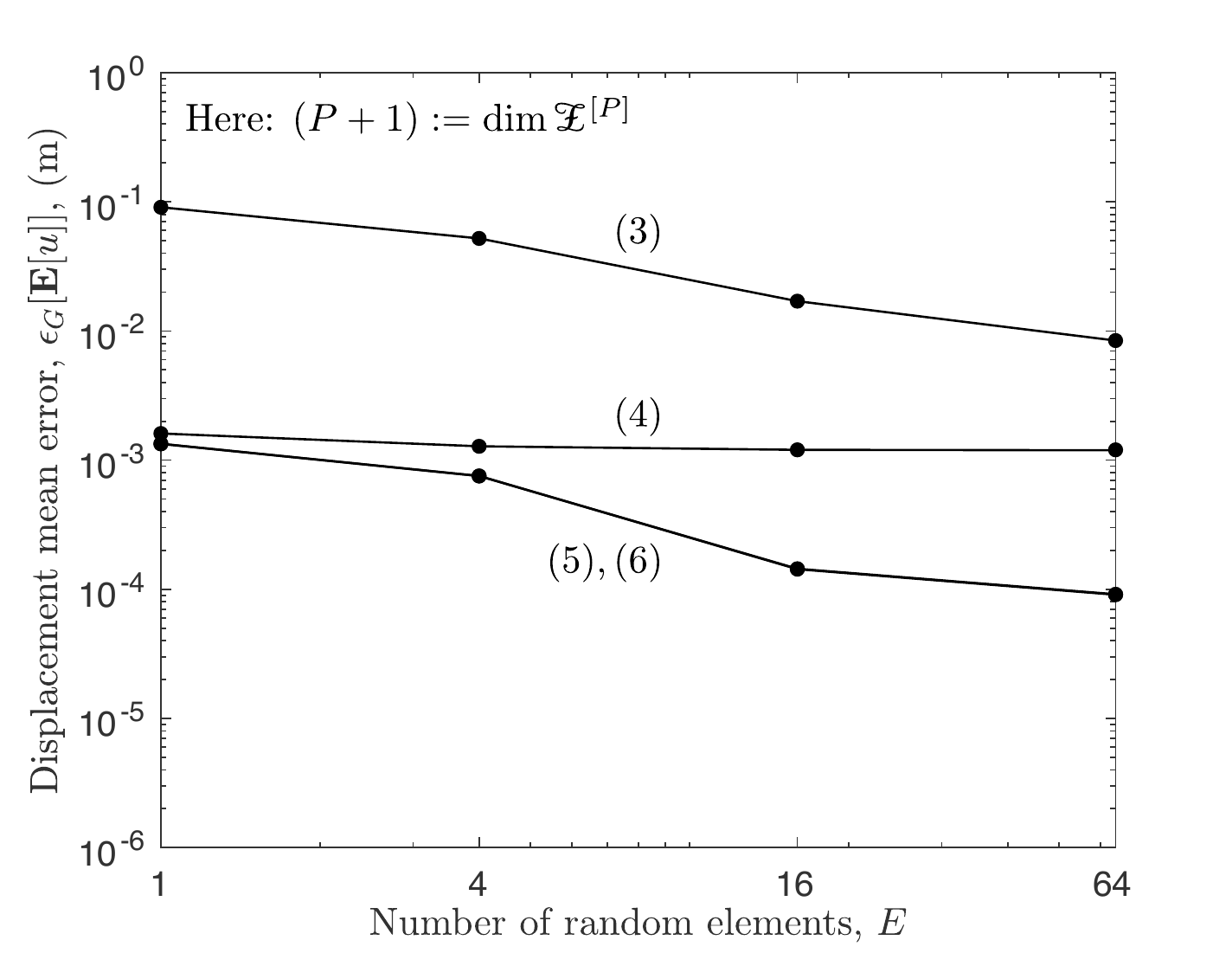}
\caption{Mean error}
\label{fig3Problem3_UniformBeta_MEFSC_Disp_Mean_GlobalError}
\end{subfigure}\hfill
\begin{subfigure}[b]{0.495\textwidth}
\includegraphics[width=\textwidth]{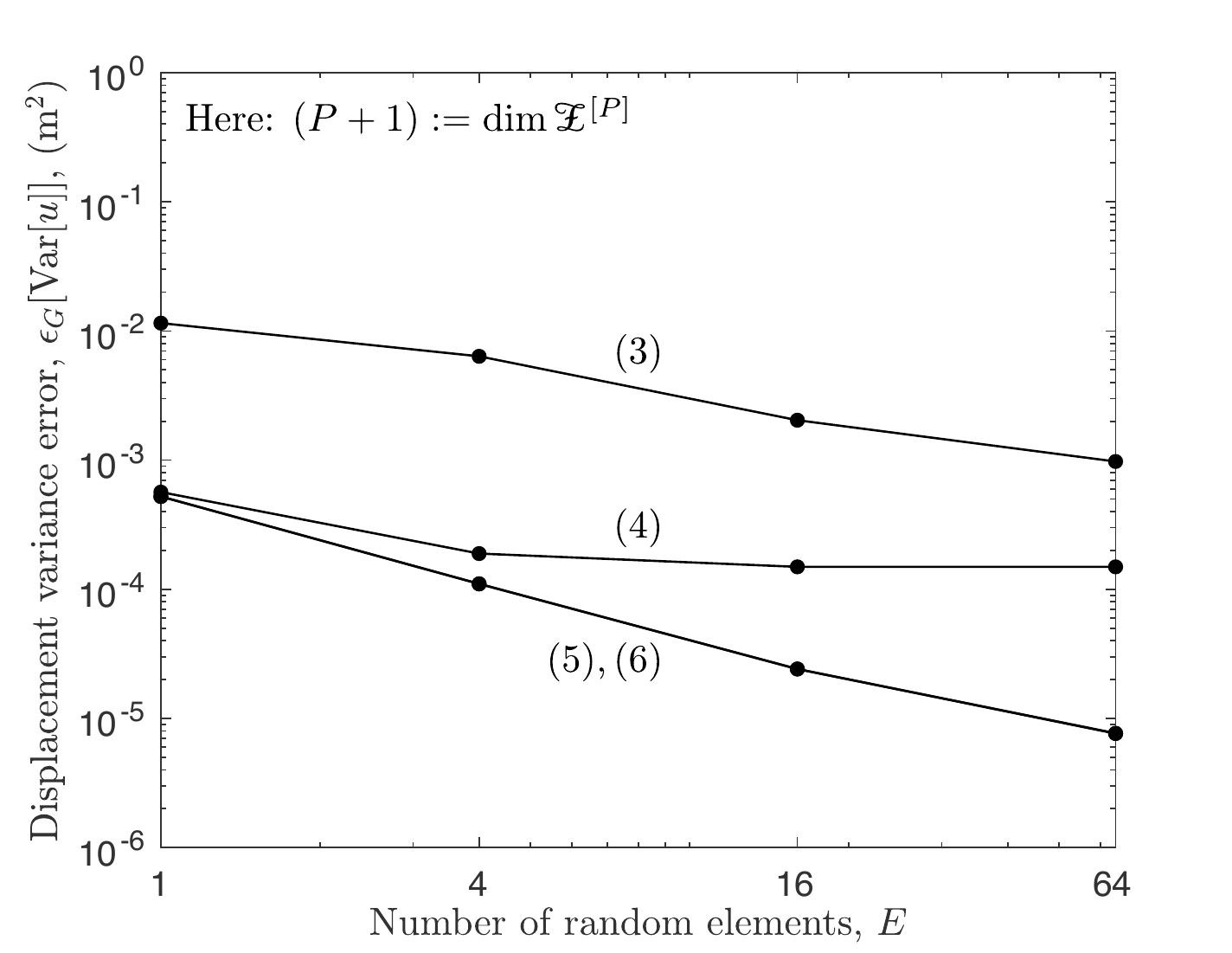}
\caption{Variance error}
\label{fig3Problem3_UniformBeta_MEFSC_Disp_Var_GlobalError}
\end{subfigure}
\caption{\emph{Problem 3} --- Global error of $\mathbf{E}[u]$ and $\mathrm{Var}[u]$ for different $(h,p)$-discretization levels of RFS and for $\mu\sim\mathrm{Uniform}\otimes\mathrm{Beta}$}
\label{fig3Problem3_UniformBeta_MEFSC_Disp_GlobalError}
\end{figure}

Fig.~\ref{fig3Problem3_UniformBeta_MEFSC_5BV_64E_Disp} depicts the evolution of the mean and variance of the system's displacement using ME-FSC and a Monte Carlo simulation with one million realizations.
The reason why a Monte Carlo simulation is used this time as the reference solution is that a closed-form solution for $u$ is not available.
One drawback of using Monte Carlo as the reference solution is that the degree of accuracy of ME-FSC will not be able to be measured adequately, since the solution obtained with ME-FSC can be far more accurate than the one obtained with Monte Carlo.
Yet, from this figure, we learn that when 64 elements and 5 basis vectors are employed to run the simulation, the ME-FSC method is capable of reproducing the Monte Carlo solution quite well.
To compare the level of accuracy of ME-FSC with respect to Monte Carlo, Fig.~\ref{fig3Problem3_UniformBeta_MEFSC_Disp_LocalError} presents the local errors in mean and variance of the system's displacement.
In general, good convergence can be observed when the number of elements and the number of basis vectors are both increased.
This observation can be better verified using Fig.~\ref{fig3Problem3_UniformBeta_MEFSC_Disp_GlobalError} which plots the global errors in mean and variance of the system's displacement.
As expected, the results only improve if the number of basis vectors is increased up to a certain number, which in this case happens to be 5.
However, contrary to what we observed in Problems 1 and 2, increasing the number of elements does not necessarily improve the accuracy of the results, as can be confirmed when $P+1$ is set to 4.
This suggests that due to the nonlinear nature of the system's ODE, higher basis vectors can play a major role in the description of the system's state over time.

% Problem 3 [1/3]:

% Problem 3 [2/3]:

% Problem 3 [3/3]:

\subsection[Problem 4 (nonlinear system)]{Problem 4: A nonlinear system governed by a system of 1st-order stochastic ODEs (the Kraichnan-Orszag three-mode problem)}

In this last problem, we explore the Kraichnan-Orszag three-mode problem \cite{orszag1967dynamical} to test the performance of ME-FSC more thoroughly.
This problem is particularly challenging for methods based on the spectral approach because the solution is known to be discontinuous over the random domain.
It has been used as a benchmark problem in various works (e.g.~\cite{wan2005adaptive,gerritsma2010time,heuveline2014hybrid}), and this is the reason why we opt to study it in this work as well.

The three-mode problem considered herein is the same as that presented in \cite{wan2005adaptive} (Page~635, Section 4.3.4).
The system's governing differential equation is defined as
\begin{align*}
\dot{u}_1&=u_1u_3,\\
\dot{u}_2&=-u_2u_3,\\
\dot{u}_3&=-u_1^2+u_2^2,
\end{align*}
where $u_1,u_2,u_3:\mathfrak{T}\times\Xi\to\mathbb{R}$ represent the three modes of the system, and $\dot{u}_1:=\partial_t u_1,\dot{u}_2:=\partial_t u_2,\dot{u}_3:=\partial_t u_3$ are the corresponding velocities.
In this problem, we take the initial conditions of the system to be stochastic and given by: $u_1(0,\xi)=\xi^1$, $u_2(0,\xi)=\xi^2$, and $u_3(0,\xi)=\xi^3$.

Two probability distributions are investigated for $\xi=(\xi^1,\xi^2,\xi^3)$.
The first one is a \emph{uniform distribution} defined by $\mathrm{Uniform}^{\otimes3}\sim\xi\in\Xi=[a,b]^3$, and the second one is a \emph{beta distribution} defined by  $\mathrm{Beta}(\alpha,\beta)^{\otimes3}\sim\xi\in\Xi=[a,b]^3$, from where $(a,b)=(-1,1)$ and $(\alpha,\beta)=(2,5)$.
The random space is thus three-dimensional.

\begin{figure}
\centering
\begin{subfigure}[b]{0.495\textwidth}
\includegraphics[width=\textwidth]{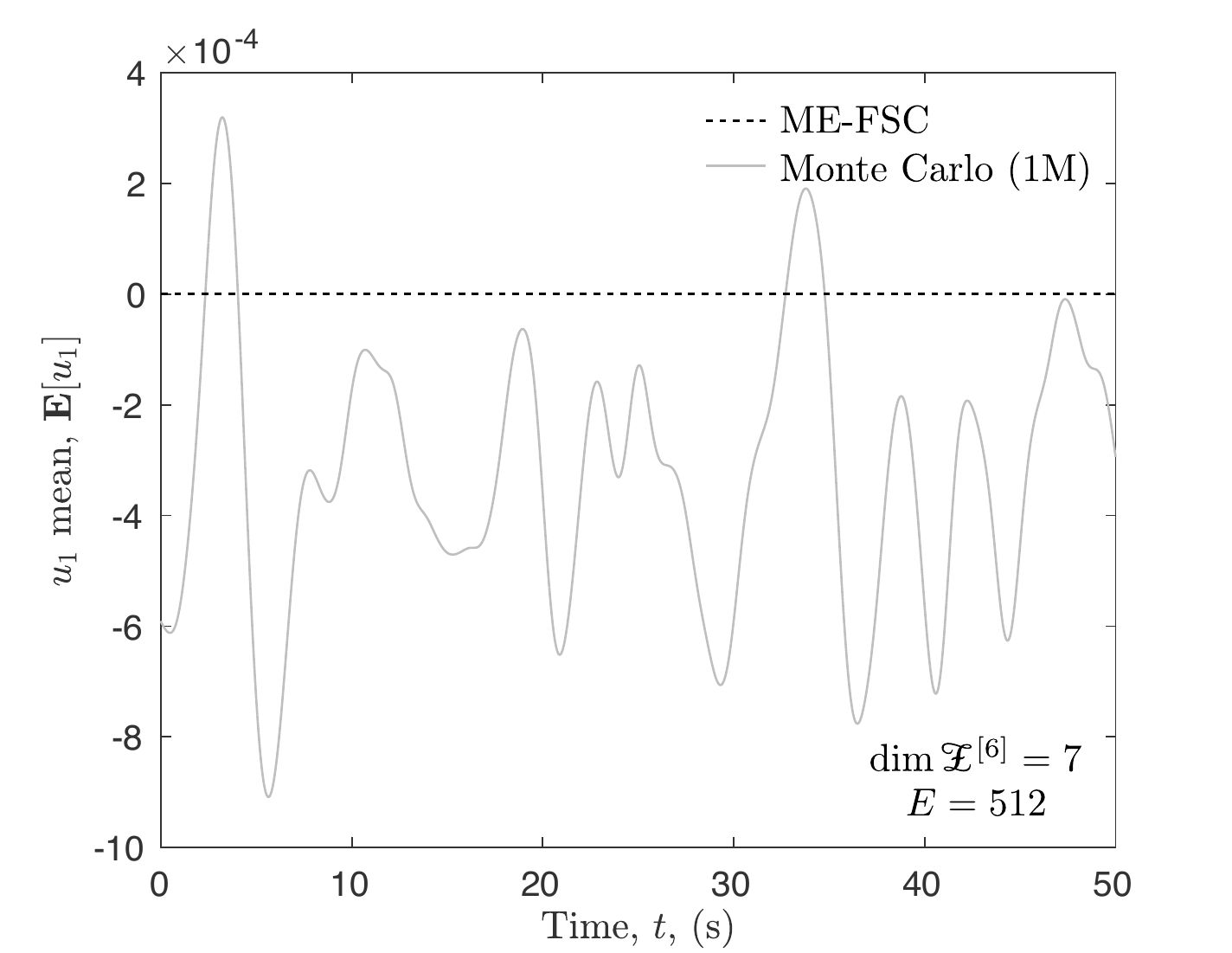}
\caption{Mean for $u_1$}
\label{fig3Problem4_Uniform3_MEFSC_7BV_512E_u1_Mean}
\end{subfigure}\hfill
\begin{subfigure}[b]{0.495\textwidth}
\includegraphics[width=\textwidth]{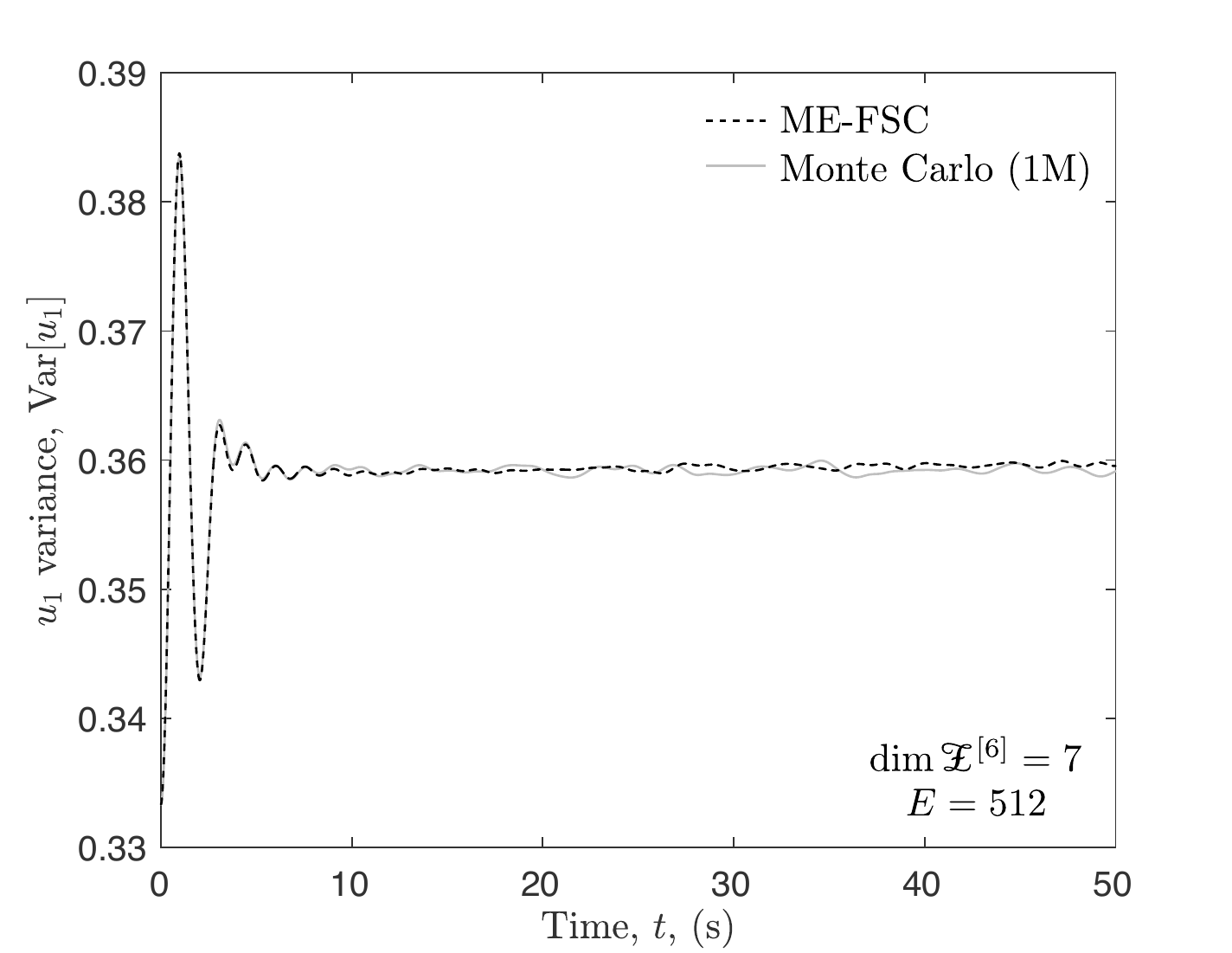}
\caption{Variance for $u_1$}
\label{fig3Problem4_Uniform3_MEFSC_7BV_512E_u1_Var}
\end{subfigure}\quad
\begin{subfigure}[b]{0.495\textwidth}
\includegraphics[width=\textwidth]{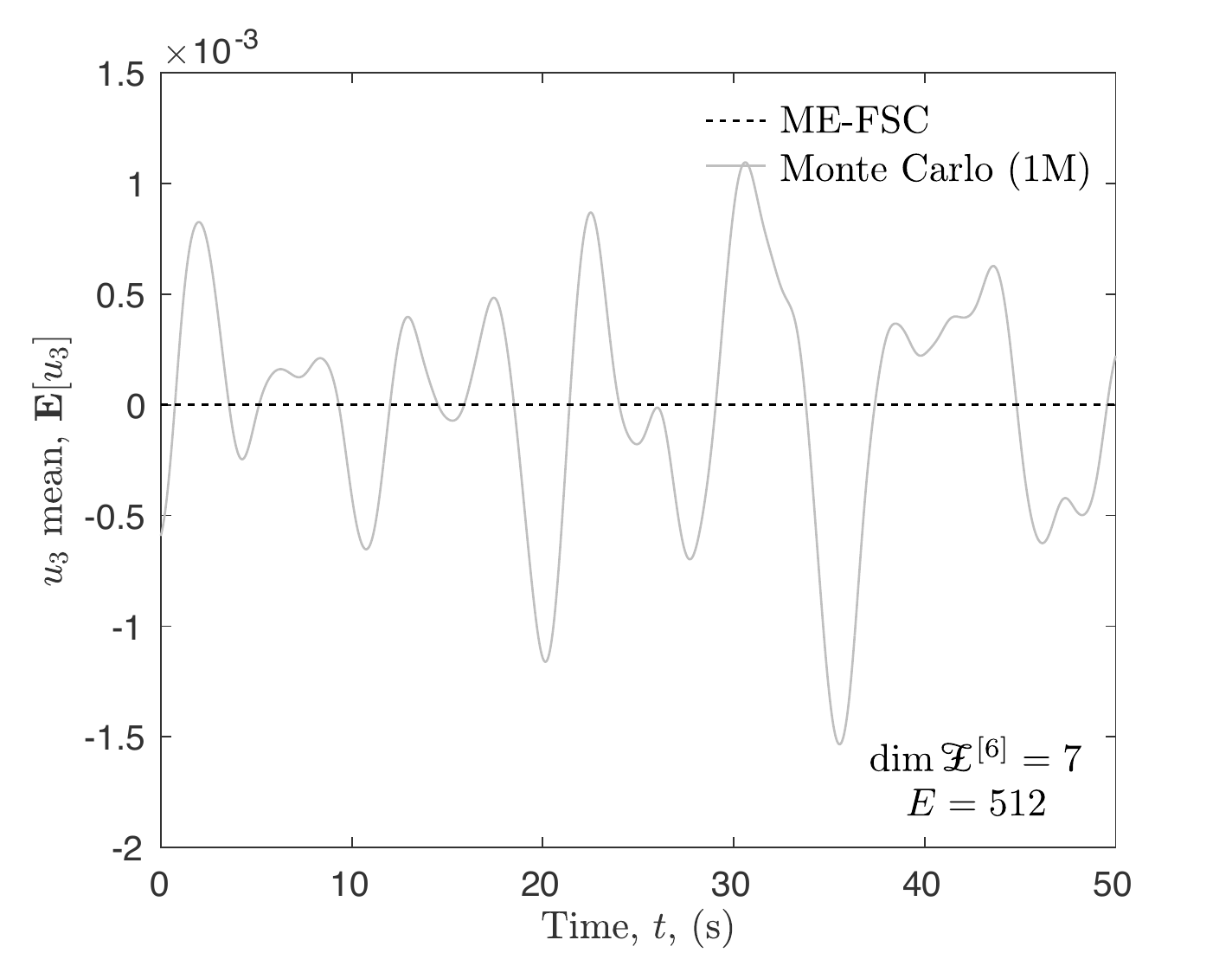}
\caption{Mean for $u_3$}
\label{fig3Problem4_Uniform3_MEFSC_7BV_512E_u3_Mean}
\end{subfigure}\hfill
\begin{subfigure}[b]{0.495\textwidth}
\includegraphics[width=\textwidth]{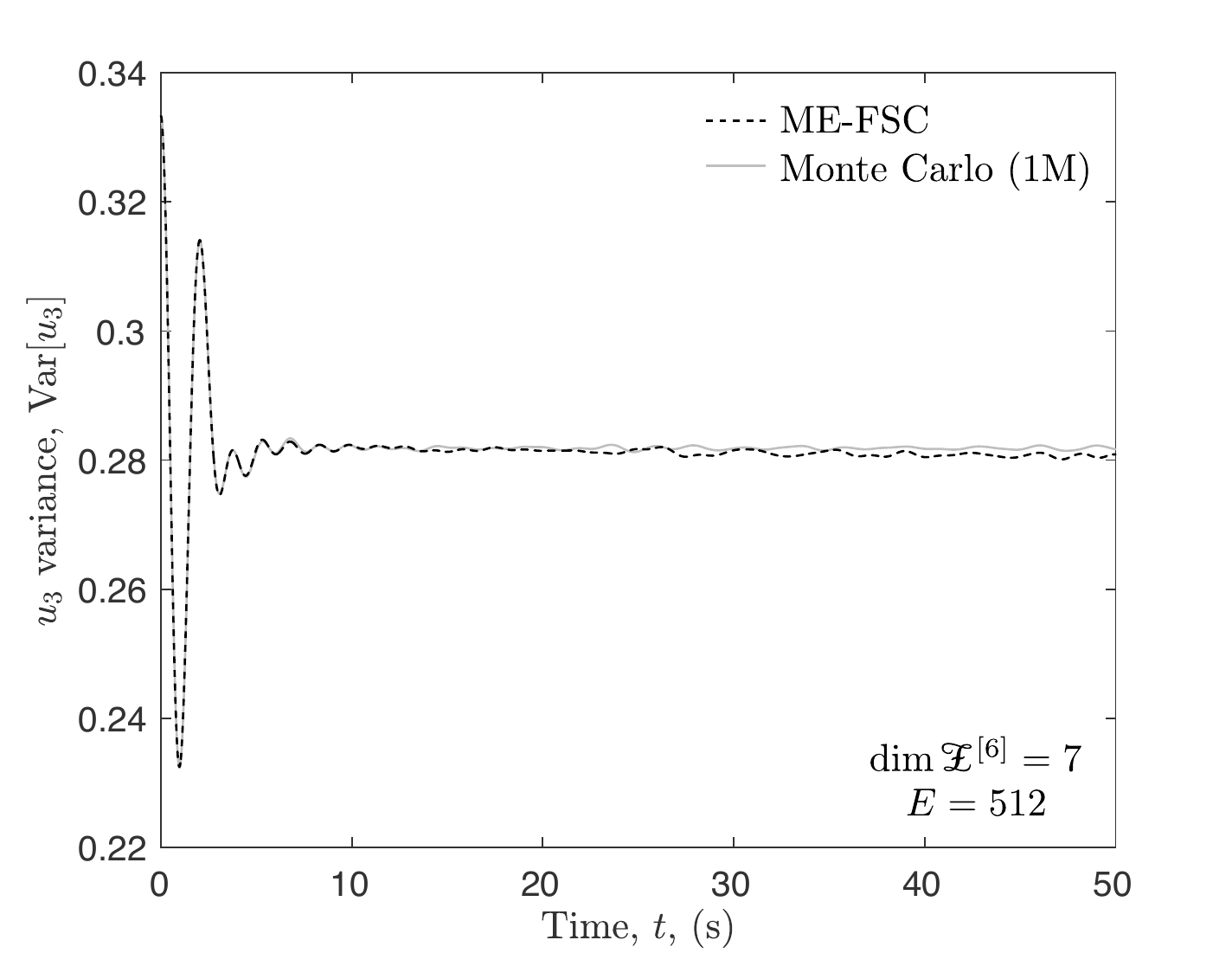}
\caption{Variance for $u_3$}
\label{fig3Problem4_Uniform3_MEFSC_7BV_512E_u3_Var}
\end{subfigure}
\caption{\emph{Problem 4} --- Evolution of $\mathbf{E}[u_1]$, $\mathrm{Var}[u_1]$, $\mathbf{E}[u_3]$ and $\mathrm{Var}[u_3]$ for the case when the $(h,p)$-discretization level of RFS is $(P,E)=(6,512)$ and $\mu\sim\mathrm{Uniform}^{\otimes3}$}
\label{fig3Problem4_Uniform3_MEFSC_7BV_512E}
\end{figure}

\begin{figure}
\centering
\begin{subfigure}[b]{0.495\textwidth}
\includegraphics[width=\textwidth]{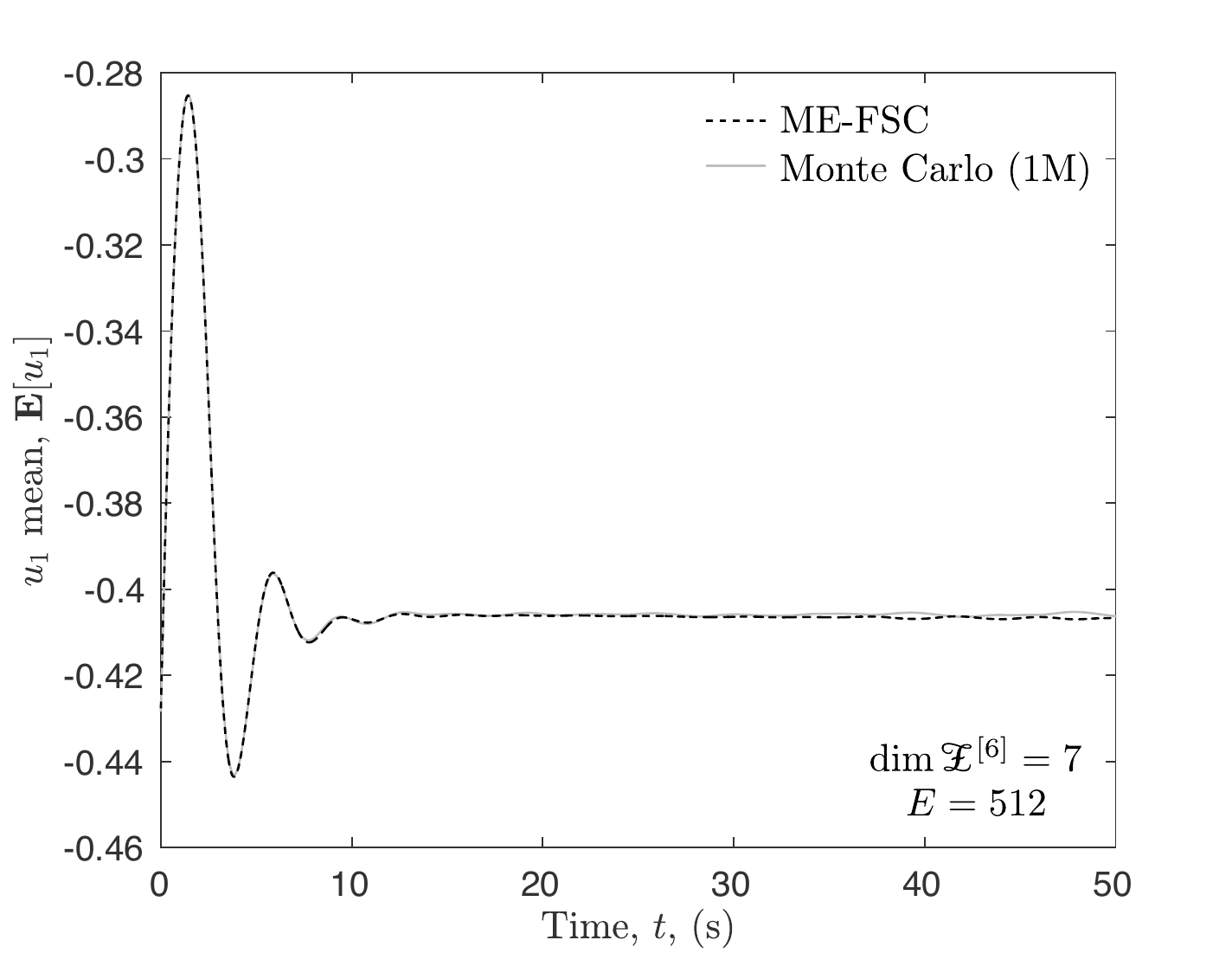}
\caption{Mean for $u_1$}
\label{fig3Problem4_Beta3_MEFSC_7BV_512E_u1_Mean}
\end{subfigure}\hfill
\begin{subfigure}[b]{0.495\textwidth}
\includegraphics[width=\textwidth]{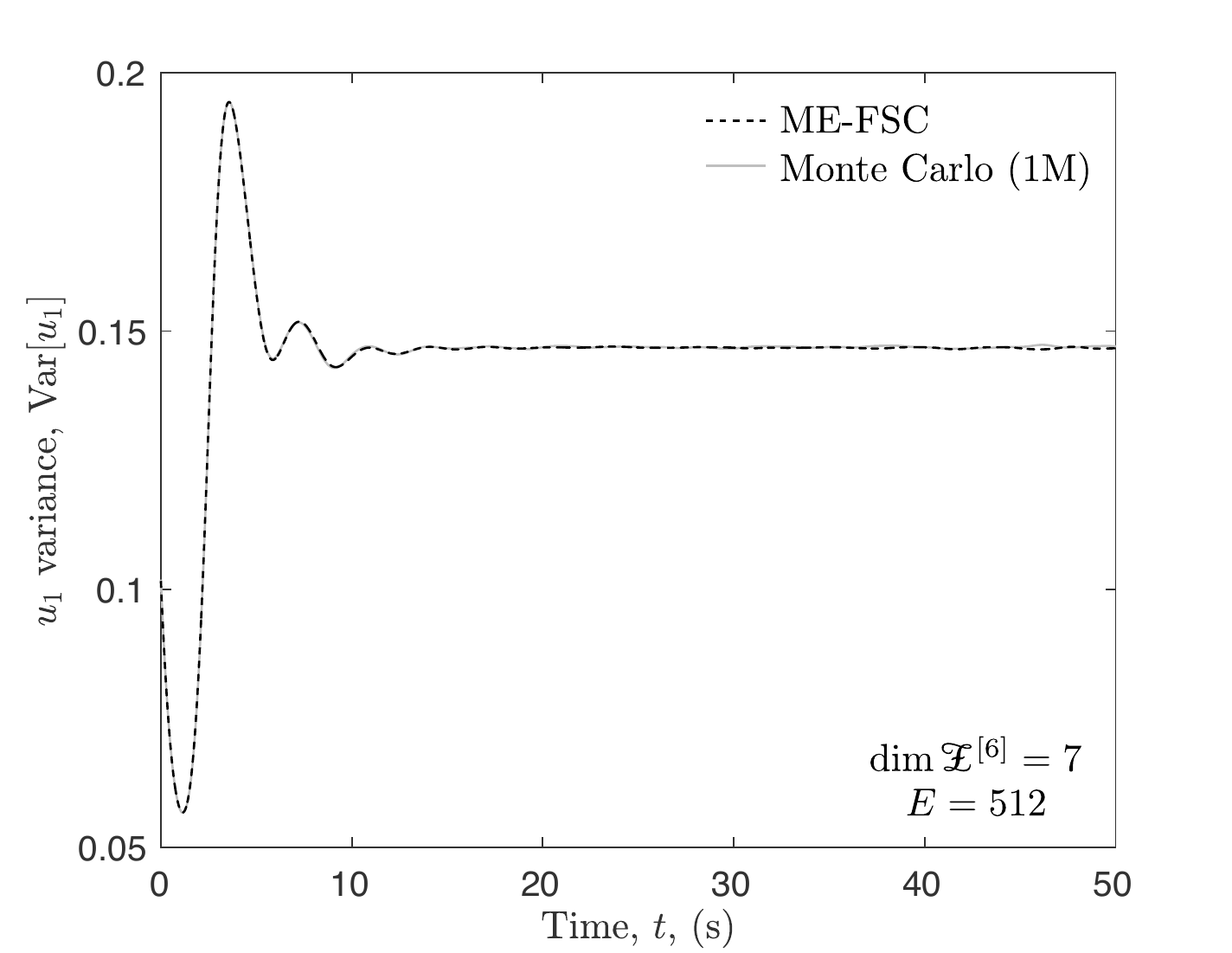}
\caption{Variance for $u_1$}
\label{fig3Problem4_Beta3_MEFSC_7BV_512E_u1_Var}
\end{subfigure}\quad
\begin{subfigure}[b]{0.495\textwidth}
\includegraphics[width=\textwidth]{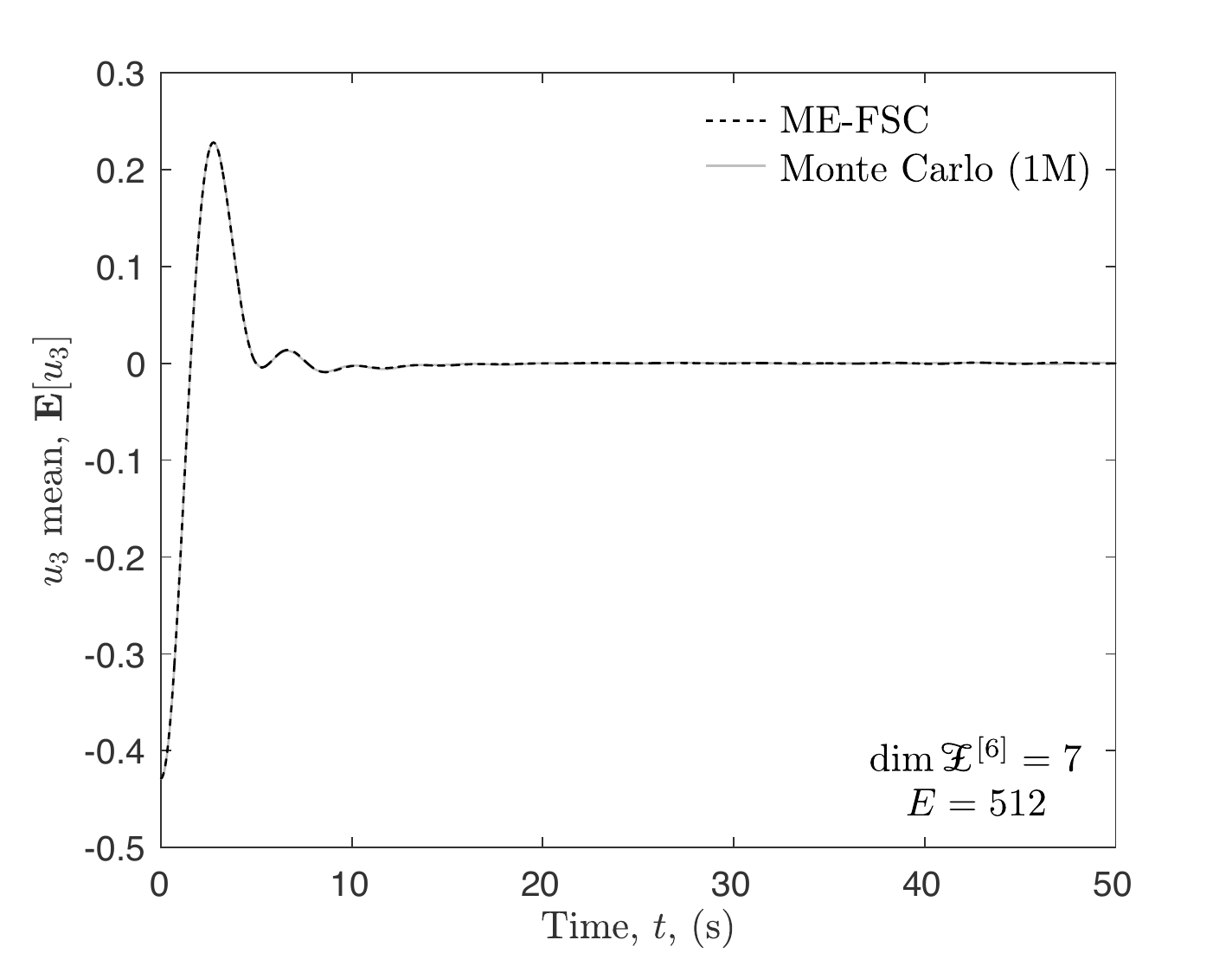}
\caption{Mean for $u_3$}
\label{fig3Problem4_Beta3_MEFSC_7BV_512E_u3_Mean}
\end{subfigure}\hfill
\begin{subfigure}[b]{0.495\textwidth}
\includegraphics[width=\textwidth]{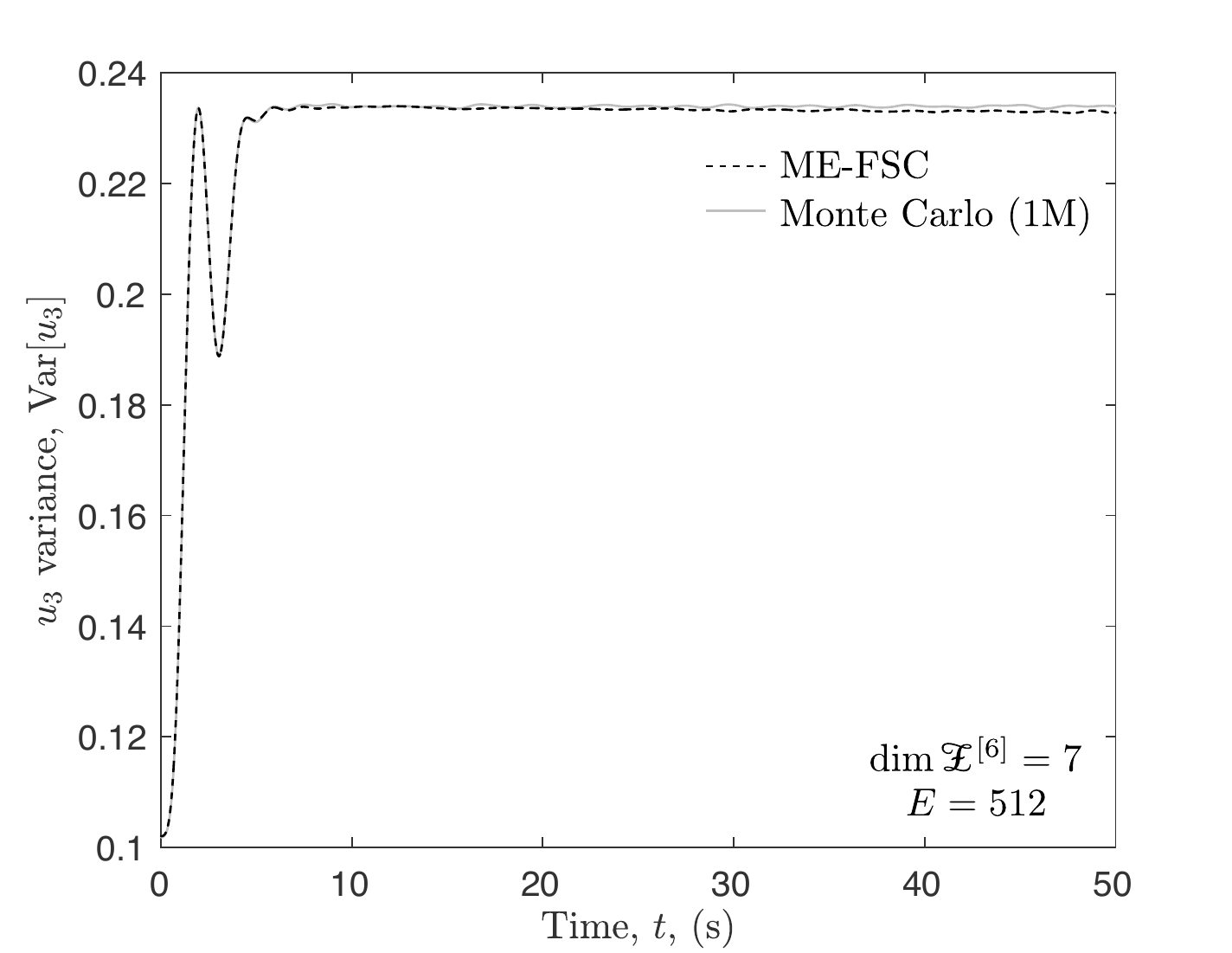}
\caption{Variance for $u_3$}
\label{fig3Problem4_Beta3_MEFSC_7BV_512E_u3_Var}
\end{subfigure}
\caption{\emph{Problem 4} --- Evolution of $\mathbf{E}[u_1]$, $\mathrm{Var}[u_1]$, $\mathbf{E}[u_3]$ and $\mathrm{Var}[u_3]$ for the case when the $(h,p)$-discretization level of RFS is $(P,E)=(6,512)$ and $\mu\sim\mathrm{Beta}^{\otimes3}$}
\label{fig3Problem4_Beta3_MEFSC_7BV_512E}
\end{figure}

In Figs.~\ref{fig3Problem4_Uniform3_MEFSC_7BV_512E} and \ref{fig3Problem4_Beta3_MEFSC_7BV_512E} we depict the evolution of the mean and variance of $u_1$ and $u_3$ for the two distributions chosen for $\xi$.
The results are obtained for the case of using 512 elements and 7 basis vectors per element, and then compared against a Monte Carlo simulation with one million realizations.
As observed, the ME-FSC solution is in good agreement with the Monte Carlo solution for all cases considered, except for the mean of $u_1$ and $u_3$ under uniform distributions (Figs.~\ref{fig3Problem4_Uniform3_MEFSC_7BV_512E_u1_Mean} and \ref{fig3Problem4_Uniform3_MEFSC_7BV_512E_u3_Mean}).
It is worth commenting that the exact mean of $u_1$, $u_2$, and $u_3$ are known to be identically equal to zero when the probability measure is uniform.
This demonstrates that the ME-FSC solution can be far more accurate than a Monte Carlo simulation with one million realizations.

% Problem 4 [1a/3]:

% Problem 4 [1b/3]:

\begin{figure}
\centering
\begin{subfigure}[b]{0.495\textwidth}
\includegraphics[width=\textwidth]{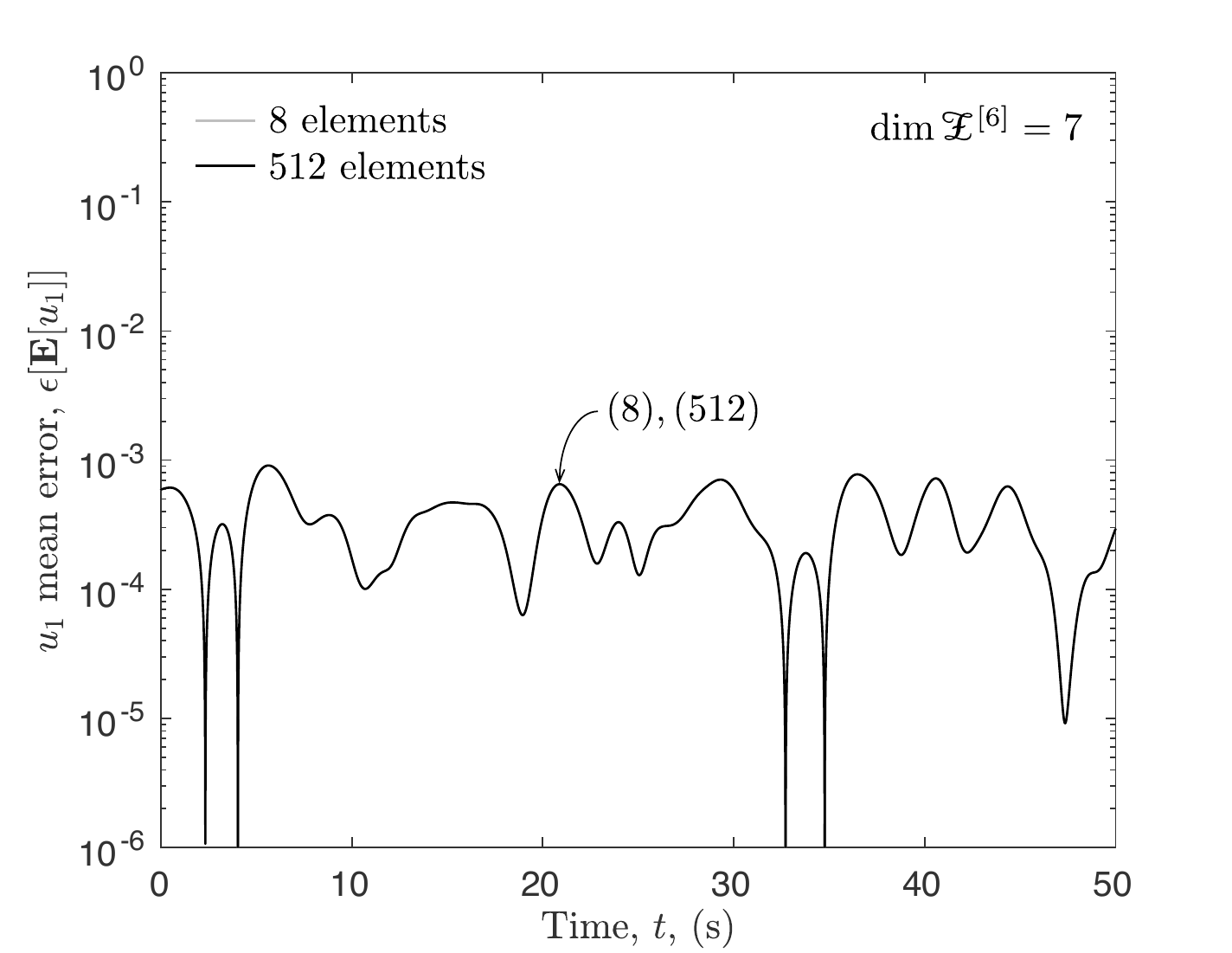}
\caption{Mean error for $u_1$}
\label{fig3Problem4_Uniform3_MEFSC_7BV_u1_Mean_LocalError}
\end{subfigure}\hfill
\begin{subfigure}[b]{0.495\textwidth}
\includegraphics[width=\textwidth]{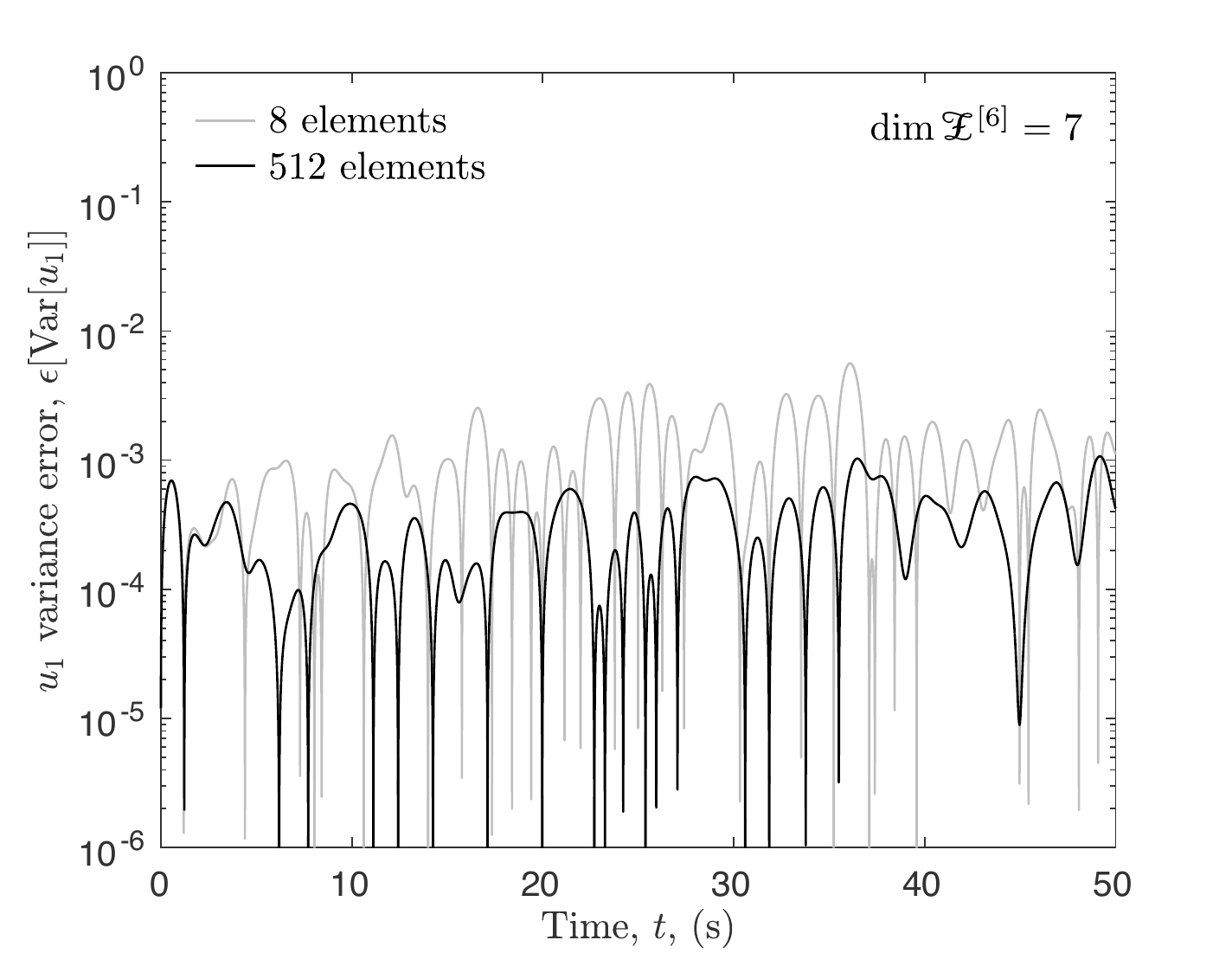}
\caption{Variance error for $u_1$}
\label{fig3Problem4_Uniform3_MEFSC_7BV_u1_Var_LocalError}
\end{subfigure}\quad
\begin{subfigure}[b]{0.495\textwidth}
\includegraphics[width=\textwidth]{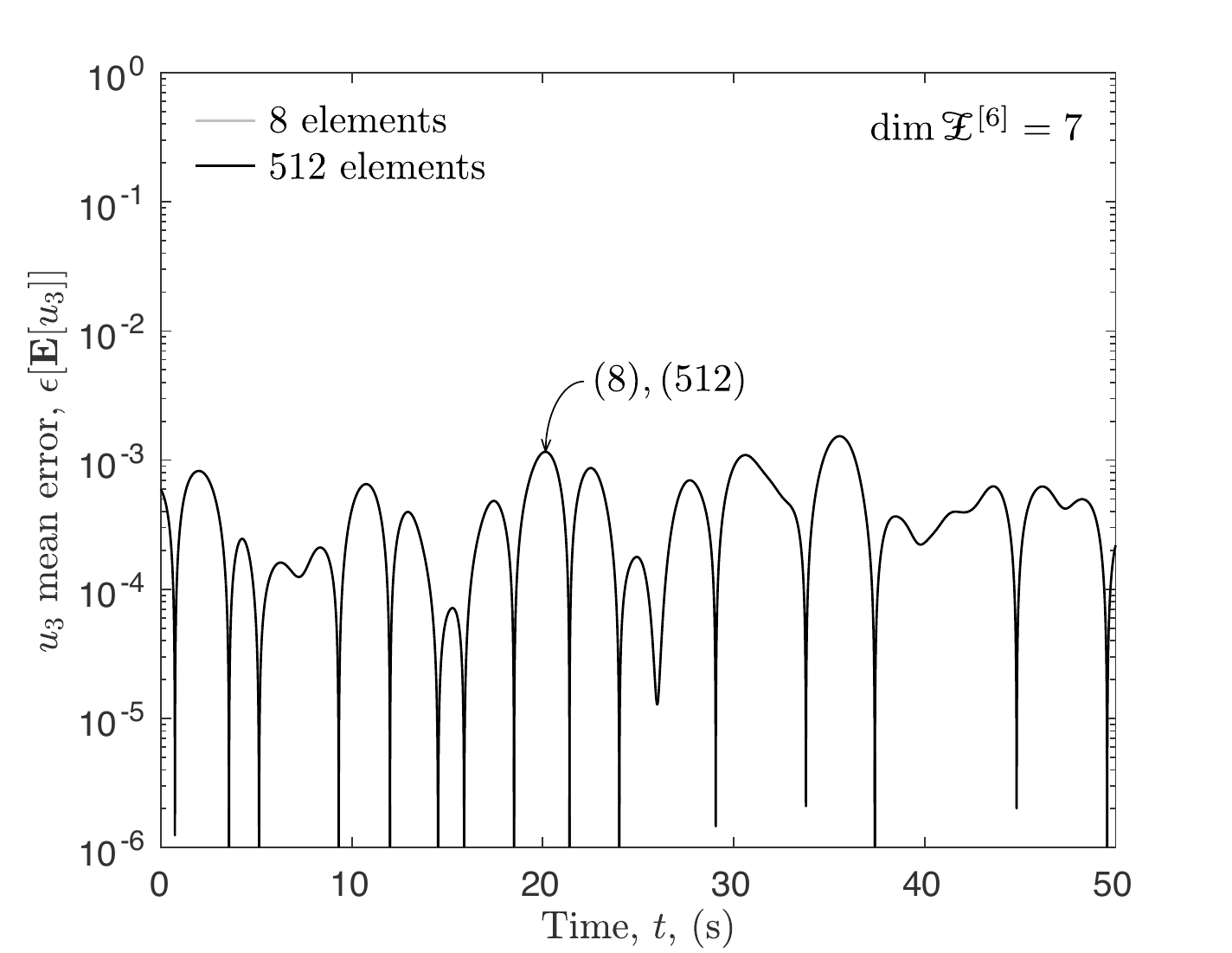}
\caption{Mean error for $u_3$}
\label{fig3Problem4_Uniform3_MEFSC_7BV_u3_Mean_LocalError}
\end{subfigure}\hfill
\begin{subfigure}[b]{0.495\textwidth}
\includegraphics[width=\textwidth]{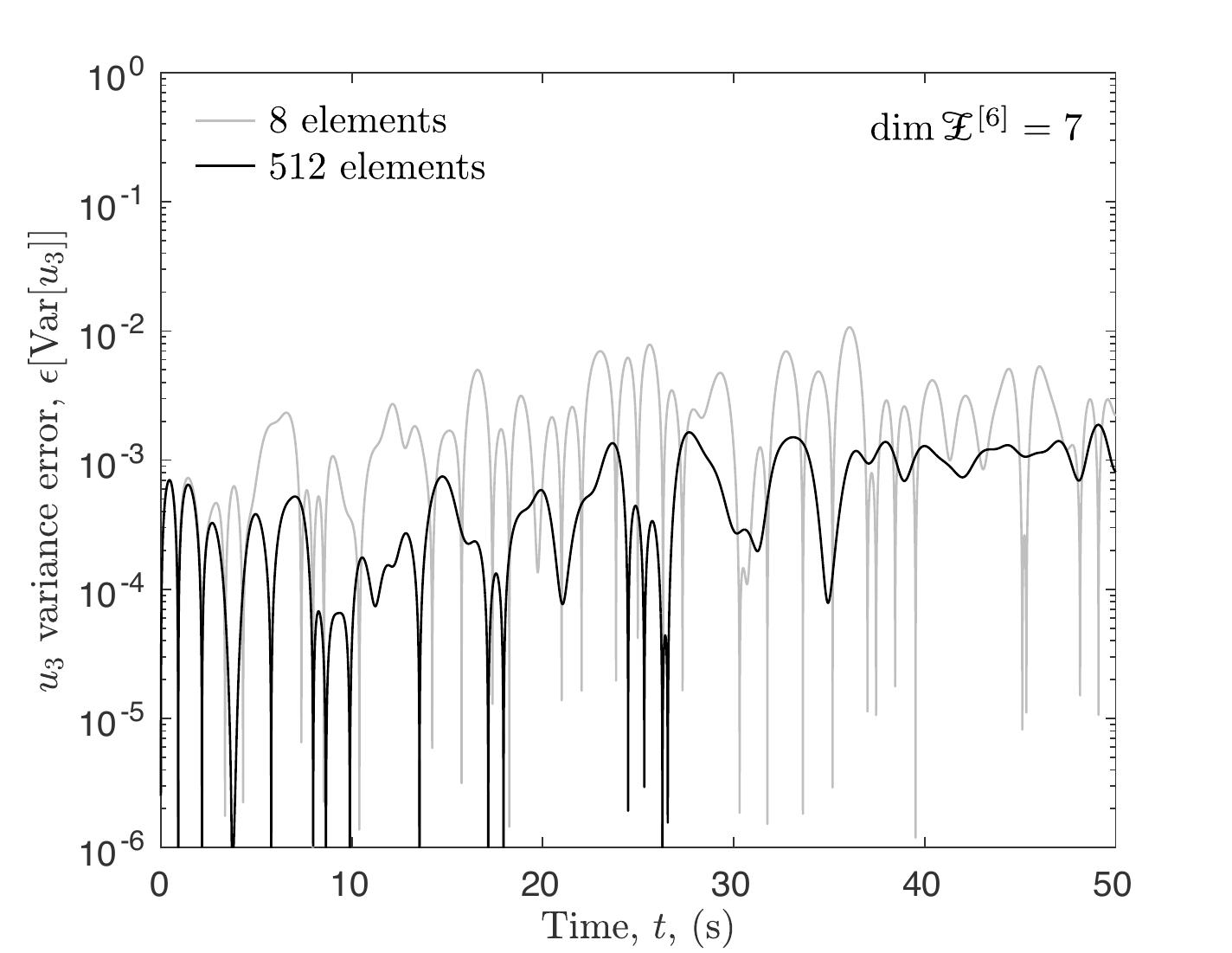}
\caption{Variance error for $u_3$}
\label{fig3Problem4_Uniform3_MEFSC_7BV_u3_Var_LocalError}
\end{subfigure}
\caption{\emph{Problem 4} --- Local error evolution of $\mathbf{E}[u_1]$, $\mathrm{Var}[u_1]$, $\mathbf{E}[u_3]$ and $\mathrm{Var}[u_3]$ for different $(h,p)$-discretization levels of RFS and for $\mu\sim\mathrm{Uniform}^{\otimes3}$}
\label{fig3Problem4_Uniform3_MEFSC_7BV_LocalError}
\end{figure}

\begin{figure}
\centering
\begin{subfigure}[b]{0.495\textwidth}
\includegraphics[width=\textwidth]{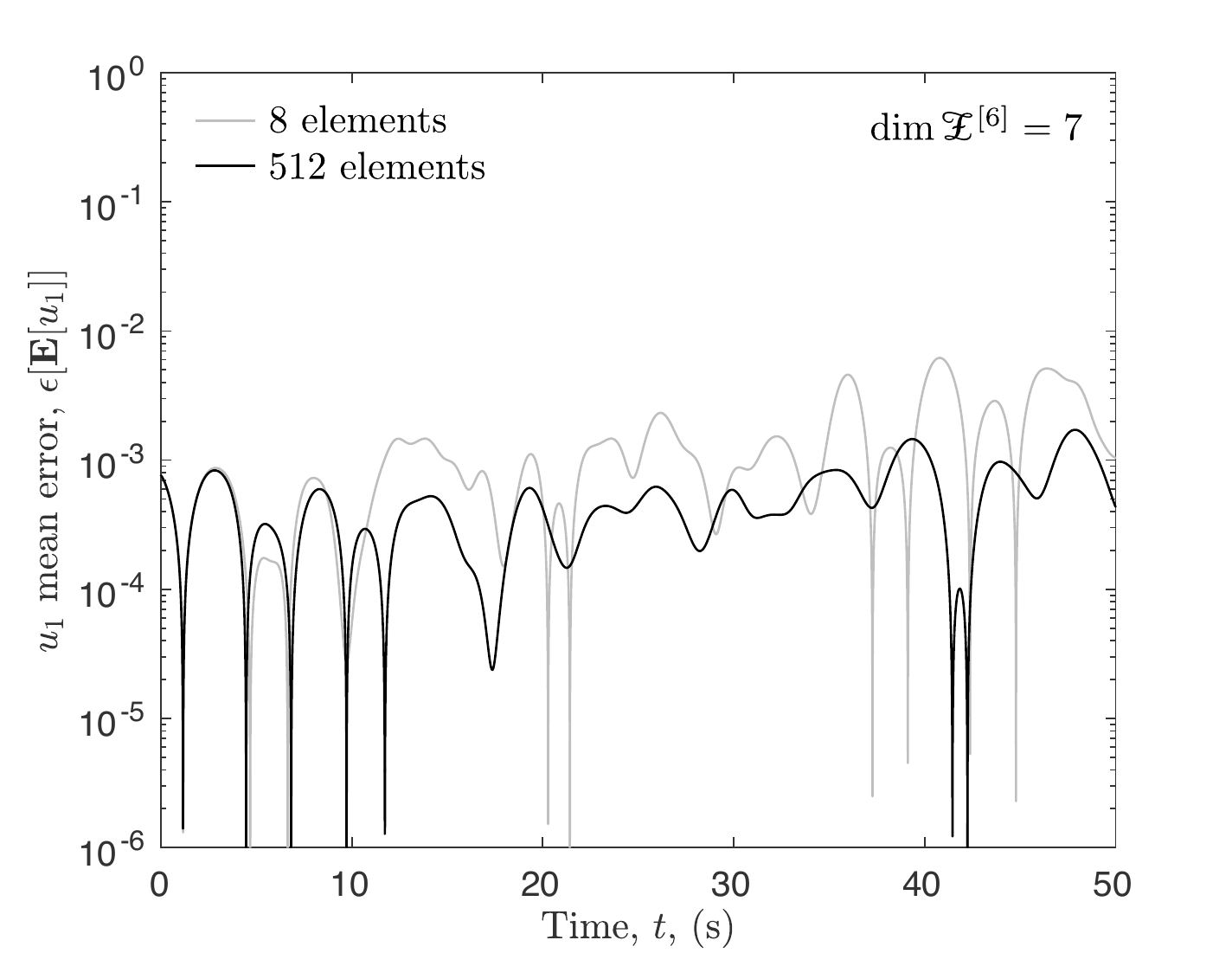}
\caption{Mean error for $u_1$}
\label{fig3Problem4_Beta3_MEFSC_7BV_u1_Mean_LocalError}
\end{subfigure}\hfill
\begin{subfigure}[b]{0.495\textwidth}
\includegraphics[width=\textwidth]{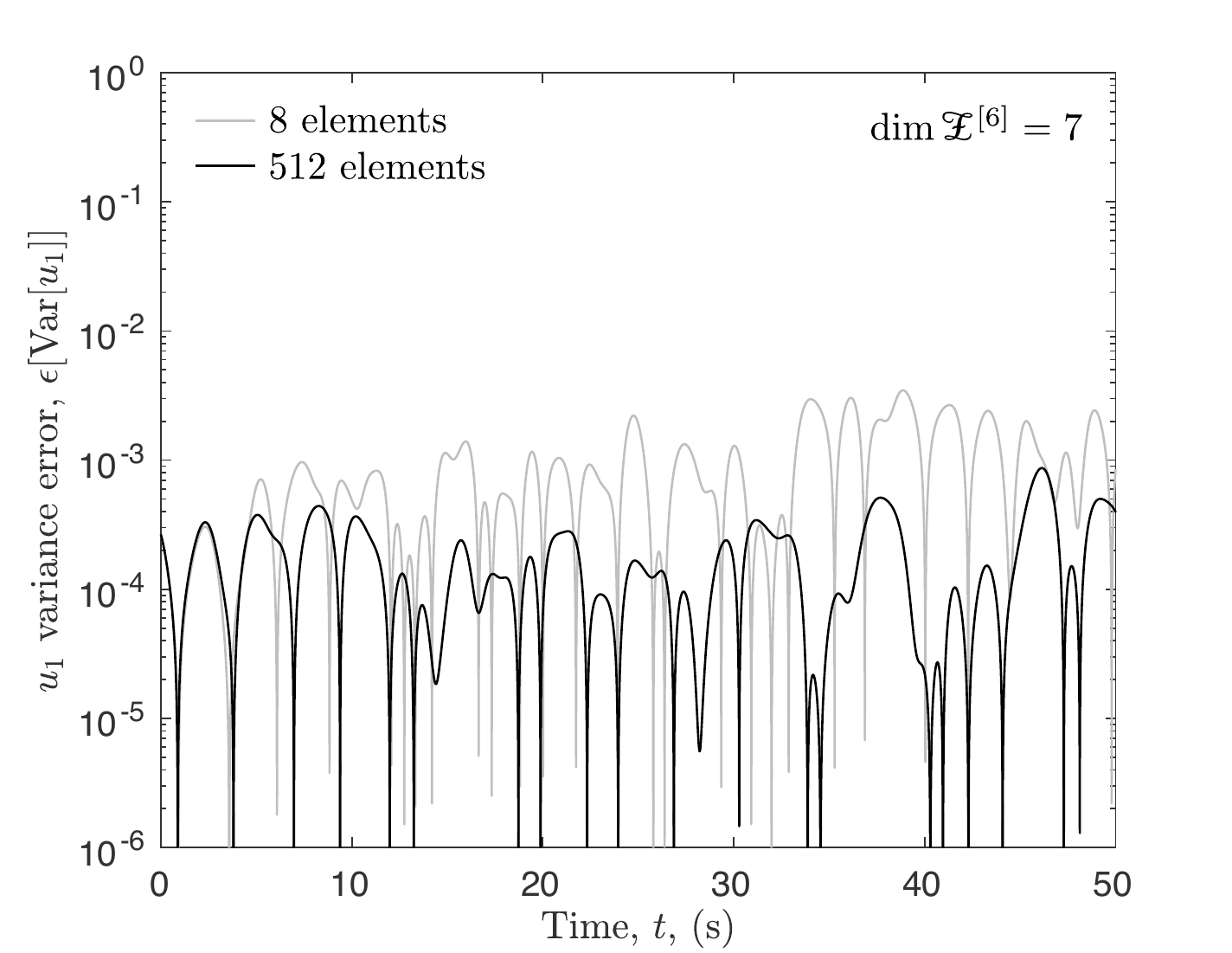}
\caption{Variance error for $u_1$}
\label{fig3Problem4_Beta3_MEFSC_7BV_u1_Var_LocalError}
\end{subfigure}\quad
\begin{subfigure}[b]{0.495\textwidth}
\includegraphics[width=\textwidth]{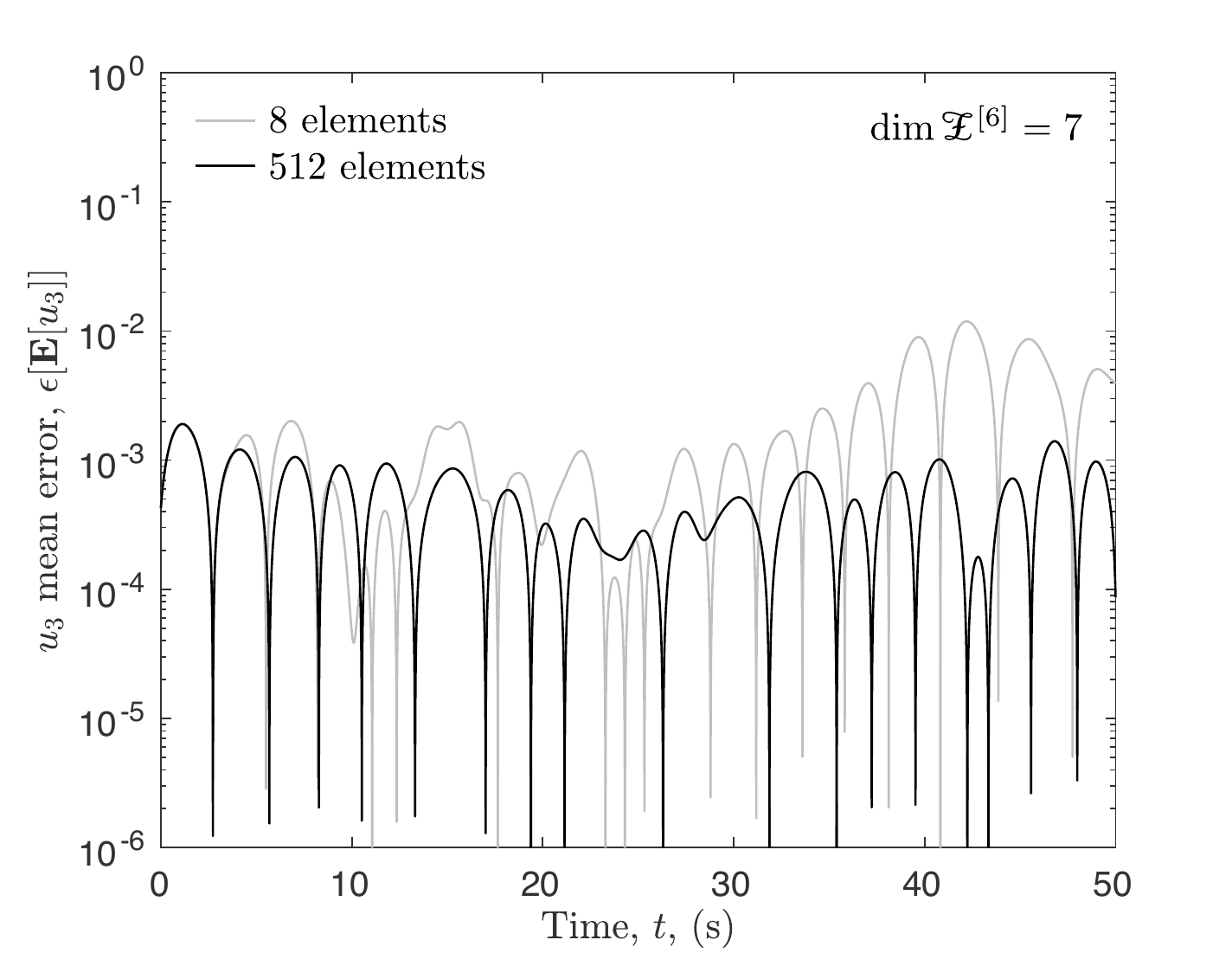}
\caption{Mean error for $u_3$}
\label{fig3Problem4_Beta3_MEFSC_7BV_u3_Mean_LocalError}
\end{subfigure}\hfill
\begin{subfigure}[b]{0.495\textwidth}
\includegraphics[width=\textwidth]{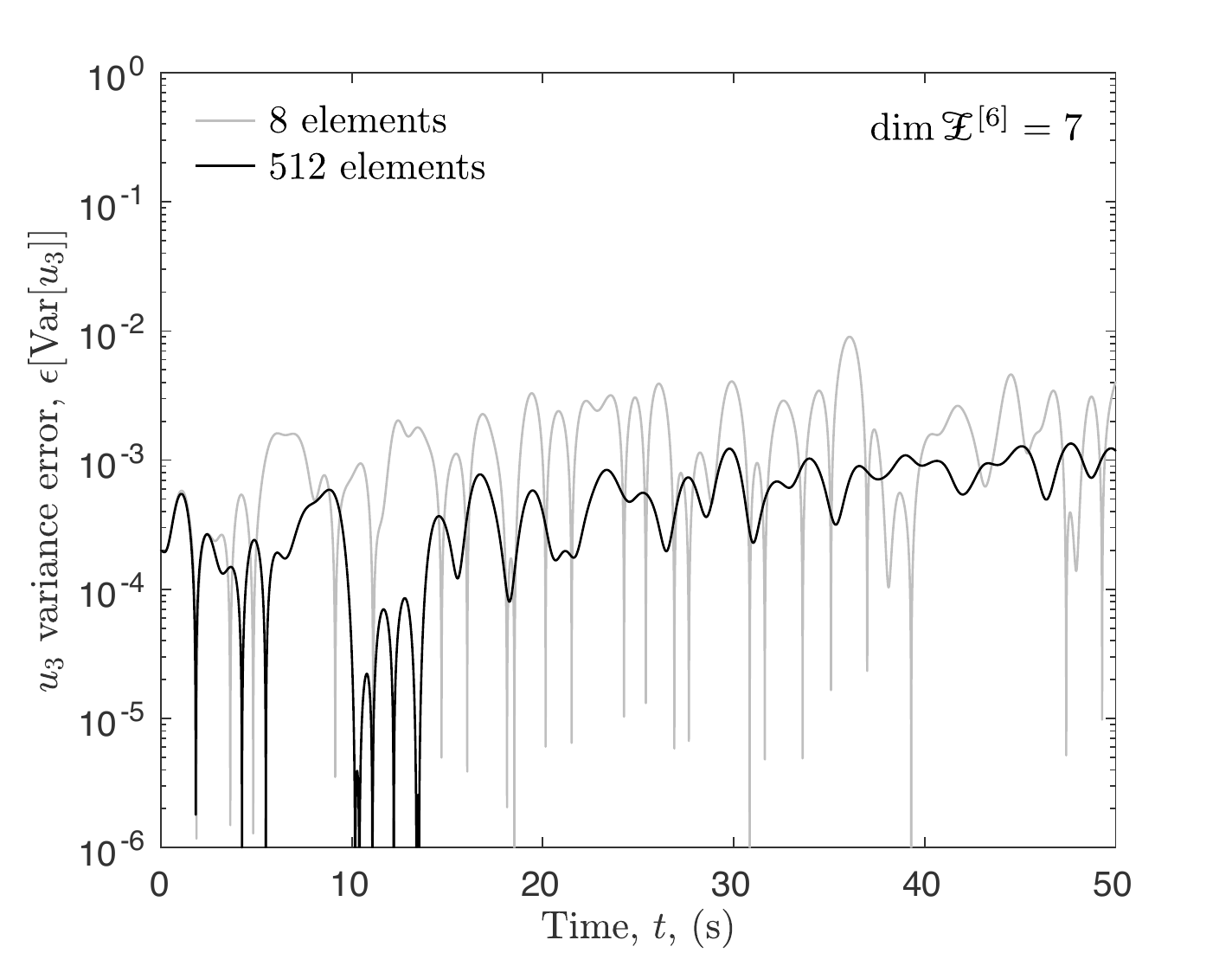}
\caption{Variance error for $u_3$}
\label{fig3Problem4_Beta3_MEFSC_7BV_u3_Var_LocalError}
\end{subfigure}
\caption{\emph{Problem 4} --- Local error evolution of $\mathbf{E}[u_1]$, $\mathrm{Var}[u_1]$, $\mathbf{E}[u_3]$ and $\mathrm{Var}[u_3]$ for different $(h,p)$-discretization levels of RFS and for $\mu\sim\mathrm{Beta}^{\otimes3}$}
\label{fig3Problem4_Beta3_MEFSC_7BV_LocalError}
\end{figure}

\begin{figure}
\centering
\begin{subfigure}[b]{0.495\textwidth}
\includegraphics[width=\textwidth]{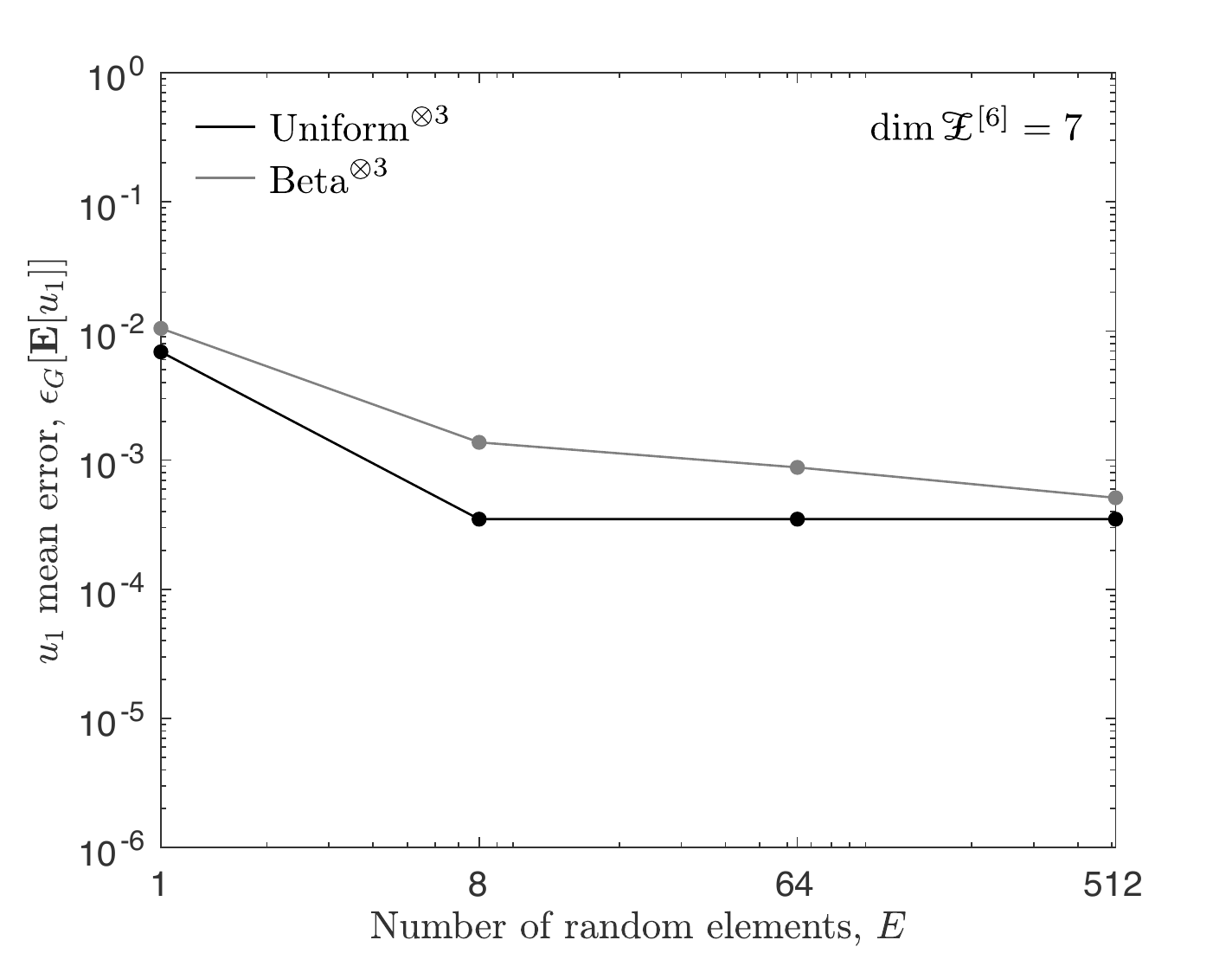}
\caption{Mean error for $u_1$}
\label{fig3Problem4_MEFSC_7BV_u1_Mean_GlobalError}
\end{subfigure}\hfill
\begin{subfigure}[b]{0.495\textwidth}
\includegraphics[width=\textwidth]{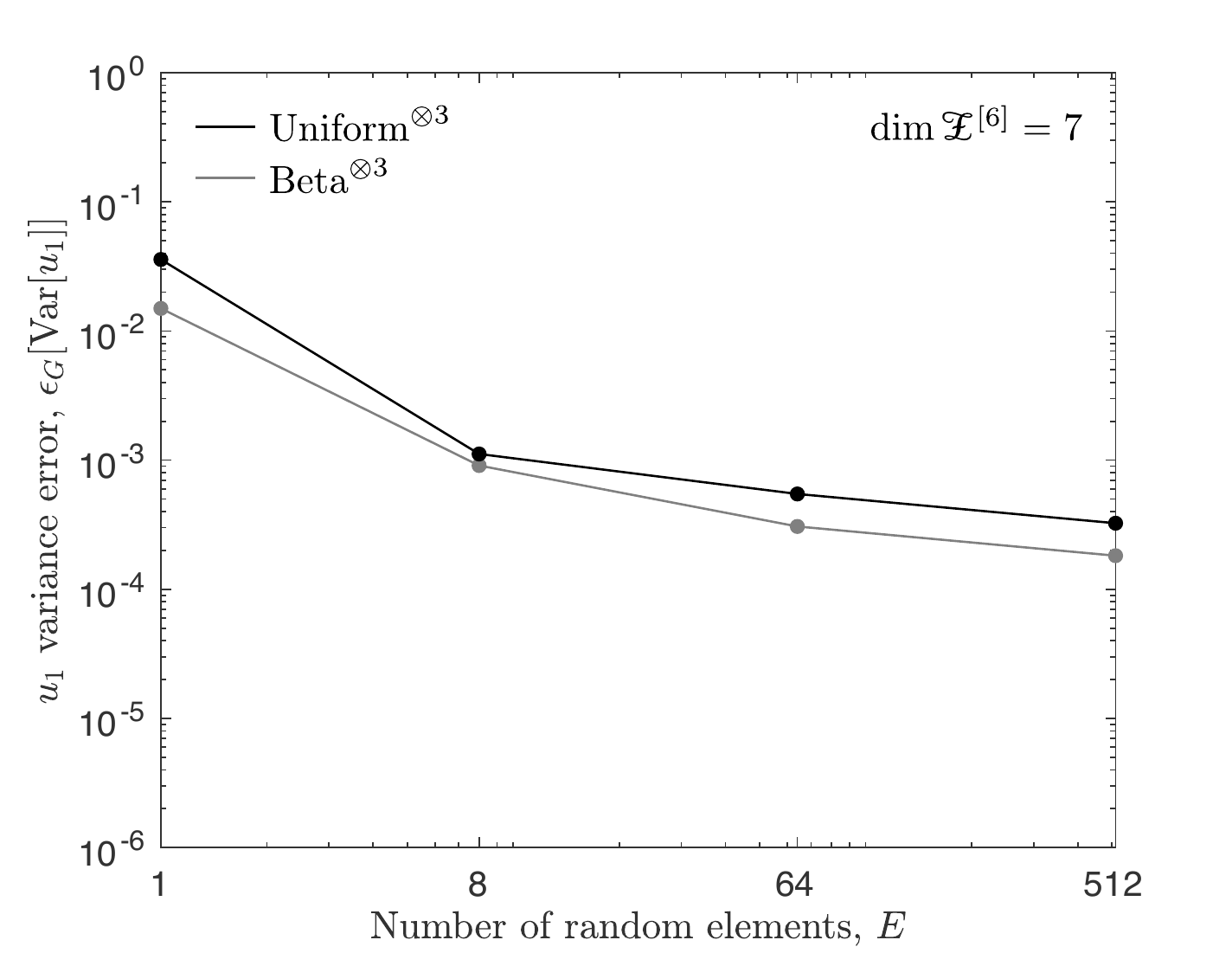}
\caption{Variance error for $u_1$}
\label{fig3Problem4_MEFSC_7BV_u1_Var_GlobalError}
\end{subfigure}\quad
\begin{subfigure}[b]{0.495\textwidth}
\includegraphics[width=\textwidth]{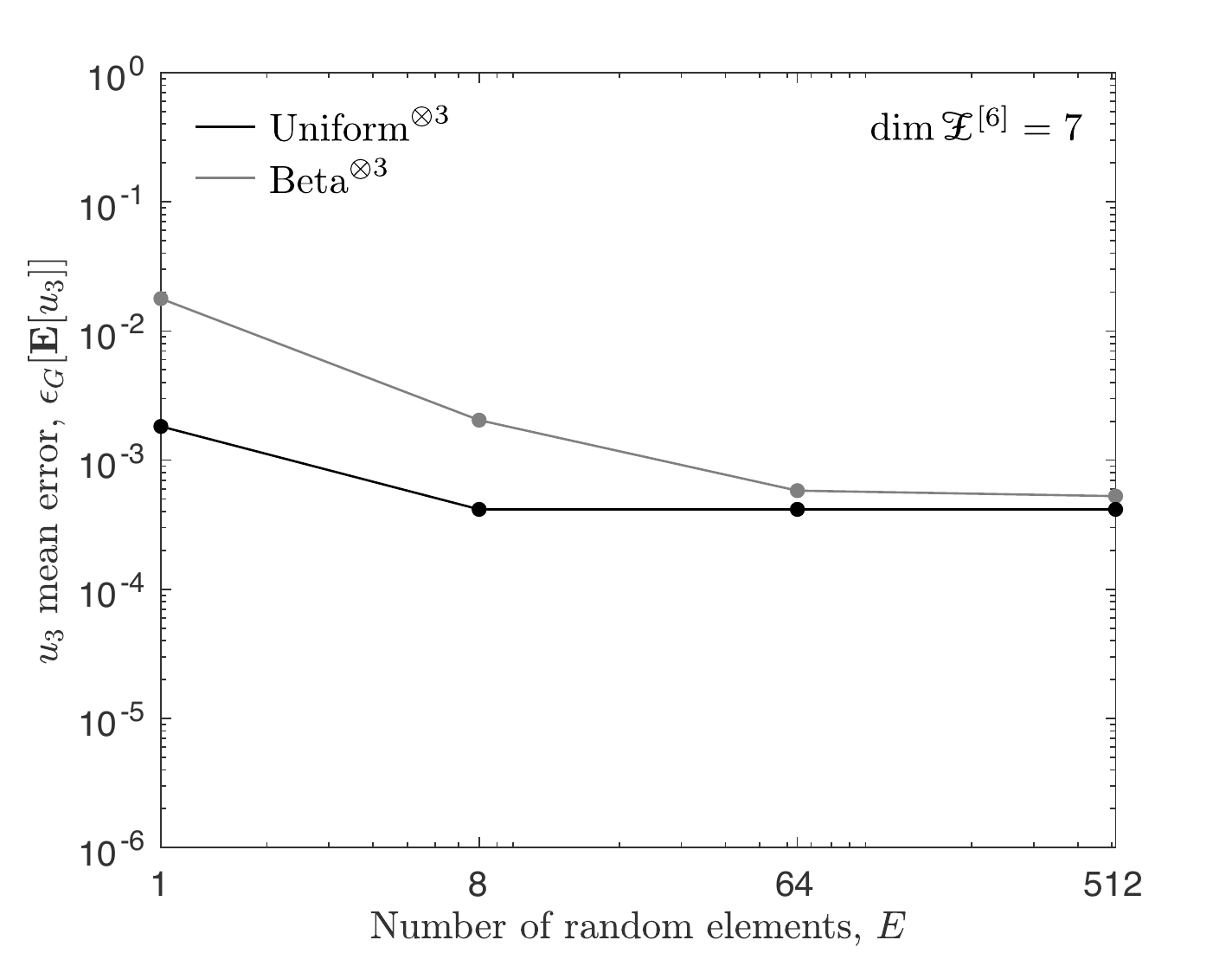}
\caption{Mean error for $u_3$}
\label{fig3Problem4_MEFSC_7BV_u3_Mean_GlobalError}
\end{subfigure}\hfill
\begin{subfigure}[b]{0.495\textwidth}
\includegraphics[width=\textwidth]{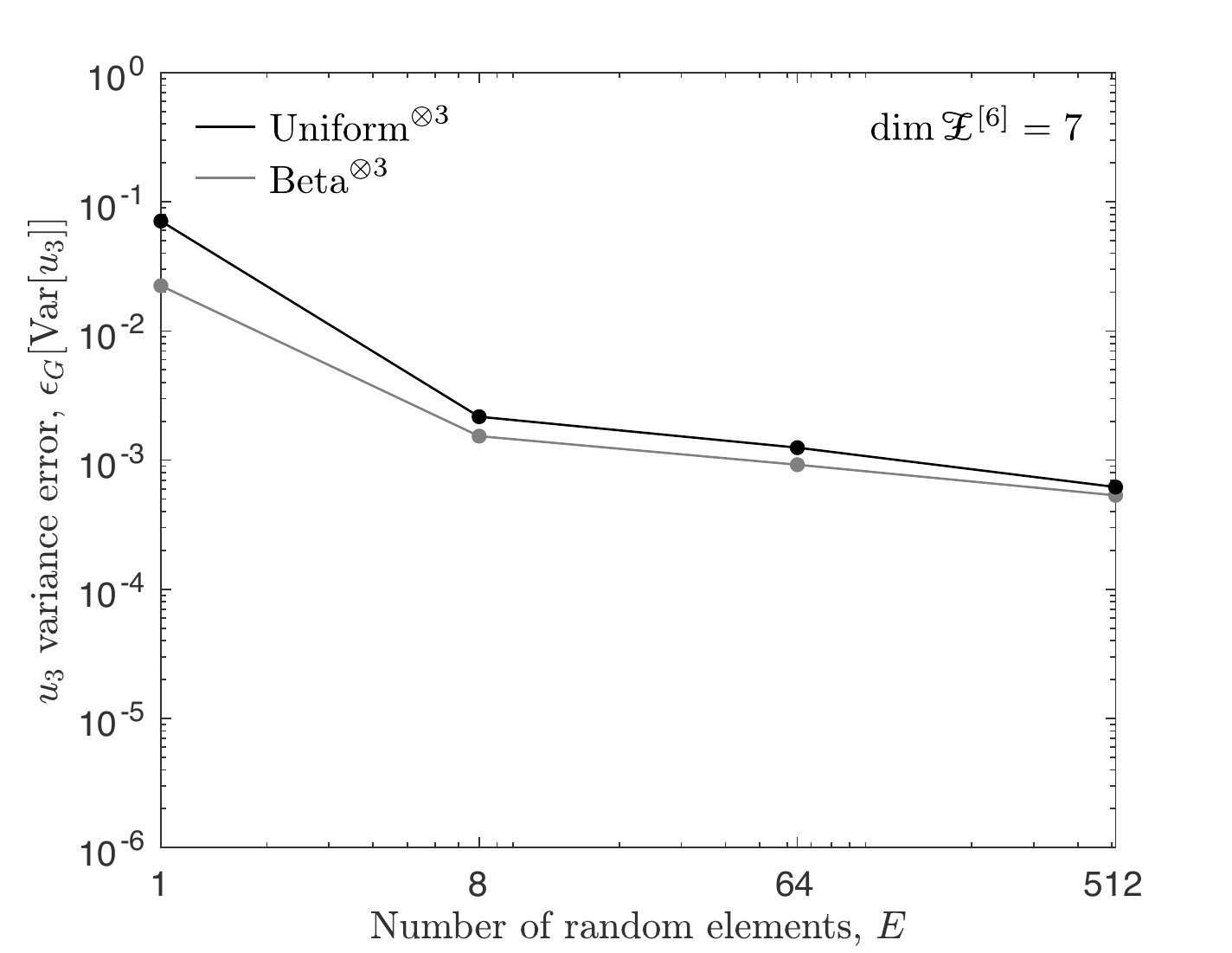}
\caption{Variance error for $u_3$}
\label{fig3Problem4_MEFSC_7BV_u3_Var_GlobalError}
\end{subfigure}
\caption{\emph{Problem 4} --- Global error of $\mathbf{E}[u_1]$, $\mathrm{Var}[u_1]$, $\mathbf{E}[u_3]$ and $\mathrm{Var}[u_3]$ for different $(h,p)$-discretization levels of RFS}
\label{fig3Problem4_MEFSC_7BV_GlobalError}
\end{figure}

Moreover, Figs.~\ref{fig3Problem4_Uniform3_MEFSC_7BV_LocalError} to \ref{fig3Problem4_MEFSC_7BV_GlobalError} present the local and global errors in mean and variance of $u_1$ and $u_3$ for the two distributions chosen for $\xi$.
From these plots, it is observed that increasing the number of elements helps improve the accuracy of the ME-FSC results.
However, because the Monte Carlo solution is not exact, the errors tend to stagnate around $10^{-3}$ and $10^{-4}$ when the ME-FSC results are compared to Monte Carlo.
This is the reason why, in the case of the mean, the accuracy of the ME-FSC results does not improve as the number of elements increases from 8 to 512; but, in the case of the variance, they do improve because the errors obtained with ME-FSC are above $10^{-4}$.
Therefore, to achieve comparable solution accuracy to Monte Carlo, around 512 elements are needed for ME-FSC.

% Problem 4 [2a/3]:

% Problem 4 [2b/3]:

% Problem 4 [3/3]:

For this problem, the computational cost associated with a simulation with 512 elements and 7 basis vectors per element was 82 seconds, whereas the computational cost associated with a Monte Carlo simulation with one million realizations was 97 seconds.
This outcome indicates that simulations conducted with ME-FSC were at least 15\% faster than those conducted with Monte Carlo.
For higher-dimensional probability spaces, it would be necessary to implement a different quadrature rule that is not sensitive to the curse-of-dimensionality issue in order to speed up the computation of the inner products and reduce the computational cost associated with ME-FSC.

\section{Conclusion}\label{sec3con}

In this paper we have presented an extension of the FSC method \cite{esquivel2020flow,esquivel2021flow}, called the \emph{multi-element flow-driven spectral chaos} (ME-FSC) method, in order to deal with discontinuities and long-time integration of stochastic dynamical systems more efficiently.
In ME-FSC, the stochastic problem is divided into multiple sub-problems so that every sub-problem can be solved individually using the FSC method.
In a subsequent step, the results are aggregated and the probability moments of interest are computed using the law of total probability.
An underlying characteristic of the FSC method makes it possible for the probability information to be transferred exactly from one random function space to another.
Furthermore, the postulated ME-FSC scheme is capable of reducing the error propagation by several orders of magnitude.

Four representative problems were presented in this paper to show the effectiveness of the ME-FSC approach.
The first two problems dealt with systems governed by a linear stochastic ODE.
They were selected because an exact solution was available for each of these problems.
The third problem dealt with a system governed by a nonlinear stochastic ODE (the Van-der-Pol oscillator).
This problem was utilized to study the effectiveness of ME-FSC in the resolution of nonlinear problems.
The last problem dealt with a nonlinear system of stochastic ODEs (the Kraichnan-Orszag three-mode problem).
It was employed to test the performance of ME-FSC when stochastic discontinuities are present in the system.
Based on our findings, we conclude that the ME-FSC method is capable of producing accurate solutions as compared to the exact solution (when available). 
For problems with no closed-form solutions, the ME-FSC method is capable of reproducing a well-resolved Monte Carlo simulation with one million realizations at a fraction of the computational cost required to do so.

The ME-FSC method is particularly useful when dealing with large stochastic dynamical systems, since it allows the analyst to decompose the problem into several smaller sub-problems and solve them separately.
In the future, the ME-FSC method can be augmented with adaptivity criteria in order to identify those regions in the probability space where the stochasticity of the system is evolving faster.
Moreover, if each sub-problem were to be solved simultaneously on multiple CPUs, the computation time needed to solve each problem could be reduced considerably.

\section*{Acknowledgements}
The first author gratefully acknowledges the scholarship granted by Colciencias and Atl\'antico Department (Colombia) under Call 673 to pursue a Ph.D.~degree in Structural Engineering at Purdue University, West Lafayette, IN.
The authors also gratefully acknowledge the support of the National Science Foundation (CNS-1136075, DMS-1555072, DMS-1736364, DMS-2053746, and DMS-2134209), Brookhaven National Laboratory Subcontract 382247, and U.S.~Department of Energy (DOE) Office of Science Advanced Scientific Computing Research program DE-SC0021142).

\section*{Credit author statement}
{\bf H.E.}~conceived the mathematical models, implemented the methods, designed the numerical experiments, interpreted the results, and wrote the paper. {\bf A.P.}~and {\bf G.L.}~supported the study, edited, and reviewed the final manuscript. All the authors gave their final approval for publication.

%\section*{References}
\bibliography{references3}

\begin{thebibliography}{10}
\expandafter\ifx\csname url\endcsname\relax
  \def\url#1{\texttt{#1}}\fi
\expandafter\ifx\csname urlprefix\endcsname\relax\def\urlprefix{URL }\fi
\expandafter\ifx\csname href\endcsname\relax
  \def\href#1#2{#2} \def\path#1{#1}\fi

\bibitem{xiu2002wiener}
D.~Xiu, G.~E. Karniadakis, The {W}iener-{A}skey polynomial chaos for stochastic
  differential equations, SIAM journal on scientific computing 24~(2) (2002)
  619--644.

\bibitem{gerritsma2010time}
M.~Gerritsma, J.-B. Van~der Steen, P.~Vos, G.~Karniadakis, Time-dependent
  generalized polynomial chaos, Journal of Computational Physics 229~(22)
  (2010) 8333--8363.

\bibitem{heuveline2014hybrid}
V.~Heuveline, M.~Schick, A hybrid generalized polynomial chaos method for
  stochastic dynamical systems, International Journal for Uncertainty
  Quantification 4~(1) (2014).

\bibitem{luchtenburg2014long}
D.~M. Luchtenburg, S.~L. Brunton, C.~W. Rowley, Long-time uncertainty
  propagation using generalized polynomial chaos and flow map composition,
  Journal of Computational Physics 274 (2014) 783--802.

\bibitem{ozen2016dynamical}
H.~C. Ozen, G.~Bal, Dynamical polynomial chaos expansions and long time
  evolution of differential equations with random forcing, SIAM/ASA Journal on
  Uncertainty Quantification 4~(1) (2016) 609--635.

\bibitem{wan2005adaptive}
X.~Wan, G.~E. Karniadakis, An adaptive multi-element generalized polynomial
  chaos method for stochastic differential equations, Journal of Computational
  Physics 209~(2) (2005) 617--642.

\bibitem{wan2006long}
X.~Wan, G.~E. Karniadakis, Long-term behavior of polynomial chaos in stochastic
  flow simulations, Computer methods in applied mechanics and engineering
  195~(41-43) (2006) 5582--5596.

\bibitem{wan2009error}
X.~Wan, G.~E. Karniadakis, Error control in multi-element generalized
  polynomial chaos method for elliptic problems with random coefficients,
  Communications in Computational physics 5~(2-4) (2009) 793--820.

\bibitem{wan2006beyond}
X.~Wan, G.~E. Karniadakis, Beyond {W}iener-{A}skey expansions: Handling
  arbitrary {PDFs}, Journal of Scientific Computing 27~(1-3) (2006) 455--464.

\bibitem{foo2008multi}
J.~Foo, X.~Wan, G.~E. Karniadakis, The multi-element probabilistic collocation
  method ({ME-PCM}): Error analysis and applications, Journal of Computational
  Physics 227~(22) (2008) 9572--9595.

\bibitem{foo2010multi}
J.~Foo, G.~E. Karniadakis, Multi-element probabilistic collocation method in
  high dimensions, Journal of Computational Physics 229~(5) (2010) 1536--1557.

\bibitem{zheng2015adaptive}
M.~Zheng, X.~Wan, G.~E. Karniadakis, Adaptive multi-element polynomial chaos
  with discrete measure: Algorithms and application to spdes, Applied Numerical
  Mathematics 90 (2015) 91--110.

\bibitem{kawai2020multi}
S.~Kawai, A.~Oyama, Multi-element stochastic galerkin method based on edge
  detection for uncertainty quantification of discontinuous responses, Journal
  of Verification, Validation and Uncertainty Quantification 5~(4) (2020)
  041004.

\bibitem{asokan2006stochastic}
B.~V. Asokan, N.~Zabaras, A stochastic variational multiscale method for
  diffusion in heterogeneous random media, Journal of Computational Physics
  218~(2) (2006) 654--676.

\bibitem{witteveen2009adaptive}
J.~A. Witteveen, A.~Loeven, H.~Bijl, An adaptive stochastic finite elements
  approach based on newton--cotes quadrature in simplex elements, Computers \&
  Fluids 38~(6) (2009) 1270--1288.

\bibitem{nouy2010extended}
A.~Nouy, A.~Clement, extended stochastic finite element method for the
  numerical simulation of heterogeneous materials with random material
  interfaces, International Journal for Numerical Methods in Engineering
  83~(10) (2010) 1312--1344.

\bibitem{ganapathysubramanian2007sparse}
B.~Ganapathysubramanian, N.~Zabaras, Sparse grid collocation schemes for
  stochastic natural convection problems, Journal of Computational Physics
  225~(1) (2007) 652--685.

\bibitem{ma2009adaptive}
X.~Ma, N.~Zabaras, An adaptive hierarchical sparse grid collocation algorithm
  for the solution of stochastic differential equations, Journal of
  Computational Physics 228~(8) (2009) 3084--3113.

\bibitem{bhaduri2018efficient}
A.~Bhaduri, L.~Graham-Brady, An efficient adaptive sparse grid collocation
  method through derivative estimation, Probabilistic Engineering Mechanics 51
  (2018) 11--22.

\bibitem{bhaduri2018stochastic}
A.~Bhaduri, Y.~He, M.~D. Shields, L.~Graham-Brady, R.~M. Kirby, Stochastic
  collocation approach with adaptive mesh refinement for parametric uncertainty
  analysis, Journal of Computational Physics 371 (2018) 732--750.

\bibitem{agarwal2009domain}
N.~Agarwal, N.~R. Aluru, A domain adaptive stochastic collocation approach for
  analysis of {MEMS} under uncertainties, Journal of Computational Physics
  228~(20) (2009) 7662--7688.

\bibitem{marzouk2007stochastic}
Y.~M. Marzouk, H.~N. Najm, L.~A. Rahn, Stochastic spectral methods for
  efficient bayesian solution of inverse problems, Journal of Computational
  Physics 224~(2) (2007) 560--586.

\bibitem{mohan2008multi}
P.~S. Mohan, P.~B. Nair, A.~J. Keane, Multi-element stochastic reduced basis
  methods, Computer Methods in Applied Mechanics and Engineering 197~(17-18)
  (2008) 1495--1506.

\bibitem{sarrouy2013piecewise}
E.~Sarrouy, O.~Dessombz, J.-J. Sinou, Piecewise polynomial chaos expansion with
  an application to brake squeal of a linear brake system, Journal of Sound and
  Vibration 332~(3) (2013) 577--594.

\bibitem{jakeman2013minimal}
J.~D. Jakeman, A.~Narayan, D.~Xiu, Minimal multi-element stochastic collocation
  for uncertainty quantification of discontinuous functions, Journal of
  Computational Physics 242 (2013) 790--808.

\bibitem{le2004uncertainty}
O.~Le~Ma{\i}tre, O.~Knio, H.~Najm, R.~Ghanem, Uncertainty propagation using
  wiener--haar expansions, Journal of computational Physics 197~(1) (2004)
  28--57.

\bibitem{le2004multi}
O.~Le~Ma{\i}tre, H.~N. Najm, R.~Ghanem, O.~Knio, Multi-resolution analysis of
  wiener-type uncertainty propagation schemes, Journal of Computational Physics
  197~(2) (2004) 502--531.

\bibitem{gittelson2014adaptive}
C.~J. Gittelson, Adaptive wavelet methods for elliptic partial differential
  equations with random operators, Numerische Mathematik 126~(3) (2014)
  471--513.

\bibitem{gittelson2013adaptive}
C.~Gittelson, An adaptive stochastic galerkin method for random elliptic
  operators, Mathematics of Computation 82~(283) (2013) 1515--1541.

\bibitem{cho2013adaptive}
H.~Cho, D.~Venturi, G.~E. Karniadakis, Adaptive discontinuous galerkin method
  for response-excitation pdf equations, SIAM Journal on Scientific Computing
  35~(4) (2013) B890--B911.

\bibitem{chen2015local}
Y.~Chen, J.~Jakeman, C.~Gittelson, D.~Xiu, Local polynomial chaos expansion for
  linear differential equations with high dimensional random inputs, SIAM
  Journal on Scientific Computing 37~(1) (2015) A79--A102.

\bibitem{li2016unified}
J.~Li, P.~Stinis, A unified framework for mesh refinement in random and
  physical space, Journal of Computational Physics 323 (2016) 243--264.

\bibitem{abdedou2019non}
A.~Abdedou, A.~Soulaimani, A non-intrusive b-splines b{\'e}zier elements-based
  method for uncertainty propagation, Computer Methods in Applied Mechanics and
  Engineering 345 (2019) 774--804.

\bibitem{eckert2020polynomial}
C.~Eckert, M.~Beer, P.~D. Spanos, A polynomial chaos method for arbitrary
  random inputs using b-splines, Probabilistic Engineering Mechanics 60 (2020)
  103051.

\bibitem{esquivel2020flow}
H.~Esquivel, A.~Prakash, G.~Lin, Flow-driven spectral chaos ({FSC}) method for
  simulating long-time dynamics of arbitrary-order non-linear stochastic
  dynamical systems, Journal of Computational Physics (2020) 110044.

\bibitem{esquivel2021flow}
H.~Esquivel, A.~Prakash, G.~Lin, Flow-driven spectral chaos (fsc) method for
  long-time integration of second-order stochastic dynamical systems, Journal
  of Computational and Applied Mathematics (2021) 113674.

\bibitem{cheney2010linear}
W.~Cheney, D.~Kincaid, Linear algebra: Theory and applications, 2nd Edition,
  Jones \& Bartlett Learning, 2010.

\bibitem{apple2021frameworks}
{A}pple {I}nc., \href{https://developer.apple.com/documentation}{{A}pple's
  frameworks through the {X}code {IDE}} (2021).
\newline\urlprefix\url{https://developer.apple.com/documentation}

\bibitem{swift2021language}
{A}pple {I}nc.~et al., \href{https://swift.org}{{S}wift programming language}
  (2021).
\newline\urlprefix\url{https://swift.org}

\bibitem{orszag1967dynamical}
S.~A. Orszag, L.~Bissonnette, Dynamical properties of truncated
  {W}iener-{H}ermite expansions, The Physics of Fluids 10~(12) (1967)
  2603--2613.

\end{thebibliography}

\end{document}